\documentclass{article}
\usepackage{CJK,CJKnumb,CJKulem,times,dsfont,ifthen,mathrsfs,latexsym,amsfonts,color}
\usepackage{amsmath,amsthm,makeidx,fontenc,amssymb,bm,graphicx,psfrag,listings,curves,extarrows}
\makeindex
\newtheorem{definition}{Definition}[section]
\newtheorem{theorem}{Theorem}[section]
\newtheorem{lemma}{Lemma}[section]
\newtheorem{corollary}{Corollary}[section]
\newtheorem{proposition}{Proposition}[section]
\newtheorem{remark}{Remark}[section]

\newcommand{\s}{\section}

\newcommand{\R}{\mathbb R}

\newcommand{\lab}{\label}
\newcommand{\bt}{\begin{theorem}}
\newcommand{\et}{\end{theorem}}
\newcommand{\bl}{\begin{lemma}}
\newcommand{\el}{\end{lemma}}
\newcommand{\bd}{\begin{definition}}
\newcommand{\ed}{\end{definition}}
\newcommand{\bc}{\begin{corollary}}
\newcommand{\ec}{\end{corollary}}
\newcommand{\bp}{\begin{proof}}
\newcommand{\ep}{\end{proof}}
\newcommand{\bx}{\begin{example}}
\newcommand{\ex}{\end{example}}
\newcommand{\bi}{\begin{exercise}}
\newcommand{\ei}{\end{exercise}}
\newcommand{\bo}{\begin{proposition}}
\newcommand{\eo}{\end{proposition}}
\newcommand{\br}{\begin{remark}}
\newcommand{\er}{\end{remark}}
\newcommand{\be}{\begin{equation}}
\newcommand{\ee}{\end{equation}}
\newcommand{\ba}{\begin{align}}
\newcommand{\ea}{\end{align}}
\newcommand{\bn}{\begin{enumerate}}
\newcommand{\en}{\end{enumerate}}
\newcommand{\bg}{\begin{align*}}
\newcommand{\bcs}{\begin{cases}}
\newcommand{\ecs}{\end{cases}}

\newcommand{\NN}{{\mathbb N}}

\newcommand{\bean}{\begin{eqnarray*}}
\newcommand{\eean}{\end{eqnarray*}}

\numberwithin{equation}{section}

\begin{document}
\begin{CJK*}{GBK}{song}
\title{\bf  On Elliptic Equations and Systems involving critical Hardy-Sobolev exponents (non-limit case) }
\date{}
\author{
{\bf  X. Zhong   \&   W. Zou}\thanks{Supported by  NSFC(11025106, 11371212, 11271386) and the Both-Side Tsinghua Fund. } \\
\footnotesize Department of Mathematical Sciences, Tsinghua University,\\
\footnotesize Beijing 100084, China\\
\footnotesize   wzou@math.tsinghua.edu.cn}

\maketitle

\vskip0.6in

\begin{center}
\begin{minipage}{120mm}
\begin{center}{\bf Abstract}\end{center}
Let   $\Omega\subset \R^N$ ($N\geq 3$)  be  an open domain (may be unbounded) with $0\in \partial\Omega$ and $\partial\Omega$  be of  $C^2$ at $0$  with  the  negative mean curvature  $H(0)$.    By using variational methods, we consider the following elliptic systems involving multiple Hardy-Sobolev critical exponents,
$$\begin{cases}
-\Delta u+\lambda^*\frac{u}{|x|^{\sigma_0}}-\lambda_1 \frac{|u|^{2^*(s_1)-2}u}{|x|^{s_1}}=\lambda \frac{1}{|x|^{s_2}}|u|^{\alpha-2}u|v|^\beta\quad &\hbox{in}\;\Omega,\\
-\Delta v+\mu^*\frac{v}{|x|^{\eta_0}}-\mu_1 \frac{|v|^{2^*(s_1)-2}v}{|x|^{s_1}}=\mu \frac{1}{|x|^{s_2}}|u|^{\alpha}|v|^{\beta-2}v\quad &\hbox{in}\;\Omega,\\
(u,v)\in D_{0}^{1,2}(\Omega)\times D_{0}^{1,2}(\Omega),
\end{cases}$$
where  $ 0\leq \sigma_0, \eta_0, s_2<2,  s_1\in (0,2);$   the parameters  $  \lambda^*\neq 0, \mu^*\neq 0, \lambda_1>0, \mu_1>0, \lambda\mu>0$;   $\alpha,\beta>1$ satisfying $\alpha+\beta \leq 2^*(s_2)$.  Here,  $2^*(s):=\frac{2(N-s)}{N-2}$ is the critical Hardy-Sobolev exponent.   We  obtain the existence and nonexistence of ground state solution under different specific assumptions.  As  the by-product, we  study
\be\lab{zou=a1}
\begin{cases}
&\Delta u+\lambda \frac{u^p}{|x|^{s_1}}+\frac{u^{2^*(s_2)-1}}{|x|^{s_2}}=0\;\quad \hbox{in}\;\Omega,\\
&u(x)>0\;\hbox{in}\;\Omega,\\
& u(x)=0\;\hbox{on}\;\partial\Omega,
\end{cases}
\ee
we  also   obtain the existence and nonexistence of   solution under different  hypotheses.  In particular,
we give   a  partial  answers to  a generalized open problem proposed by  Y. Y. Li  and  C. S. Lin  (ARMA, 2012).
Around the above two types of  equation or systems, we  systematically   study  the  elliptic equations  which
have multiple singular terms and are defined on any open domain.  We establish some fundamental results.
\vskip0.23in

{\it   Key  words:}  Elliptic system, Ground state, Hardy-Sobolev exponent.

\end{minipage}
\end{center}
\vskip0.26in

%
%
%
%
%

\newpage
\s{Introduction}
\renewcommand{\theequation}{1.\arabic{equation}}
\renewcommand{\theremark}{1.\arabic{remark}}
\renewcommand{\thedefinition}{1.\arabic{definition}}
Let   $\Omega\subset \R^N$ ($N\geq 3$)  be  an open domain with $0\in \partial\Omega$ and $\partial\Omega$  be of  $C^2$ at $0$  with  the  negative mean curvature  $H(0)$.  We study the following nonlinear elliptic systems
\be\lab{P}
\begin{cases}
-\Delta u+\lambda^*\frac{u}{|x|^{\sigma_0}}-\lambda_1 \frac{|u|^{2^*(s_1)-2}u}{|x|^{s_1}}=\lambda \frac{1}{|x|^{s_2}}|u|^{\alpha-2}u|v|^\beta\quad &\hbox{in}\;\Omega,\\
-\Delta v+\mu^*\frac{v}{|x|^{\eta_0}}-\mu_1 \frac{|v|^{2^*(s_1)-2}v}{|x|^{s_1}}=\mu \frac{1}{|x|^{s_2}}|u|^{\alpha}|v|^{\beta-2}v\quad &\hbox{in}\;\Omega,\\
(u,v)\in  H_0^1(\Omega)\times H_0^1(\Omega),
\end{cases}
\ee
Note that the parameters  $\lambda^*, \mu^*, \lambda_1, \mu_1, \lambda, \mu$ may different from each other.  Hence, the system has no variational structure.

The interest in  nonlinear Schr\"odinger systems is motivated by applications
to nonlinear optics,plasma physics, condensed matter physics, etc.
For example, the coupled nonlinear Schr\"odinger systems arise in the description
of several physical phenomena such as the propagation of pulses in birefringent
optical fibers and Kerr-like photorefractive media, see \cite{AkhmedievAnkiewicz.1999,EvangelidesMollenauerGordonBergano.1992,Kaminow.1981,Menyuk.1987,Menyuk.1989,WaiMenyukChen.1991}, etc.
The problem comes from the physical phenomenon with a clear practical significance. Research on solutions under different situations, not only correspond to different
Physical interpretation, but also has a pure mathematical theoretical significance, with a high scientific value.
Hence, the coupled nonlinear Schr\"odinger systems are widely
studied in recently years, we refer the readers to \cite{AbdellaouiFelliPeral.2009,AmbrosettiColorado.2006,LinWei.2005,MaiaMontefuscoPellacci.2006,Sirakov.2007,ZhongZou.2015, ZhongZou.2015arXiv:1503.08917v1[math.AP]31Mar2015}
and the
references therein.

For any $s\in [0,2]$, we define $\displaystyle\|u\|_{s,p}^{p}=\int_{\Omega}\frac{|u|^p}{|x|^s}dx.$ We also use the notations $\displaystyle\|u\|_p=\|u\|_{0,p}.$
The Hardy-Sobolev inequality \cite{CaffarelliKohnNirenberg.1984,CatrinaWang.2001,GhoussoubYuan.2000}
asserts that $\displaystyle D_{0}^{1,2}(\R^N)\hookrightarrow L^{2^*(s)}\left(\R^N,\frac{dx}{|x|^s}\right)$ is an embedding for $s\in [0,2]$. For a general open domain $\Omega$, there exists a positive constant $C(s,\Omega)$ such that
$$\int_\Omega |\nabla u|^2\geq C(s, \Omega)\left(\int_\Omega \frac{|u|^{2^*(s)}}{|x|^s}dx\right)^{\frac{2}{2^*(s)}},\quad u\in D_{0}^{1,2}(\Omega).$$
Define
$\mu_{s_1}(\Omega)$ is defined as
\be\lab{2014-2-27-e2}
\mu_{s_1}(\Omega)=\inf\left\{\frac{\int_\Omega |\nabla u|^2 dx}{\left(\int_\Omega \frac{|u|^{2^*(s_1)}}{|x|^{s_1}}dx\right)^{\frac{2}{2^*(s_1)}}};\;u\in D_{0}^{1,2}(\Omega)\backslash\{0\}\right\}.
\ee
Consider the case of $\Omega=\R_+^N$, it is well known that
the extremal function of $\mu_{s_1}(\R_+^N)$ is  parallel to the ground state solution of the following problem:
\be\lab{BP}
\begin{cases}
-\Delta u=\frac{|u|^{2^*(s_1)-2}u}{|x|^{s_1}}\quad &\hbox{in}\;\R_+^N,\\
u=0\;\hbox{on}\;\partial\R_+^N.
\end{cases}
\ee
We note that the existence of ground state solution of (\ref{BP}) for $0<s_1<2$ is solved by Ghoussoub and Robert \cite{GhoussoubRobert.2006}.
They also gave out some properties about the regularity, symmetry and decay estimates.
And the instanton $\displaystyle U(x):=C\left(\kappa+|x|^{2-s_2}\right)^{-\frac{N-2}{2-s_2}}$  for  $0\leq s_2<2$  is a ground state solution to (\ref{BPP}) below (see \cite{Lieb.1983} and \cite{Talenti.1976a}):
\be\lab{BPP}
\begin{cases}
\Delta u+\frac{u^{2^*(s_2)-1}}{|x|^{s_2}}=0\quad &\hbox{in}\;\R^N,\\
u>0\;\;\hbox{in}\;\R^N\;\;\;\hbox{and}& u\rightarrow 0\;\hbox{as}\;|x|\rightarrow +\infty.
\end{cases}
\ee
That
$0\in \partial\Omega$ and
$\partial\Omega$ is smooth at $0$ has become a hot research topic in recent years as the curvature at $0$ plays an important role, see \cite{ChernLin.2010,GhoussoubKang.2004,GhoussoubRobert.2006,HsiaLinWadade.2010} .etc.
When $0\in \partial\Omega$ with $H(0)<0$, $\mu_{s_1}(\Omega)<\mu_{s_1}(\R_+^N)$  was proved in \cite{GhoussoubKang.2004} for $N\geq 3$ and \cite{LinWadadeothers.2012} for $N\geq 4$. It is standard to apply the blow-up analysis to obtain that $\mu_{s_1}(\Omega)$ can be achieved by some positive $u\in H_0^1(\Omega)$ (see\cite[Corollary 3.2]{GhoussoubKang.2004}), which is a ground state solution of $$\begin{cases}
-\Delta u=\frac{|u|^{2^*(s_1)-2}u}{|x|^{s_1}}\quad &\hbox{in}\;\Omega,\\
u=0\;\hbox{on}\;\partial\Omega.
\end{cases}$$ and the ground state value equals $\displaystyle \left(\frac{1}{2}-\frac{1}{2^*(s_1)}\right)\mu_{s_1}(\Omega)^{\frac{N-s_1}{2-s_1}}$.

\vskip 0.2in
Although \eqref{P} has no variational structure, for $\kappa\in \R$ such that $\kappa \lambda>0$, we can take $t_1>0,t_2>0$ properly  such that
$$\lambda \left(\frac{1}{t_1}\right)^{\alpha-2}\left(\frac{1}{t_2}\right)^\beta=\kappa\cdot \alpha;  \quad \mu\left(\frac{1}{t_2}\right)^{\beta-2}\left(\frac{1}{t_1}\right)^\alpha=\kappa\cdot\beta.$$
Denote $\tilde{u}=t_1u, \tilde{v}=t_2v$, then
a direct calculation  show that $(u,v)$ is a solution to (\ref{P}) if and only if $(\tilde{u}, \tilde{v})$ is a solution to the following problem:
\be
\begin{cases}
-\Delta \tilde{u}+\lambda^*\frac{u}{|x|^{\sigma_0}}-\tilde{\lambda}_1 \frac{|\tilde{u}|^{2^*(s_1)-2}\tilde{u}}{|x|^{s_1}}=\kappa\alpha \frac{1}{|x|^{s_2}}|\tilde{u}|^{\alpha-2}\tilde{u}|\tilde{v}|^\beta\quad &\hbox{in}\;\Omega,\\
-\Delta \tilde{v}+\mu^*\frac{v}{|x|^{\eta_0}}-\tilde{\mu}_1 \frac{|\tilde{v}|^{2^*(s_1)-2}\tilde{v}}{|x|^{s_1}}=\kappa\beta \frac{1}{|x|^{s_2}}|\tilde{u}|^{\alpha}|\tilde{v}|^{\beta-2}\tilde{v}\quad &\hbox{in}\;\Omega,\\
(\tilde{u},\tilde{v})\in D_{0}^{1,2}(\Omega)\times D_{0}^{1,2}(\Omega),
\end{cases}
\ee
where $\tilde{\lambda}_1=\lambda_1 t_{1}^{2-2^*(s_1)}, \tilde{\mu}_1=\mu_1 t_{2}^{2-2^*(s_1)}$.
Hence, without loss of generality, we will study the following problem:
\be\lab{PT}
\begin{cases}
-\Delta u+\lambda^*\frac{u}{|x|^{\sigma_0}}-\lambda \frac{|u|^{2^*(s_1)-2}u}{|x|^{s_1}}=\kappa\alpha \frac{1}{|x|^{s_2}}|u|^{\alpha-2}u|v|^\beta\quad &\hbox{in}\;\Omega,\\
-\Delta v+\mu^*\frac{v}{|x|^{\eta_0}}-\mu \frac{|v|^{2^*(s_1)-2}v}{|x|^{s_1}}=\kappa\beta \frac{1}{|x|^{s_2}}|u|^{\alpha}|v|^{\beta-2}v\quad &\hbox{in}\;\Omega,\\
(u,v)\in D_{0}^{1,2}(\Omega)\times D_{0}^{1,2}(\Omega),
\end{cases}
\ee
where $\kappa$ is a parameter.

There seems to has no article involves the system case with  $0\in \partial\Omega$. Hence, in present paper, our aim is to establish related results for the elliptic systems on the domain with $0\in \partial\Omega$. There also seems to have no article studies the unbounded domain which is not a limit domain (i.e., $\displaystyle \Omega\neq \lim_{t\rightarrow 0}\frac{\Omega}{t}$, where $\displaystyle \frac{\Omega}{t}:=\left\{x\in \R^N\;:\;tx\in \Omega\right\}$). We note that  the limit equation of \eqref{PT} has been studied in \cite{ZhongZou.2015arXiv:1504.01005v1[math.AP]4Apr2015,ZhongZou.2015arXiv:submit/1229202[math.AP]12Apr2015}.

We are interested in the existence of nontrivial ground state solutions, hence we need a well study in a single equation, especially for the unbounded domain case.
We obtain the following results:
\bt\lab{2015-4-24-th1}
Assume that $\lambda^*>0, \Omega\subset \R^N$ is $C^1$ domain such that $|\Omega|<\infty, 0\in \partial\Omega$ and the boundary $\partial\Omega$ is $C^2$ at $0$ with the mean curvature $H(0)<0$. We also assume that $0\leq s_1<1$ or $1\leq s_1<2$ with $\lambda^*$ small enough. Then problem
\be\lab{2015-4-27-le1}
\begin{cases}
-\Delta u+\lambda^*\frac{u}{|x|^{s_1}}=\lambda \frac{|u|^{2^*(s_2)-2}u}{|x|^{s_2}}\quad&\hbox{in}\;\Omega\\
u(x)=0&\hbox{on}\;\partial\Omega,
\end{cases}
\ee
 possesses a positive ground state solution $u$.
\et

\bt\lab{2015-4-24-th2}
Let  $\Omega\subset \R^N$ be a  $C^1$ domain such that $|\Omega|<\infty, 0\in \partial\Omega$ and the boundary $\partial\Omega$ is $C^2$ at $0$ with the mean curvature $H(0)<0$. We also assume that $-\lambda_{1,s_1}(\Omega)<\lambda^*<0$ and $0\leq s_1<2$, then problem \eqref{2015-4-27-le1} possesses a positive ground state solution $u$.
\et

\bt\lab{2014-6-20-wth1}
Suppose that $\Omega\subset \R^N$ is a $C^1$ open domain with $0\in \partial \Omega$ and $\partial\Omega$ is $C^2$ at $0$, $H(0)<0$. We also assume that
$$
\begin{cases}
\lambda>0, 1\leq p< \min\{2^*(s_2)-1, 2^*(s_1)-1\}\quad &\hbox{if}\;|\Omega|=\infty;\\\\
either \lambda< 0, 1<p< 2^*(s_1)-1\;\hbox{ or }\;\\
\lambda> 0, 1<p< \min\{2^*(s_2)-1, 2^*(s_1)-1\}&\hbox{if}\;|\Omega|<\infty.
\end{cases}
$$
Moreover, if $\lambda>0$, we require that $p<\frac{N-2s_1}{N-2}$ or $p=\frac{N-2s_1}{N-2}$ with $|\lambda|$ small enough.
Then problem
\be\lab{2015-4-27-le2}
\begin{cases}
&-\Delta u+\lambda \frac{u^p}{|x|^{s_1}}=\frac{u^{2^*(s_2)-1}}{|x|^{s_2}}\;\quad \hbox{in}\;\Omega,\\
&u(x)>0\;\hbox{in}\;\Omega,\\
& u(x)=0\;\hbox{on}\;\partial\Omega,
\end{cases}
\ee
 possesses a positive ground state solution.
\et

As a product, basing on some lemmas established in present paper for unbounded domain, we can also give a partial answers to Li-Lin's open problem (see \cite[Remark 1.2]{LiLin.2012}) in the unbounded domain case, see Theorem \ref{2014-6-24-th1} and Theorem \ref{2014-6-24-th2} in section 3,  which are parallel to \cite[Theorem 1.5 and 1.6]{CeramiZhongZou.2015.DOI10.1007/s00526-015-0844-z.} which are studied in bounded domain.

\vskip 0.2in
We shall prove that when $|\Omega|<\infty$ and $0\leq \sigma<2$, $\displaystyle -\Delta u=\lambda \frac{u}{|x|^\sigma}, u\in D_{0}^{1,2}(\Omega)$ possesses a sequence of eigenvalues such that
$$0<\lambda_{1,\sigma}<\lambda_{2,\sigma}\leq \cdots\leq \lambda_{k,\sigma}\leq \lambda_{k+1,\sigma}\leq \cdots.$$
Define $\displaystyle E:=D_{0}^{1,2}(\Omega)\cap L^2\left(\Omega,   \frac{dx}{|x|^{\sigma_0}}\right)$ and $\displaystyle F:=D_{0}^{1,2}(\Omega)\cap L^2\left((\Omega,   \frac{dx}{|x|^{\eta_0}})\right)$. Then for some proper range of $\lambda^*$ and $\mu^*$, we can define the following norms:
$$\|u\|_E:=\left(\int_{\Omega}\left(|\nabla u|^2+\lambda^* \frac{|u|^2}{|x|^{\sigma_0}}\right)dx\right)^{\frac{1}{2}}, \|v\|_F:=\left(\int_{\Omega}\left(|\nabla v|^2+\mu^* \frac{|v|^2}{|x|^{\eta_0}}\right)dx\right)^{\frac{1}{2}}.$$
Our framework space for the system case is $\mathscr{D}:=E\times F$.

Collect the results about a single equation, we have that equation
\be\lab{2015-4-27-cle1}
\begin{cases}
-\Delta u+\lambda^*\frac{u}{|x|^{\sigma_0}}=\lambda \frac{|u|^{2^*(s_1)-2}u}{|x|^{s_1}}\quad&\hbox{in}\;\Omega\\
u(x)=0&\hbox{on}\;\partial\Omega,
\end{cases}
\ee
possesses a positive ground state solution when the following condition {\bf $(H_1)$} is satisfied:\\
{\bf $(H_1)$:}one of the following holds:
\begin{itemize}
\item[$(a.1)$]$|\Omega|<\infty, -\lambda_{1,\sigma_0}(\Omega)<\lambda^*<0, 0\leq \sigma_0<2$ (see Theorem \ref{2015-4-24-th2}).
\item[$(a.2)$]$|\Omega|<\infty, \lambda^*>0,0\leq \sigma_0<1$ (see Theorem \ref{2015-4-24-th1}).
\item[$(a.3)$]$|\Omega|<\infty,\lambda^*>0$ small enough, $1\leq \sigma_0<2$ (see Theorem \ref{2015-4-24-th1}).
\item[$(a.4)$]$|\Omega|=\infty,\lambda^*>0,0\leq \sigma_0<1$ (see Theorem \ref{2014-6-20-wth1}).
\item[$(a.5)$]$|\Omega|=\infty,\lambda^*>0$ small enough,$1\leq \sigma_0<2$ (see Theorem \ref{2014-6-20-wth1}).
\end{itemize}

 Let $\tilde{u}=\lambda^{\frac{1}{2^*(s_1)-2}}u$, then $\tilde{u}$ is a ground state solution of
\be\lab{2014-4-11-we9}
\begin{cases}
-\Delta \tilde{u}+\lambda^*\frac{\tilde{u}}{|x|^{\sigma_0}}=\frac{|\tilde{u}|^{2^*(s_1)-2}\tilde{u}}{|x|^{s_1}}\quad &\hbox{in}\;\Omega,\\
 \tilde{u}(x)=0\quad &\hbox{on}\;\partial\Omega.
\end{cases}
\ee

Similarly, equation
\be\lab{2015-4-25-xe1}
\begin{cases}
-\Delta v+\mu^*\frac{v}{|x|^{\eta_0}}=\mu\frac{|v|^{2^*(s_1)-2}v}{|x|^{s_1}}\;\hbox{in}\;\Omega\\
\mu^*\neq 0, v=0\;\hbox{on}\;\partial\Omega
\end{cases}
\ee
possesses a positive ground state solution when the following condition {\bf $(H_2)$} is satisfied:\\
{\bf $(H_2)$:}one of the following holds:
\begin{itemize}
\item[$(b.1)$]$|\Omega|<\infty, -\lambda_{1,\eta_0}(\Omega)<\mu^*<0, 0\leq \eta_0<2$ (see Theorem \ref{2015-4-24-th2}).
\item[$(b.2)$]$|\Omega|<\infty, \mu^*>0,0\leq \eta_0<1$ (see Theorem \ref{2015-4-24-th1}).
\item[$(b.3)$]$|\Omega|<\infty,\mu^*>0$ small enough, $1\leq \eta_0<2$ (see Theorem \ref{2015-4-24-th1}).
\item[$(b.4)$]$|\Omega|=\infty,\mu^*>0,0\leq \eta_0<1$ (see Theorem \ref{2014-6-20-wth1}).
\item[$(b.5)$]$|\Omega|=\infty,\mu^*>0$ small enough,$1\leq \eta_0<2$ (see Theorem \ref{2014-6-20-wth1}).
\end{itemize}

Define
\be\lab{2015-4-27-le3}
\tilde{\eta}_1:=\inf_{v\in F\backslash \{0\}}\frac{\|v\|_F^2}{\int_\Omega \frac{|U|^{\alpha}|v|^2}{|x|^{s_2}}dx},
\ee
where $0<U\in E$ is a positive ground state solution of \eqref{2015-4-27-cle1}.
Similarly define
\be\lab{2015-4-27-le4}
\tilde{\eta}_2:=\inf_{u\in E\backslash \{0\}}\frac{\|u\|_E^2}{\int_\Omega \frac{|V|^{\beta}|u|^2}{|x|^{s_2}}dx},
\ee
where $0<V\in F$ is a positive ground state solution of \eqref{2015-4-25-xe1}.

We denote
\begin{itemize}
\item[{\bf $(A_{\sigma_0}^{*})$}]  either $\sigma_0<s_2$ or $\begin{cases}0<s_2\leq \sigma_0,\\ 2+\frac{4}{\sigma_0}\frac{\sigma_0-s_2}{N-2}< \alpha+\beta\end{cases}$;\\
\item[{\bf $(A_{\eta_0}^{*})$}]  either $\eta_0<s_2$ or $\begin{cases} 0<s_2\leq \eta_0,\\ 2+\frac{4}{\eta_0}\frac{\eta_0-s_2}{N-2}< \alpha+\beta\end{cases}$.
\end{itemize}

Here comes our first main result for the system case with $\alpha+\beta<2^*(s_2)$.
\bt\lab{2015-4-26-xth1}({\bf The case of $\alpha+\beta<2^*(s_2)$})
Assume $(H_1),(H_2), \kappa>0,1<\alpha,1<\beta, \alpha+\beta<2^*(s_2)$, and if $|\Omega|=\infty$, we assume further that $(A_{\sigma_0}^{*})$ and $(A_{\eta_0}^{*})$ hold. If $\beta<2$ or $\beta=2$ with $\kappa>\tilde{\eta}_1$, and if $\alpha<2$ or $\alpha=2$ with $\kappa>\tilde{\eta}_2$, then problem
\be\lab{Problem-System}
\begin{cases}
-\Delta u+\lambda^*\frac{u}{|x|^{\sigma_0}}-\lambda \frac{|u|^{2^*(s_1)-2}u}{|x|^{s_1}}=\kappa\alpha \frac{1}{|x|^{s_2}}|u|^{\alpha-2}u|v|^\beta\quad &\hbox{in}\;\Omega,\\
-\Delta v+\mu^*\frac{v}{|x|^{\eta_0}}-\mu \frac{|v|^{2^*(s_1)-2}v}{|x|^{s_1}}=\kappa\beta \frac{1}{|x|^{s_2}}|u|^{\alpha}|v|^{\beta-2}v\quad &\hbox{in}\;\Omega,\\
(u,v)\in \mathscr{D},
\end{cases}
\ee
possesses a nontrivial ground state solution $(u,v)$.
\et

\br
When $|\Omega|=\infty$ and $\alpha+\beta<2^*(s_2)$, the assumptions of $(A_{\sigma_0}^{*})$ and $(A_{\eta_0}^{*})$ play an crucial role to guarantee that $\int_\Omega \frac{|u|^\alpha |v|^\beta}{|x|^{s_2}}$ is well defined for all $(u,v)\in \mathscr{D}$ (see Proposition \ref{2014-6-5-xl1}).
\er

To study the critical couple case, for the technical reasons, we need the following assumptions $(\tilde{H}_1)$ and $(\tilde{H}_2)$, which are stronger than $(H_1)$ and $(H_2)$:\\
{\bf $(\tilde{H}_1)$:}one of the following holds:
\begin{itemize}
\item[$(\tilde{a}.1)$] $|\Omega|<\infty, -\lambda_{1,\sigma_0}(\Omega)<\lambda^*<0, 0\leq \sigma_0<2$.
\item[$(\tilde{a}.2)$] $|\Omega|\leq \infty, \lambda^*>0, 0\leq \sigma_0<1$.
\item[$(\tilde{a}.3)$] $|\Omega|=\infty,\lambda^*>0$ small enough and $\sigma_0=1$.
\end{itemize}
{\bf $(\tilde{H}_2)$:}one of the following holds:
\begin{itemize}
\item[$(\tilde{b}.1)$] $|\Omega|<\infty, -\lambda_{1,\eta_0}(\Omega)<\mu^*<0, 0\leq \eta_0<2$.
\item[$(\tilde{b}.2)$] $|\Omega|\leq \infty, \mu^*>0, 0\leq \eta_0<1$.
\item[$(\tilde{b}.3)$] $|\Omega|=\infty,\mu^*>0$ small enough and $\eta_0=1$.
\end{itemize}
\br
When $(\tilde{H}_1)$ and $(\tilde{H}_2)$ hold, it is easy to see that $(H_1)$ and $(H_2)$ are satisfied.
\er

\bt\lab{2015-4-26-xth2}({\bf The case of $\alpha+\beta=2^*(s_2)$})
Assume $(\tilde{H}_1),(\tilde{H}_2), \kappa>0, 1<\alpha,1<\beta, \alpha+\beta=2^*(s_2)$.
Especially, when $s_1=s_2=s\in (0,2)$, we assume further that one of the following holds:
\begin{itemize}
\item[(a.1)]$\lambda>\mu, 1<\beta<2$ or $\begin{cases}\beta=2\\ \kappa>\max\{\lambda,\tilde{\eta}_1\} \end{cases}$;
\item[(a.2)] $\lambda=\mu, \min\{\alpha,\beta\}<2$ or $\begin{cases}\min\{\alpha,\beta\}=2,\kappa>\max\{\lambda,\tilde{\eta}_1\} \end{cases}$;
\item[(a.3)] $\lambda<\mu, 1<\alpha<2$ or $\begin{cases}\alpha=2\\ \kappa>\max\{\mu,\tilde{\eta}_2\} \end{cases}$.
\end{itemize}
When $s_1,s_2\in (0,2)$ but $s_1\neq s_2$, we assume further that one of the following holds:
\begin{itemize}
\item[(b.1)]$\lambda>\mu, 1<\beta<2$ or $\begin{cases}\beta=2\\ \kappa>\max\{\eta_1,\tilde{\eta}_1\} \end{cases}$;
\item[(b.2)] $\lambda=\mu, \min\{\alpha,\beta\}<2$ or $\begin{cases}\min\{\alpha,\beta\}=2,\kappa>\max\{\eta_1,\eta_2,\tilde{\eta}_1,\tilde{\eta}_2\} \end{cases}$;
\item[(b.3)] $\lambda<\mu, 1<\alpha<2$ or $\begin{cases}\alpha=2\\ \kappa>\max\{\eta_2,\tilde{\eta}_2\} \end{cases}$.
\end{itemize}
Where the constant $\eta_1$ is that one given by \cite[Lemma 6.5]{ZhongZou.2015arXiv:1504.01005v1[math.AP]4Apr2015} and $\eta_2$ is the constant given by \cite[Corollary 6.4]{ZhongZou.2015arXiv:1504.01005v1[math.AP]4Apr2015}.

Then problem \eqref{Problem-System} possesses a nontrivial ground state solution $(u,v)$.
\et

\s{Preliminaries}
\renewcommand{\theequation}{2.\arabic{equation}}
\renewcommand{\theremark}{2.\arabic{remark}}
\renewcommand{\thedefinition}{2.\arabic{definition}}
\subsection{Eigenvalues and compact embedding}
Assume $ 0\leq \sigma<2, |\Omega|<\infty$, then $\frac{1}{|x|^\sigma}\in L^{\frac{N}{2}}(\Omega)$. By \cite[Lemma 2.13]{Willem.1996}, the functional $$\chi:D_{0}^{1,2}(\Omega)\rightarrow \R\;   \hbox{ with }\;  \chi(u)= \int_\Omega \frac{u^2}{|x|^{\sigma}}dx$$ is weakly continuous.
Define $\langle u, v\rangle_\sigma :=\int_{\Omega}\frac{uv}{|x|^\sigma}dx$, then it is easy to check that $\langle \cdot, \cdot\rangle_\sigma$ is an inner product. We say that $u$ and $v$ are  orthogonal if and only if $\langle u, v\rangle_\sigma=0$.

\bo\lab{2014-6-4-wl1}
Assume that $\Omega\subset \R^N (N\geq 3), 0\leq \sigma<2$ and $|\Omega|<\infty$, then the following eigenvalue problem
\be\lab{2014-6-4-we1}
-\Delta u=\lambda \frac{u}{|x|^\sigma}, \;u\in D_{0}^{1,2}(\Omega).
\ee
possesses a sequence of eigenvalues such that
$$0<\lambda_{1,\sigma}<\lambda_{2,\sigma}\leq \cdots\leq \lambda_{k,\sigma}\leq \lambda_{k+1,\sigma}\leq \cdots,$$
where each eigenvalue is repeated according to its multiplicity.
Moreover, $\lambda_{k,\sigma}\rightarrow \infty$ as $k\rightarrow \infty$. Let $e_{1,\sigma}, e_{2,\sigma}, e_{3,\sigma}, \cdots$ be the corresponding orthonormal eigenfunctions, then $e_1$ is positive and others are sign-changing.
\eo
\bp
It is well known for the case of $\sigma=0$. When $0<\sigma<2$, we refer to \cite[Corollary 4.2]{ZhongZou.2015arXiv:1504.00433v1[math.AP]2Apr2015}.
\ep

Define $E:=D_{0}^{1,2}(\Omega)\cap L^2(\Omega,   \frac{dx}{|x|^{\sigma_0}})$, if $|\Omega|<\infty$, we assume that $\lambda^*>-\lambda_{1,\sigma_0}$, then it is easy to see that
\be\lab{normE}
\|u\|_E:=\Big(\int_{\Omega}\big(|\nabla u|^2+\lambda^* \frac{|u|^2}{|x|^{\sigma_0}}\big)dx\Big)^{\frac{1}{2}}
\ee
is a norm, and it is equivalent to the norm of $\|u\|:=\big(\int_\Omega |\nabla u|^2dx\big)^{\frac{1}{2}}$ due to Proposition\ref{2014-6-4-wl1}. Similarly, we define $F:=D_{0}^{1,2}(\Omega)\cap L^2(\Omega,   \frac{dx}{|x|^{\eta_0}})$ and if $|\Omega|<\infty$, we assume that $\mu^*>-\lambda_{1,\eta_0}$, then
\be\lab{normF}
\|v\|_F:=\Big(\int_{\Omega}\big(|\nabla v|^2+\mu^* \frac{|v|^2}{|x|^{\eta_0}}\big)dx\Big)^{\frac{1}{2}}
\ee
is a norm which is also equivalent to the norm of $\|v\|$.
In present paper, if $|\Omega|=\infty$, we always assume that $\lambda^*>0, \mu^*>0$.


 \bd\lab{2013-8-4-def1}  (see \cite{Lions.1985a, Lions.1985} and  also\cite{ZhongZou.2014a})
Assume $\{\rho_k\}$ is a bounded sequence in $L^1(\R^N)$ and $\rho_k\geq 0$ satisfies
\be\lab{2013-8-4-e1}
\|\rho_k\|_{L^1}=\lambda+o(1), \;\lambda\geq 0.
\ee
Then we call this sequence $\{\rho_k\}$ is a tight sequence if  $ \forall\;\varepsilon>0,  \exists\;R>0$ such that
\be\lab{2013-8-4-e2}
\int_{|x|\geq R}\rho_k(x)dx<\varepsilon,\;\forall\;k\geq 1.
\ee
We call $u_k$ is a $L^p$ tight sequence, if $|u_k|^p$ is a tight sequence.
\ed

\vskip 0.2in

\br\lab{2014-6-9-r1}
It is well known that the embedding $H^1(\R^N)\hookrightarrow L^{p}(\R^N)$ is not compact for $2<p<2^*$. However,   we have the following compact embedding  result  for proper $p$.
\er

\bo\lab{2014-6-5-xl1}
Assume   $ \sigma,  s\in [0, 2)$.  Furthermore,  suppose  either

\begin{itemize}
\item[(i)]   $ 2<p<2^*(s)$ if $\sigma\leq s$  or
\item[(ii)]  $ 2+\frac{4}{\sigma}\frac{\sigma-s}{N-2}\leq p<2^*(s)$ if $\sigma>s$.
\end{itemize}
Set
$\displaystyle E:=D_{0}^{1,2}(\R^N)\cap L^2(\R^N, \frac{dx}{|x|^\sigma})$  with the norm
$$\|u\|_E:=\Big(\int_{\R^N} \big(|\nabla u|^2 +\frac{|u|^2}{|x|^\sigma}\big)dx\Big)^{\frac{1}{2}}.$$
Then \be\label{zou8=20} (E, \|\cdot\|_E) \hookrightarrow L^{p}\left(\R^N,   \frac{dx}{|x|^s} \right)\ee is a  continuous embedding. Moreover, when $\sigma<s$ with $ 2<p<2^*(s)$ or $0<s\leq \sigma$ with $2+\frac{4}{\sigma}\frac{\sigma-s}{N-2}< p<2^*(s)$, the embedding
(\ref{zou8=20})  is compact.
\eo
\bp
It can be proved by a modification of \cite[Lemma 3.1 and Lemma 3.2]{ZhongZou.2015arXiv:1504.00433v1[math.AP]2Apr2015}. However, for the convenience of the readers, we prefer to give the details.

Define
\be\lab{2014-6-5-xe1}
\begin{cases}
p_1:=\frac{2N+2p-Np-2s}{2-\sigma},\\
p_2:=\frac{Np+2s-p\sigma-2N}{2-\sigma},\\
\sigma_1:=\frac{(2N+2p-Np-2s)\sigma}{4-2\sigma},\\
\sigma_2:=\frac{4s+Np\sigma-2p\sigma-2N\sigma}{4-2\sigma}.
\end{cases}
\ee
Then we have the following results:
\be\lab{2014-6-5-xe2}
\begin{cases}
0<p_1<2,0<p_2, p_1+p_2=p,\\
\sigma_1=\frac{p_1}{2}\sigma, \sigma_1+\sigma_2=s,\\
\bar{\sigma}:=\frac{2\sigma_2}{2-p_1}\in [0,2),\\
\frac{2p_2}{2-p_1}=2^*(\bar{\sigma}).
\end{cases}
\ee
Then, by  the H\"{o}lder inequality, we have
\begin{align}\lab{2014-6-5-xe3}
\int_{\R^N}\frac{|u|^p}{|x|^s}dx=&\int_{\R^N} \frac{|u|^{p_1}}{|x|^{\sigma_1}} \frac{|u|^{p_2}}{|x|^{\sigma_2}}dx\nonumber\\
\leq&\Big(\int_{\R^N}\frac{|u|^2}{|x|^{\sigma}}dx\Big)^{\frac{p_1}{2}} \Big(\int_{\R^N}\frac{|u|^{2^*(\bar{\sigma})}}{|x|^{\bar{\sigma}}}dx\Big)^{\frac{2-p_1}{2}}.
\end{align}
By the Hardy-Sobolev inequality, there exists $C>0$ independent of $u$ such that
\begin{align}\lab{2014-6-5-xe4}
|u|_{s,  p}\leq C|u|_{\sigma,2}^{\frac{p_1}{p}} \|u\|^{\frac{p_2}{p}},
\end{align}
here $\displaystyle\|u\|=\big(\int_{\R^N}|\nabla u|^2dx\big)^{\frac{1}{2}}$ and the best constant $C$ satisfies
$$C\leq \mu_{\bar{\sigma}}(\R^N)^{-\frac{p_2}{2p}}.$$
When $s=\sigma=0$, it is just the well known
Gagliardo-Nirenberg inequality \cite{Gagliardo.1958, Nirenberg.1959}.
Next, we assume that $\{u_n\}\subset L^2(\R^N,    \frac{dx}{|x|^\sigma})\cap D_{0}^{1,2}(\R^N)$ is  a bounded sequence. Without loss of generality, we may assume that $u_n\rightharpoonup u$ in $L^2(\R^N,\frac{dx}{|x|^\sigma})\cap D_{0}^{1,2}(\R^N)$ and $u_n\rightarrow u$ a.e. in $\R^N$.

\vskip0.1in

\noindent{\bf Case 1:  $\sigma<s, 2<p<2^*(s)$.} Fix $\sigma<\bar{s}<s$, then we see that $2<p<2^*(\bar{s})$, by the embedding result, we obtain that there exists some $C>0$ such that
\be\lab{2014-6-5-we1}
|u_n|_{\bar{s},p}\leq C\;\hbox{for all}\;n.
\ee
Then we can see that
\be\lab{2014-6-5-we2}
\int_{|x|>R} \frac{|u_n|^p}{|x|^s}dx\leq \frac{1}{R^{s-\bar{s}}} \int_{|x|>R} \frac{|u_n|^p}{|x|^{\bar{s}}}dx\rightarrow 0\;\hbox{as}\;R\rightarrow \infty.
\ee
Hence $\frac{|u|^p}{|x|^s}$ is a tight sequence. On the other hand, since $p<2^*(s)$, for any fixed $R>0$, we have that
\be\lab{2014-6-5-we3}
\int_{|x|\leq R} \frac{|u_n-u|^p}{|x|^s}dx\rightarrow 0\;\hbox{as}\;n\rightarrow \infty.
\ee
It follows  from  (\ref{2014-6-5-we2}) and (\ref{2014-6-5-we3}) that
\be\lab{2014-6-5-we4}
\int_{\R^N}\frac{|u_n-u|^p}{|x|^s}dx\rightarrow 0\;\hbox{as}\;n\rightarrow \infty.
\ee
Hence,  the embedding $L^2(\R^N,  \frac{dx}{|x|^\sigma})\cap D_{0}^{1,2}(\R^N)\hookrightarrow L^{p}(\R^N, \frac{dx}{|x|^s})$ is compact.

\vskip0.1in

\noindent{\bf Case 2:  $0<s\leq \sigma, 2+\frac{4}{\sigma}\frac{\sigma-s}{N-2}< p<2^*(s)$.}   In this case, we can fix some $ \bar{s}\in (0,  s)$   such that
$$2+\frac{4}{\sigma}\frac{\sigma-{\bar{s}}}{N-2}< p<2^*(\bar{s}),$$
then by the embedding result, there exists some $C>0$ such that
(\ref{2014-6-5-we1}) holds. Then (\ref{2014-6-5-we2}) is also valid. Hence, in this case, we also obtain that $\frac{|u|^p}{|x|^s}$ is a tight sequence. Then similar to the previous case, we can also prove  the same conclusion.
\ep

\subsection{Applications on a single equation without critical terms}


Based  on the compact embedding result in Proposition \ref{2014-6-5-xl1}, we obtain the following results.

\bt\lab{2014-6-6-xcro1}
Let $\lambda>0$. Assume that  either
\begin{itemize}
\item [(i)]  $0\leq \sigma<s<2, 2<p<2^*(s)$ or
\item [(ii)] $0<s\leq \sigma<2, 2+\frac{4}{\sigma}\frac{\sigma-s}{N-2}< p<2^*(s)$,
\end{itemize}
then there exists a ground state solution to the following problem
\be\lab{2014-6-6-xe1}
\begin{cases}
-\Delta u+\lambda\frac{u}{|x|^\sigma}=\frac{u^{p-1}}{|x|^s}\;\hbox{in}\;\R^N,\\

u>0\;\hbox{in}\;\R^N, u\in D_{0}^{1,2}(\R^N)\cap L^2(\R^N,  \frac{dx}{|x|^\sigma}).
\end{cases}
\ee
\et
\bp
One can prove it by a modification of \cite[Lemma3.4]{ZhongZou.2015arXiv:1504.00433v1[math.AP]2Apr2015}.
\ep

\br\lab{2014-6-9-r2}
For  the equation (\ref{2014-6-6-xe1}), when $\sigma=s=0$, it has been studied in \cite{Kwong.1989}. The result of  Theorem \ref{2014-6-6-xcro1}  above covers the case of $\sigma\geq 0, s>0$.  But, the case of $\sigma>0, s=0$ is still open.
\er


\bt\lab{2014-6-6-xcro2}
Assume that  either
\begin{itemize}
\item [(1)] $0\leq \sigma<s<2, 2<p<2^*(s)$ or
\item [(2)] $0<s\leq \sigma<2, 2+\frac{4}{\sigma}\frac{\sigma-s}{N-2}< p<2^*(s)$,
\end{itemize}
then the best constant in (\ref{2014-6-5-xe4}) can be achieved.
 \et
\bp
This is a kind of CKN inequality with interpolation term. And the result is covered by \cite[Theorem 1.1]{ZhongZou.2015arXiv:1504.00433v1[math.AP]2Apr2015}.
\ep

In order to obtain multiple solutions,   we now recall  the the   Krasnoselskii Genus  defined  for studying  the  even functional $\Psi$ on a Banach space $E$. Set  $$\mathcal{A}: =\{A\subset E| A\;\hbox{is \, closed}, A=-A\}.$$
  For $A\subset \mathcal{A}, A\neq \emptyset$,  define
$$\gamma(A):=\begin{cases}\inf\{m:   \;\exists\;h\in C^0(A;\R^m\backslash \{0\}), h(-u)=-h(u)\}\\ \infty,\quad \hbox{if}\;\{\cdot \cdot\}=\emptyset,\;\hbox{in particular, if}\;A\ni 0, \end{cases}$$
and define $\gamma(\emptyset)=0$.
\vskip 0.2in
\noindent
{\bf Theorem A  (cf. \cite[Theomre 5.7]{Struwe.2008})} {\it Suppose $\Psi\in C^1(M)$ is an even functional on a complete symmetric $C^{1,1}-$manifold $M\subset E\backslash\{0\}$ in some Banach space $E$ and suppose $\Psi$ satisfies $(PS)$ condition and is bounded from below on $M$. Let $\tilde{\gamma}(M)=\sup\{\gamma(K):  K\subset M\;\hbox{compact and symmetric}\}$. Then the functional $\Psi$ possesses at least $\tilde{\gamma}(M)\leq \infty$ pairs of critical points.} \hfill$\Box$

\vskip0.11in

As a simple application of Theorem A, the following results are well known. Consider the problem
\be\lab{2014-6-16-e1}
\begin{cases}
-\Delta u+\lambda u=f(x, u)\;\quad&\hbox{in}\;\Omega,\\
u=0&\hbox{on}\;\partial\Omega,
\end{cases}
\ee
where $\Omega$ is a bounded domain in $\R^N$ and $f$ satisfies the following assumptions:
\begin{itemize}
\item[$(f_1)$]$f(x,-t)=-f(x, t)$;
\item[$(f_2)$]$f(x,u)=o(u)$ uniformly in $x$ as $u\rightarrow 0$;
\item[$(f_3)$]$f$ is continuous and $|f(x, u)|\leq a(1+|u|^{p-1})$ for some $a>0$ and $2<p<2^{*}$, where $2^{*}=2N/(N-2)$ for $N\geq 3$ and $2^{*}:=+\infty$ for $N=1, 2$.
\end{itemize}
Then for any $\lambda\geq 0$, problem (\ref{2014-6-16-e1}) admits infinitely many distinct pairs of solutions (see \cite[Themre 5.8]{Struwe.2008} for instance). We note that when $0\in \bar{\Omega}$ and $f(x, u)=\frac{|u|^{p-2}u}{|x|^{s}}$ with $s>0$, $(f_3)$ fails. We also note that the case of $|\Omega|=\infty$ will bring some truble. However, based   on  the  Proposition \ref{2014-6-5-xl1},  we can obtain the following result which seems to be new.

\bt\lab{2014-6-16-th1}
Let $\Omega\subset \R^N$ be an open set with $0\in \bar{\Omega}$ and $|\Omega|<\infty$ or $|\Omega|=\infty$. Furthermore,

\begin{itemize}
\item [(1)]  if $|\Omega|<\infty$, we assume that $0\leq s<2, 2<p<2^*(s);$
\item [(2)]  if $|\Omega|=\infty$, we assume that  either $0\leq \sigma<s<2, 2<p<2^*(s)$ or $0<s\leq \sigma<2, 2+\frac{4}{\sigma}\frac{\sigma-s}{N-2}< p<2^*(s)$.
\end{itemize}
Then equation
\be\lab{2014-6-16-e2}
\begin{cases}
-\Delta u+\lambda\frac{u}{|x|^\sigma}=\frac{|u|^{p-2}u}{|x|^s}\;\hbox{in}\;\Omega,\\
u=0\;\hbox{on}\;\partial\Omega, \quad u\in D_{0}^{1,2}(\Omega)\cap L^2(\Omega, \frac{dx}{|x|^\sigma}).
\end{cases}
\ee
possesses infinitely many distinct pairs of solutions for any $\lambda>-\lambda_{1,\sigma}(\Omega)$ if $|\Omega|<\infty$ and for any $\lambda>0$ if $|\Omega|=\infty$.
\et
\bp
On  $E= D_{0}^{1,2}(\Omega)\cap L^2(\Omega, \frac{dx}{|x|^\sigma})$,  let
$$\Psi(u):=\frac{1}{2}\int_\Omega (|\nabla u|^2+\lambda |u|^2)dx,  \quad  M:=\{u\in E:   \|u\|_{s,p}=1\}.$$
By Proposition \ref{2014-6-5-xl1}, we see that the continuous embedding $E\hookrightarrow L^p(\Omega,  \frac{dx}{|x|^\sigma})$ is compact. Hence, it is easy to see that $\Psi$ satisfies $(PS)$ condition and bounded from below on $M$.
By Theorem A, $\Psi$ admits infinitely many distinct pairs of critical points on $M$. After scaling, we obtain infinitely many distinct pairs of solutions for (\ref{2014-6-16-e2}).
\ep

\vskip0.3in

\s{Single equation involves critical exponent}
\renewcommand{\theequation}{3.\arabic{equation}}
\renewcommand{\theremark}{3.\arabic{remark}}
\renewcommand{\thedefinition}{3.\arabic{definition}}


\subsection{Splitting results for  a scalar  equation with unbounded domain}
It is well known that, when the nonlinear terms are critical, the corresponding functionals do not satisfy  the $(PS)$ condition. Usually, we need a better description of the behavior of  the $(PS)$ sequences of the corresponding functionals. For the case of a scalar equation with a bounded domain, some splitting results have been established  by Cerami-Zhong-Zou in  \cite[Theorem 3.1]{CeramiZhongZou.2015.DOI10.1007/s00526-015-0844-z.}, Ghoussoub-Kang in \cite[Theorem 3.1]{GhoussoubKang.2004}.
In this subsection, we are trying to extend this kind of results to the cases of unbounded domain in this subsection.  We consider the problem
\be\lab{2014-6-19we1}
\begin{cases}
&-\Delta u+\lambda \frac{|u|^{p-1}u}{|x|^{s_1}}=\frac{|u|^{2^*(s_2)-2}u}{|x|^{s_2}}\;\quad \hbox{in}\;\Omega,\\
& u(x)=0\;\hbox{on}\;\partial\Omega,
\end{cases}
\ee
 where $\Omega$ is an open unbounded domain in $\R^N$  $(N\geq 3).$
 When $1\leq p<2^*(s_1)-1$, if $\Omega=\R^N$ or a cone, it is easy to apply the Pohozaev identity to obtain the nonexistence of nontrivial solutions (see Corollary \ref{2014-6-21-cro1}). However, when $\Omega$ is unbounded but not a limit domain(i.e., $\Omega$ is not a cone), there seems to have no related results. Hence, we are aim to extend this work to the unbounded domain case. Firsly, we will prepare some results which will be useful.

\vskip0.1in

When $|\Omega|=\infty$ and $1\leq p<2^*(s_1)-1$, we can not ensure that $\displaystyle \int_\Omega \frac{|u|^{p+1}}{|x|^{s_1}}dx$ is well defined for all $u\in D_{0}^{1,2}(\Omega)$.   Therefore,  problem (\ref{2014-6-19we1})  have to be  posed in the framework of the Sobolev space
\be\lab{zouwm+b}E:=D_{0}^{1,2}(\Omega)\cap L^{p+1}\left(\Omega, \frac{dx}{|x|^{s_1}}\right).\ee  It is easy to see that $E$ is weak close due to the Fatou's Lemma. We take the norm in $E$
\be\lab{zouwm+a}\|u\|_E:=\left(\int_{\Omega}|\nabla u|^2dx\right)^{\frac{1}{2}}+\left(\int_\Omega \frac{|u|^{p+1}}{|x|^{s_1}}dx\right)^{\frac{1}{p+1}}.\ee
It is easy to see that when $|\Omega|<\infty$, $E=D_{0}^{1,2}(\Omega)$.
When $p=1$, it has been defined by \eqref{normE}, but for the convenience, even $p\neq 1$, we still prefer to adopt the denotation $E$.
To (\ref{2014-6-19we1}) there corresponds  energy  functional, defined in $E$,  is
\be\lab{zouwm+c} \varphi(u):=\frac{1}{2}\int_\Omega |\nabla u|^2 dx+\frac{\lambda}{p+1}\int_\Omega \frac{|u|^{p+1}}{|x|^{s_1}}dx-\frac{1}{2^*(s_2)}\int_\Omega \frac{|u|^{2^*(s_2)}}{|x|^{s_2}}dx,\ee
which is of class $C^2(E,\R)$.  Precisely, we will obtain the following splitting result. Compared with the bounded case, they are consistent in the presentation.

\bt\lab{2014-6-19-th1}(Splitting  theorem for unbounded domain $\Omega$)
Assume that $\Omega\subset \R^N$ is a $C^1$ open unbounded domain with $0\in \partial \Omega$ and $\partial\Omega$ is $C^2$ at $0$,  the mean curvature is negative, i.e., $H(0)<0$.   Suppose  that
$0\leq s_1<2, 0< s_2<2, 1\leq p<2^*(s_1)-1$.    Let   $\{u_n\}\subset E$  be  a bounded $(PS)_c$ sequence of the functional $\varphi$,  i.e., $\varphi(u_n)\rightarrow c$ and $\varphi'(u_n)\rightarrow 0$ strongly in $E^*$ as $n\rightarrow \infty$. Then there exists a critical point $U^0$ of $\varphi$, a number $k\in \NN$, $k$ functions $U^1, \cdots, U^k$ and $k$ sequences of radius $(r_n^j)_n, r_n^j>0, 1\leq j\leq k$ such that, up  to a subsequence if necessary,  the following properties are satisfied. Either
\begin{itemize}
\item[(a)] $u_n\rightarrow U^0$ in $E$ or
\item[(b)]  the following  items are true:
\begin{itemize}
\item[(b1)]$U^j\in D_{0}^{1,2}(\R_+^N)\subset D_{0}^{1,2}(\R^N)$ are nontrivial solutions of
\be\lab{2014-6-19-e1}
\begin{cases}
-\Delta u=\frac{|u|^{2^*(s_2)-2}u}{|x|^{s_2}}\quad &\hbox{in}\;\R_+^N,\\
u=0\;\hbox{on}\;\partial\R_+^N;
\end{cases}
\ee
\item[(b2)]$r_n^j\rightarrow 0$ as $n\rightarrow \infty$;\\
\item[(b3)] $\left\|u_n-U^0-\sum_{j=1}^{k} (r_n^j)^{\frac{2-N}{2}}U^j(\frac{\cdot}{r_n^j})\right\|\rightarrow 0$, where $\|\cdot\|$ is the norm in $D^{1,2}(\R^N)$;\\
\item[(b4)]$\|u_n\|^2\rightarrow \sum_{j=0}^{k}\|U^j\|^2$;\\
\item[(b5)] $\varphi(u_n)\rightarrow \varphi(U^0)+\sum_{j=1}^{k}I(U^j)$ with $$I(U^j)\geq \left(\frac{1}{2}-\frac{1}{2^*(s_2)}\right)\mu_{s_2}(\R_+^N)^{\frac{2^*(s_2)}{2^*(s_2)-2}},$$
    and
    $$I(u):=\frac{1}{2}\int_{\R_+^N}|\nabla u|^2dx-\frac{1}{2^*(s_2)}\int_{\R_+^N}\frac{|u|^{2^*(s_2)}}{|x|^{s_2}}dx \;\; \hbox{ for }\;u\in D_{0}^{1,2}(\R_+^N).$$
\end{itemize}
\end{itemize}
\et

\br\lab{2015-4-24-r1}
What surprising us is that the representation of the splitting result for unbounded domain are exactly the same as that one for bounded domain case (see \cite[Theorem 3.1]{CeramiZhongZou.2015.DOI10.1007/s00526-015-0844-z.}). However, on the applications,the assumptions we required  are distinct to ensure the ``bounded $PS$ sequence" exists. For example, in Theorem , when $|\Omega|=\infty$, we always assume that $\lambda>0$.
\er

In order to prove the above theorem, we need to prepare  beforehand several lemmas.

\bl\lab{2014-3-19-l1}(See \cite[Lemma 3.2, Remark 3.2]{CeramiZhongZou.2015.DOI10.1007/s00526-015-0844-z.} and \cite[Lemma 3.3]{GhoussoubKang.2004})
 Let $\Omega$ be an open subset of $\R^N$.
 \begin{itemize}
 \item[1)]Assume that $\{u_n\}\subset L^{p}(\Omega, \frac{dx}{|x|^s}), 1\leq p<\infty$. If $\{u_n\}$ is bounded in $L^{p}(\Omega, \frac{dx}{|x|^s})$ and $u_n\rightarrow u$ almost everywhere on $\Omega$, then
 $$ \int_\Omega \frac{|u_n|^p}{|x|^s}dx-\int_\Omega \frac{|u_m-u|^p}{|x|^s}dx\rightarrow \int_\Omega \frac{|u|^p}{|x|^s}dx\;\hbox{as}\; n\rightarrow \infty.$$
 \item[2)]Assume $\{u_n\}\subset D_{0}^{1,2}(\Omega)$ is such that
 $u_n\rightharpoonup u$ in $D_{0}^{1,2}(\Omega)$, then
 $$\displaystyle \int_\Omega |\nabla u_n|^2dx-\int_\Omega |\nabla u_n-\nabla u|^2dx\rightarrow \int_\Omega |\nabla u|^2dx
 \hbox{   as   }   n\rightarrow \infty.  $$
 \item[3)]  If $u_n\rightarrow u$ weakly in $D_{0}^{1,2}(\R^N)$, then
     $$\frac{|u_n|^{2^*(s)-2}u_n}{|x|^s}-\frac{|u_n-u|^{2^*(s)-2}(u_n-u)}{|x|^s}\rightarrow \frac{|u|^{2^*(s)-2}u}{|x|^s}  \hbox{   as   }   n\rightarrow \infty $$ in $H^{-1}(\R^N)$.
 \end{itemize}
\el\hfill$\Box$
\bl\lab{2014-6-20-l1}(cf. \cite[Lemma  2.1]{CeramiZhongZou.2015.DOI10.1007/s00526-015-0844-z.})
 Let   $\Omega\subset \R^N$ be a set having Lebesgue measure $|\Omega|<\infty.$ Let  $0\leq s< 2$, $1\leq p<2^*(s)$. Then the embedding $H_0^1(\Omega)\hookrightarrow L^{p}(\Omega,\frac{dx}{|x|^s})$ is compact.
\el

\bl\lab{2014-6-19-l1}
Consider $\{u_n\} \subset E$ such that $I(u_n)\rightarrow c$ and $I'(u_n)\rightarrow 0$ in $E^*$. For $r_n \in (0,\infty)$ with $r_n\rightarrow 0$, we assume that the rescaled sequence
$v_n(x):=r_{n}^{\frac{N-2}{2}}u_n(r_nx)$ is such that $v_n\rightarrow v$ weakly in $D_{0}^{1,2}(\R^N)$ and $v_n\rightarrow v$ a.e. on $\R^N$.
Then, $I'(v)=0$ and the sequence
$$w_n(x):=u_n(x)-r_{n}^{\frac{2-N}{2}}v(\frac{x}{r_n})$$
satisfies $I(w_n)\rightarrow c-I(v), I'(w_n)\rightarrow 0$ in $E^*$ and $\|w_n\|^2=\|u_n\|^2-\|v\|^2+o(1)$.
\el
\bp
We note that when $\Omega$ is a bounded domain, this lemma is given in \cite[Lemma 3.3]{CeramiZhongZou.2015.DOI10.1007/s00526-015-0844-z.}, a similar result we refer to \cite[Lemma 3.4]{GhoussoubKang.2004}. Following the procedure  carefully, we see that the proof of  \cite[Lemma 3.3]{CeramiZhongZou.2015.DOI10.1007/s00526-015-0844-z.} is still valid when  $\Omega$ is unbounded. Hence, we omit the details. We shall  give the detailed proof  about this type of result for the system case, see the forthcoming Lemma \ref{2014-6-10-l1}.
\ep

Indeed, compared   with the bounded case, the biggest difficulty is that how to ensure $r_n\rightarrow 0$ which is a necessary condition to apply Lemma \ref{2014-6-19-l1} to continue the splitting iteration.   When $|\Omega|=\infty$, let $u_n$ be a bounded $(PS)_c$ sequence for $\varphi$, and $u_n\rightharpoonup u$ in $E$,  it is not trivial that $I'(u_n-u)\rightarrow 0$ in $E^*$. To overcome these difficulties, we  require the following results:
\bl\lab{2014-6-18-xl4}  Consider the equation (\ref{2014-6-19we1}).
Assume that $0\leq s_1<2, 0<s_2< 2, p<2^*(s_1)-1$, then
$\displaystyle \left\{\frac{|u_n|^{2^*(s_2)}}{|x|^{s_2}}\right\}$ is a tight sequence if
$\{u_n\}\subset E$ is a bounded sequence.
\el
\bp
Since $p+1<2^*(s_1)$, we have
$$2N-2s_1-N(p+1)+2(p+1)>0.$$
Since $s_2>0$, we can take $ \sigma\in (0, s_2)$   such that
\be\lab{2014-6-18-e5}
\begin{cases}
\frac{2-\sigma}{s_2-\sigma}\geq \frac{4-2s_1}{2N-2s_1-N(p+1)+2(p+1)},\\
\sigma-\frac{2s_1(s_2-\sigma)}{2N-2s_1-N(p+1)+2(p+1)}>0,\\
1-\frac{2(s_2-\sigma)}{2N-2s_1-N(p+1)+2(p+1)}>0,\\
\frac{2(p+1)(s_2-\sigma)}{2N-2s_1-N(p+1)+2(p+1)}<2^*(s_2)-1.
\end{cases}
\ee
Define
\be\lab{2014-6-18-e6}
\begin{cases}
p_1:=\frac{2(p+1)(s_2-\sigma)}{2N-2s_1-N(p+1)+2(p+1)},\\
p_2:=2^*(s_2)-(p_1+1)-1,\\
\sigma_1:=\frac{2s_1(s_2-\sigma)}{2N-2s_1-N(p+1)+2(p+1)},\\
\sigma_2:=\sigma-\sigma_1,\\
\bar{\sigma}:=\frac{(p+1)\sigma_2}{p-p_1}
\end{cases}
\ee
then we have
\be\lab{2014-6-18-e7}
\begin{cases}
0<p_1+1<\min\{p+1,2^*(s_2)-1\},\\
0<p_2+1, p_1+p_2=2^*(s_2)-2,\\
\sigma_1>0, \sigma_2>0, \sigma_1+\sigma_2=\sigma,\\
0<\bar{\sigma}\leq 2,\\
\frac{(p+1)(p_2+1)}{p-p_1}=2^*(\bar{\sigma}).
\end{cases}
\ee
Then by  the  H\"{o}lder inequality, for any $\Lambda\subset \Omega$, we have
\begin{align*}
\int_{\Lambda}\frac{|u|^{2^*(s_2)}}{|x|^{\sigma}}dx=&\int_{\Lambda} \frac{|u|^{p_1+1}}{|x|^{\sigma_1}} \frac{|u|^{p_2+1}}{|x|^{\sigma_2}}dx\\
\leq&\left(\int_\Lambda \frac{|u|^{p+1}}{|x|^{\sigma_1 \frac{p+1}{p_1+1}}}dx\right)^{\frac{p_1+1}{p+1}}\left(\int_{\Lambda} \frac{|u|^{(p_2+1)\frac{p+1}{p-p_1}}}{|x|^{\sigma_2\frac{p+1}{p-p_1}}}dx\right)^{\frac{p-p_1}{p+1}}\\
=&\left(\int_\Lambda \frac{|u|^{p+1}}{|x|^{s_1}}dx\right)^{\frac{p_1+1}{p+1}}\left(\int_{\Lambda} \frac{|u|^{2^*(\bar{\sigma})}}{|x|^{\bar{\sigma}}}dx\right)^{\frac{p-p_1}{p+1}}.
\end{align*}
Hence, by the Hardy-Sobolev inequality, we see that
$|u_n|_{\sigma, 2^*(s_2)}$ is bounded due to the boundedness of $\{u_n\}$ in $E$.
Hence,
\begin{align*}
\int_{\mathbb{B}_R^c\cap \Omega} \frac{|u_n|^{2^*(s_2)}}{|x|^{s_2}}dx=&\int_{\mathbb{B}_R^c\cap \Omega} \frac{1}{|x|^{s_2-\sigma}} \frac{|u_n|^{2^*(s_2)}}{|x|^{\sigma}}dx\\
\leq&\frac{1}{R^{s_2-\sigma}}\int_{\mathbb{B}_R^c\cap \Omega} \frac{|u_n|^{2^*(s_2)}}{|x|^{\sigma}}dx\\
\rightarrow&0\;\hbox{as}\;R\rightarrow \infty,
\end{align*}
which implies that $\displaystyle\left\{\frac{|u_n|^{2^*(s_2)}}{|x|^{s_2}}\right\}$ is a tight sequence.
\ep

\bl\lab{2014-6-20-l2}   Consider the equation (\ref{2014-6-19we1}).
Assume that $0\leq s_1<2, 0<s_2< 2, p<2^*(s_1)-1$  and that  $\{u_n\}\subset E$ is a bounded $(PS)_c$ sequence for $\varphi(u)$, then $\displaystyle\left\{|\nabla u_n|^2+\lambda \frac{|u_n|^{p+1}}{|x|^{s_1}}\right\}$ is a tight sequence.
\el
\bp
Let $\chi_R(x)\in C^\infty(\R^N)$ be a cutoff function such that $0\leq \chi_R\leq 1$, $\chi_R(x)\equiv 0$ in $\mathbb{B}_R$, $\chi_R(x)\equiv 1$ in $\mathbb{B}_{2R}^c$ and $|\nabla \chi_R|\leq \frac{2}{R}$. Then it is easy to see that $u_n\chi_R\in E$. When $\{u_n\}$ is a $(PS)_c$ sequence, we have
\be\lab{2014-6-20-e1}
\int_{\Omega}\nabla u_n\cdot \nabla (u_n\chi_R)+\lambda \int_\Omega \frac{|u_n|^{p+1}\chi_R}{|x|^{s_1}}-\int_{\Omega}\frac{|u_n|^{2^*(s_2)}\chi_R}{|x|^{s_2}}=o(1)\|u_n\|_E.
\ee
Since $\{u_n\}$ is bounded, going to a subsequence if necessary, we may assume that $u_n\rightharpoonup u$ in $E$, $u_n\rightarrow u$ a.e in $\R^N$, $u_n\rightarrow u$ in $L_{loc}^{\gamma}(\R^N)$ for $\gamma\in [1,2^*)$. Then for any fixed $R>0$, we see that
\be\lab{2014-6-20-e2}
\lim_{n\rightarrow \infty} \int_\Omega u_n\nabla u_n\cdot\nabla \chi_R=\int_\Omega u\nabla u\cdot \nabla\chi_R.
\ee
By the absolute continuity of the integral, we have
\be\lab{2014-6-20-e3}
\lim_{R\rightarrow \infty}\int_\Omega u\nabla u\cdot \nabla\chi_R=0,
\ee
because of that
\begin{align*}
\left|\int_\Omega u\nabla u\cdot \nabla\chi_R\right|\leq& 2\int_{R\leq|x|\leq 2R}\frac{|u|}{R}|\nabla u|
\leq 8\left(\int_{R\leq|x|\leq 2R}\frac{u^2}{|x|^2}\right)^{\frac{1}{2}}\left(\int_{R\leq |x|\leq 2R}|\nabla u|^2\right)^{\frac{1}{2}}\\
\leq&C \int_{R\leq|x|\leq 2R}|\nabla u|^2dx.
\end{align*}
By (\ref{2014-6-20-e1}), we have
\be\lab{2014-6-20-e4}
\int_\Omega \left(|\nabla u_n|^2+\lambda \frac{|u_n|^{p+1}}{|x|^{s_1}}\right)\chi_R=o(1)\|u_n\|_E+\int_\Omega \frac{|u_n|^{2^*(s_2)}}{|x|^{s_2}}\chi_R-\int_\Omega u_n\nabla u_n\cdot\nabla \chi_R.
\ee
Recalling that $\{u_n\}$ is bounded in $E$, by Lemma \ref{2014-6-18-xl4}, we obtain that $\displaystyle \left\{\frac{|u_n|^{2^*(s_2)}}{|x|^{s_2}}\right\}$ is a tight sequence, which implies that
\be\lab{2014-6-20-e5}
\lim_{R\rightarrow \infty} \int_\Omega \frac{|u_n|^{2^*(s_2)}}{|x|^{s_2}}\chi_R=0\;\hbox{uniformly for all}\; n.
\ee
By (\ref{2014-6-20-e2})-(\ref{2014-6-20-e5}), we have
\be\lab{2014-6-20-e6}
\lim_{R\rightarrow \infty}\limsup_{n\rightarrow \infty}\int_\Omega \left(|\nabla u_n|^2+\lambda \frac{|u_n|^{p+1}}{|x|^{s_1}}\right)\chi_R=0
\ee
and it following that
$\displaystyle\left\{|\nabla u_n|^2+\lambda \frac{|u_n|^{p+1}}{|x|^{s_1}}\right\}$ is a tight sequence.
\ep

\vskip0.1in

\bl\lab{2014-6-20-l3}
Assume that $0\leq s_1<2, 0<s_2< 2, 1\leq p<2^*(s_1)-1$, and  that $\{u_n\}\subset E$ is a bounded $(PS)_c$ sequence for $\varphi(u)$ such that $u_n\rightharpoonup u$ in $E$, then up to a subsequence,
$$\frac{|u_n|^{p-1}u_n}{|x|^{s_1}}\rightarrow \frac{|u|^{p-1}u}{|x|^{s_1}}\;\hbox{in}\;E^*\;\hbox{as}\;n\rightarrow \infty.$$
\el
\bp
For any $h\in E$ and $\varepsilon>0$, by Lemma \ref{2014-6-20-l2} and H\"{o}lder inequality, there exists some $R_\varepsilon>0$ such that
\begin{align}\lab{2014-6-20-e7}
&\left|\int_{\Omega\cap \mathbb{B}_{R_\varepsilon}^{c}} \left(\frac{|u_n|^{p-1}u_n-|u|^{p-1}u}{|x|^{s_1}}\right)hdx\right|\nonumber\\
& \leq\int_{\Omega\cap \mathbb{B}_{R_\varepsilon}^{c}}\left|\frac{|u_n|^{p-1}u_n}{|x|^{s_1}}h\right|dx +\int_{\Omega\cap\mathbb{B}_{R_\varepsilon}^{c}}\left|\frac{|u|^{p-1}u}{|x|^{s_1}}h\right|dx\nonumber\\
& \leq \left(\int_{\Omega\cap \mathbb{B}_{R_\varepsilon}^{c}}\frac{|u_n|^{p+1}}{|x|^{s_1}}\right)^{\frac{p}{p+1}} \|h\|_{E}+\left(\int_{\Omega\cap \mathbb{B}_{R_\varepsilon}^{c}}\frac{|u|^{p+1}}{|x|^{s_1}}\right)^{\frac{p}{p+1}} \|h\|_{E}\nonumber\\
& < \varepsilon \|h\|_{E}.
\end{align}
On the other hand, by Lemma \ref{2014-6-20-l1} and H\"{o}lder inequality, we have
\be\lab{2014-6-20-e8}
\left|\int_{\Omega\cap \mathbb{B}_{R_\varepsilon}} \left(\frac{|u_n|^{p-1}u_n-|u|^{p-1}u}{|x|^{s_1}}\right)hdx\right|=o(1)\|h\|_E\;\hbox{as}\;n\rightarrow \infty.
\ee
By (\ref{2014-6-20-e7}) and (\ref{2014-6-20-e8}), we prove that
$$\frac{|u_n|^{p-1}u_n}{|x|^{s_1}}\rightarrow \frac{|u|^{p-1}u}{|x|^{s_1}}\;\hbox{in}\;E^*\;\hbox{as}\;n\rightarrow \infty.$$
\ep

\bc\lab{2014-6-20-cro1}
Assume that $0\leq s_1<2, 0<s_2<2, 1\leq p<2^*(s_1)-1$.  Let $\{u_n\}\subset E$ be  a bounded $(PS)_c$ sequence for $\varphi(u)$ such that $u_n\rightharpoonup u$ in $E$, then
$$\int_\Omega \frac{|u_n|^{p+1}}{|x|^{s_1}}dx\rightarrow \int_\Omega \frac{|u|^{p+1}}{|x|^{s_1}}dx\;\hbox{as}\;n\rightarrow \infty$$
and
$$u_n\rightarrow u\;\hbox{in}\;L^{p+1}\left(\Omega, \frac{dx}{|x|^{s_1}}\right).$$
\ec
\bp
Rewrite
$$\int_\Omega \frac{|u_n|^{p+1}-|u|^{p+1}}{|x|^{s_1}}=\int_\Omega \frac{\left(|u_n|^{p-1}u_n-|u|^{p-1}u\right)u_n}{|x|^{s_1}}+\int_\Omega \frac{|u|^{p-1}u(u_n-u)}{|x|^{s_1}}.$$
Since $\{u_n\}\subset E$ is a bounded $(PS)_c$ sequence for $\varphi(u)$, by Lemma \ref{2014-6-20-l3}, we have
$$\int_\Omega \frac{|u_n|^{p+1}}{|x|^{s_1}}dx\rightarrow \int_\Omega \frac{|u|^{p+1}}{|x|^{s_1}}dx\;\hbox{as}\;n\rightarrow \infty.$$
Combined with Lemma  \ref{2014-3-19-l1}, we obtain that
$$u_n\rightarrow u\;\hbox{in}\;L^{p+1}\left(\Omega, \frac{dx}{|x|^{s_1}}\right).$$
\ep

\vskip 0.3in
\noindent{\bf Proof of Theorem \ref{2014-6-19-th1}. }
The proof is very like to \cite[The proof of Theorem 3.1]{CeramiZhongZou.2015.DOI10.1007/s00526-015-0844-z.}. However, since $\Omega$ is unbounded, there will be some differences in the details. Hence, for the convenience of the readers, we will give the details.

Let $\{u_n\}\subset E$ be a bounded $(PS)_c$ sequence of $\varphi(u)$, that is, $\varphi(u)\rightarrow c$ and $\varphi'(u_n)\rightarrow 0$ in $E^*.$ Recalling that $E\subset D_{0}^{1,2}(\Omega)$ is weak close,
then  there exists $U^0\in H_0^1(\Omega)$ so that, up to a subsequence,   $u_n\rightharpoonup U^0$ in $E$   and $ u_n\rightarrow  U^0$ a.e. on $\R^N$. For any $h\in E$, by Lemma \ref{2014-6-20-l3} we have
$$\int_\Omega \frac{|u_n|^{p-1}u_nh}{|x|^{s_1}}dx\rightarrow \int_\Omega \frac{|U^0|^{p-1}U^0 h}{|x|^{s_1}}dx\;\hbox{as}\;n\rightarrow \infty.$$
Combined with Lemma \ref{2014-3-19-l1}, we see that $\langle \varphi'(U^0), h\rangle=0$ since $\varphi'(u_n)\rightarrow 0$ in $E^*$. Thus, $\varphi'(U^0)=0$.
Let $u_n^1:=u_n-U^0$, if $|\Omega|<\infty$, by Lemma \ref{2014-3-19-l1} and Lemma \ref{2014-6-20-l1},
it is easy to see that the sequence $u_n^1:=u_n-U^0$ satisfies
\be\lab{2014-3-21-e1}
\begin{cases}
\|u_n^1\|=\|u_n^1\|_E+o(1)=\|u_n\|_E-\|U^0\|_E+o(1)=\|u_n\|-\|U^0\|+o(1),\\
I'(u_n^1)\rightarrow 0\;\hbox{in}\;H^{-1}(\Omega),\\
I(u_n^1)\rightarrow c_1:=c-\varphi(U^0),
\end{cases}
\ee
where $\|u\|:=\big(\int_\Omega |\nabla u|^2 dx\big)^{\frac{1}{2}}$.
If $|\Omega|=\infty$, by Lemma \ref{2014-6-20-l3}, Corollary \ref{2014-6-20-cro1} and Lemma \ref{2014-3-19-l1}, we can also obtain the results of (\ref{2014-3-21-e1}).

If $u_n^1\rightarrow 0$ in $D_{0}^{1,2}(\Omega)$, by Corollary \ref{2014-6-20-cro1}, we obtain that $u_n\rightarrow U^0$ in $E$ and we are done. If $u_n^1\not\rightarrow 0$ in $D_{0}^{1,2}(\Omega)$, we claim that
\be\lab{2014-3-21-e2}
\liminf_{n\rightarrow \infty} \int_\Omega \frac{|u_n^1|^{2^*(s_2)}}{|x|^{s_2}}dx>0.
\ee
Otherwise, the facts  that  $\{u_n^1\}\subset E\subset D_{0}^{1,2}(\Omega)$ is bounded and $I'(u_n^1)\rightarrow 0$ in $H^{-1}(\Omega)$ would bring to a contradiction because
$$\langle I'(u_n^1), u_n^1\rangle =\int_\Omega |\nabla u_n^1|^2 dx -\int_\Omega \frac{|u_n^1|^{2^*(s_2)}}{|x|^{s_2}}dx\rightarrow 0,$$
and, hence, $u_n^1\rightarrow 0$ in $D_{0}^{1,2}(\Omega)$.  Thereby, (\ref{2014-3-21-e2}) holds true.

Define an analogue of Levy's concentration function
$$Q_n(r):=\int_{\mathbb{B}_r}\frac{|u_n^1|^{2^*(s_2)}}{|x|^{s_2}}dx.$$
Since $Q_n(0)=0$ and, due to (\ref{2014-3-21-e2}), $Q_n(\infty)>0$ , we can choose  $$\displaystyle\delta_1<\frac{1}{2}\liminf_{n\rightarrow \infty} \int_\Omega \frac{|u_n^1|^{2^*(s_2)}}{|x|^{s_2}}dx$$ so small that $0<\delta_1<\left(\frac{\mu_{s_2}(\R^N)}{2}\right)^{\frac{N-s_2}{2-s_2}}$
and find a positive sequence $\{r_n^1\}$ such that $Q_n(r_n^1)=\delta_1$.
Set $v_n^1(x):=(r_n^1)^{\frac{N-2}{2}}u_n^1(r_n^1x)$. Since $\|v_n^1\|=\|u_n^1\|$ is bounded, we may assume $v_n^1\rightharpoonup U^1$ in $D^{1,2}(\R^N)$, $v_n^1\rightarrow U^1$ a.e. on $\R^N$ and
\be\lab{2014-6-21-xe1}
\delta_1=\int_{\mathbb{B}_1}\frac{|v_n^1|^{2^*(s_2)}}{|x|^{s_2}}dx.
\ee
We claim that $U^1\not\equiv 0$.

 Set $\Omega_n=\frac{1}{r_n^1}\Omega$ and let $f_n\in D_{0}^{1,2}(\Omega)$ be such that   $\displaystyle \langle I'(u_n^1), h\rangle=\int_\Omega \nabla f_n\cdot \nabla h$
for any $h\in D_{0}^{1,2}(\Omega)$. Then $g_n(x):=(r_n^1)^{\frac{N-2}{2}}f_n(r_n^1x)$ satisfies
$\displaystyle \int_{\Omega_n}|\nabla g_n|^2 =\int_\Omega |\nabla f_n|^2$ and $\displaystyle \langle I'(v_n^1), h\rangle=\int_{\Omega_n}\nabla g_n\cdot \nabla h$ for any $h\in D_{0}^{1,2}(\Omega_n)$.

 If $U^1\equiv 0$ would be true, then we could deduce that,
for any $h\in C_0^\infty(\R^N)$ with $supp\;h\subset \mathbb{B}_1,$
\begin{align*}
\int_{\mathbb{B}_1}|\nabla (hv_n^1)|^2=&\int_{\mathbb{B}_1}\nabla v_n^1\cdot \nabla (h^2v_n^1)+o(1)\\
=&\int_{\mathbb{B}_1}\frac{h^2|v_n^1|^{2^*(s_2)}}{|x|^{s_2}}+\int_{\mathbb{B}_1}\nabla g_n\cdot \nabla(h^2v_n^1)+o(1)\\
\leq& \mu_{s_2}(\R^N)^{-1}\left(\int_{supp\; h}\frac{|v_n^1|^{2^*(s_2)}}{|x|^{s_2}}\right)^{\frac{2^*(s_2)-2}{2^*(s_2)}}
\int_{\mathbb{B}_1} \left|\nabla (hv_n^1)\right|^2+o(1)\\
=&\mu_{s_2}(\R^N)^{-1}\delta_1^{\frac{2-s_2}{N-s_2}}\int_{\mathbb{B}_1} |\nabla (hv_n^1)|^2+o(1)\\
<&\frac{1}{2}\int_{\mathbb{B}_1} |\nabla (hv_n^1)|^2+o(1),
\end{align*}
here we are  using the inequality:
$$\int_{\R^N} \frac{h^2|u|^{2^*(s)}}{|x|^s}\leq \mu_s(\R^N)^{-1}\left(\int_{supp\;h}\frac{|u|^{2^*(s)}}{|x|^s}\right)^{\frac{2^*(s)-2}{2^*(s)}}\int_{\R^N} |\nabla (hu)|^2$$
that holds for all $u\in D^{1,2}(\R^N)$ and $h\in C_0^\infty(\R^N)$ (see \cite[Lemma 3.5.]{GhoussoubKang.2004}).

\vskip0.1in

Hence, $\nabla v_n^1\rightarrow 0$ in $L^{2}_{loc}\left(\mathbb{B}_1\right)$ and, by the Hardy-Sobolev inequality,  $v_n^1\rightarrow 0$ in $L^{2^*(s_2)}_{loc}\left(\mathbb{B}_1, \frac{dx}{|x|^{s_2}}\right)$, which follows a  contradiction with  \eqref{2014-6-21-xe1}.  Thus, the claim is proved.

\vskip0.1in

Now, let us prove that $r_n^1\rightarrow 0$. We proceed by contradiction.
Otherwise, by the choice of $\delta_1$ and Lemma \ref{2014-6-18-xl4}, we obtain that $\{r_n^1\}$ is bounded. Hence,
we may assume that $r_n^1\rightarrow r_\infty^1>0$, then the fact that $u_n^1\rightharpoonup 0$ in $D_{0}^{1,2}(\Omega)$ would imply  $v_n^1(x):=(r_n^1)^{\frac{N-2}{2}}u_n^1(r_n^1x)\rightharpoonup 0$ in $D_{0}^{1,2}(\R^N)$,  in contradiction with  $U^1\not\equiv 0$. Therefore,  $r_n^1\rightarrow 0$.

Next, we prove that $supp\;U^1\subset \R_+^N$. Without loss of generality, we assume  that   $\partial R_+^N=\{x_N=0\}$ is tangent to $\partial\Omega$ at $0$, and that $-e_N=(0,\cdots, -1)$ is the outward normal to $\partial\Omega$ at that point. For any compact $K\subset \R_-^N$, we have for $n$ large enough, that $\frac{\Omega}{r_n^1}\cap K=\emptyset$   as $r_n^1\rightarrow 0$. Since $supp\;v_n^1\subset \frac{\Omega}{r_n^1}$ and $v_n^1\rightarrow U^1$ a.e. in $\R^N$,  $U^1=0$ a.e. on $K$ follows, and, therefore, $supp\;U^1\subset \R_+^N$.
\vskip0.1in

By (\ref{2014-3-21-e1}) and Lemma \ref{2014-6-19-l1}, we obtain that $I'(U^1)=0$. Hence, $U^1$ is a weak solution of (\ref{BP}).  Moreover, by Lemma \ref{2014-6-19-l1} again, we see  that the sequence $u_n^2(x):=u_n^1(x)-(r_n^1)^{\frac{2-N}{2}}U^1(\frac{x}{r_n^1})$ also satisfies
\be\lab{2014-3-21-xe2}
\begin{cases}
\|u_n^2\|^2=\|u_n\|^2-\|U^0\|^2-\|U^1\|^2+o(1), \\
I(u_n^2)\rightarrow c_2:=c-\varphi(U^0)-I(U^1),\\
I'(u_n^2)\rightarrow 0\;\hbox{in}\;H^{-1}(\Omega).
\end{cases}
\ee
We also note that
$$I(U^1)\geq \left(\frac{1}{2}-\frac{1}{2^*(s_2)}\right)
\mu_{s_2}(\R_+^N)^{\frac{2^*(s_2)}{2^*(s_2)-2}}.$$
By iterating the above procedure, we construct critical points $U^j$ of $I(u)$ and sequences $(r_n^j)$ with the above properties. Here we note that since $\int_\Omega \frac{|u_n^i|^{p+1}}{|x|^{s_1}}dx=o(1), \|u_n^i\|=\|u_n^i\|_E+o(1)$, we see that $\{u_n^i\}$ is also can be viewed as a bounded $(PS)_{c^i}$ for $\varphi(u)$ with $c_1=c-\varphi(U^0)$ and $c_i=c-\varphi(U^0)-I(U^1)-\cdots-I(U^{i-1}), i\geq 2$. Hence, by Lemma \ref{2014-6-18-xl4}, we can guarantee the boundedness of $\{r_n^i\}$.
 And since $\varphi(u_n) \rightarrow c,$  the iteration must terminate after a finite number of steps.\hfill$\Box$

\vskip0.2in
\subsection{Existence results on a single equation}
\subsubsection{A supplement to a single equation on the domain with finite Lebesgue measure}
Firstly we will study the case of $p=1$ basing on the result of Proposition \ref{2014-6-4-wl1} as a supplement to the results of \cite{CeramiZhongZou.2015.DOI10.1007/s00526-015-0844-z.}(see \cite[Remark 2.1]{CeramiZhongZou.2015.DOI10.1007/s00526-015-0844-z.}). That is, we consider the following problem
\be\lab{2015-4-24-we1}
\begin{cases}
-\Delta u+\lambda^*\frac{u}{|x|^{s_1}}=\lambda \frac{|u|^{2^*(s_2)-2}u}{|x|^{s_2}}\quad&\hbox{in}\;\Omega\\
u(x)=0&\hbox{on}\;\partial\Omega,
\end{cases}
\ee
After rescaling properly, we may assume that $\lambda=1$.
We obtain Theorem \ref{2015-4-24-th1} and Theorem \ref{2015-4-24-th2}.

\noindent
{\bf Proof of Theorem \ref{2015-4-24-th1}:}
When $\Omega$ is a bounded domain, it is just a special case of \cite[Theorem 1.4]{CeramiZhongZou.2015.DOI10.1007/s00526-015-0844-z.}. We only need to prove that $u>0$. Since the process is very like, we prefer to give the detail of this part in the proof of Theorem \ref{2015-4-24-th2}. When $\Omega$ is unbounded, we  only need to replace \cite[Theorem 3.1]{CeramiZhongZou.2015.DOI10.1007/s00526-015-0844-z.} by Theorem \ref{2014-6-19-th1} above, by the similar arguments as the bounded case, we can obtain the results.
\hfill$\Box$

\noindent
{\bf Proof of Theorem \ref{2015-4-24-th2}:}
Without loss of generality, we may assume that $\lambda=1$. Define
$$\varphi(u):=\int_{\Omega}\left[\frac{|\nabla u|^2}{2}+\lambda^*\frac{u^2}{2}-\frac{|u|^{2^*(s_1)}}{2^*(s_1)}\right]dx,$$
which is of class $C^2(H_0^1(\Omega), \R)$.  Since $\lambda^*>-\lambda_{1,s_1}(\Omega)$, by Proposition   \ref{2014-6-4-wl1}, we can choose the norm
$$\|u\|:=\left(\int_\Omega |\nabla u|^2+\lambda^*\frac{|u|^2}{|x|^{s_1}}\right)^{\frac{1}{2}}$$
which is equivalent to the usual one in $H_0^1(\Omega)$.
Recalling that $2^*(s_2)>2$, it is easy to verify that $\varphi$ has the mountain pass geometry. Now, we define the mountain pass value
\be\lab{2014-6-18-e1}
c^*:=\inf_{\gamma\in \Gamma}\max_{t\in [0,1]}\varphi(\gamma(t)),
\ee
where $\Gamma:=\{\gamma(t)\in C([0,1],  H_0^1(\Omega)): \gamma(0)=0, \varphi(\gamma(1))<0\}.$ Then there exists a $(PS)_{c^*}$ sequence, that is, a sequence $\{u_n\}_n,\ \ u_n \in H_0^1(\Omega),$ such that
$$\begin{cases}
 \Phi(u_n)\rightarrow c^*,\\
 \Phi'(u_n)\rightarrow 0\quad \hbox{in}\;H_{0}^{-1}(\Omega).
 \end{cases}$$
By the Hardy-Sobolev inequality, we see that $c^*>0$. On the other hand, following the processes of \cite[Lemma 2.4.]{CeramiZhongZou.2015.DOI10.1007/s00526-015-0844-z.}, we can prove that there exists some $v_0\in H_0^1(\Omega)\backslash\{0\}$ such that $\varphi(v_0)<0$ and
$$\max_{0\leq t\leq 1}\varphi(tv_0)<\left(\frac{1}{2}-\frac{1}{2^*(s_1)}\right)\mu_{s_1}(\R_+^N)^{\frac{2^*(s_1)}{2^*(s_1)-2}}.$$
Hence, we have
$$0<c^*<\left(\frac{1}{2}-\frac{1}{2^*(s_1)}\right)\mu_{s_1}(\R_+^N)^{\frac{2^*(s_1)}{2^*(s_1)-2}}.$$
It is easy to check that $\{u_n\}$ is bounded in $H_0^1(\Omega)$. Then by \cite[Theorem 3.1]{CeramiZhongZou.2015.DOI10.1007/s00526-015-0844-z.} if $\Omega$ is bounded and by Theorem \ref{2014-6-19-th1} if $\Omega$ is unbounded, there exists some $U^0\in H_0^1(\Omega)$ such that $u_n\rightarrow U^0$ in $H_0^1(\Omega)$ and $\varphi'(U^0)=0, \varphi(U^0)=c^*$. Hence, $U^0$ is a solution of  \eqref{2015-4-24-we1}.
Now, let $$\mathscr{A}:=\{u\neq 0, \varphi'(u)=0\},$$
then $\mathscr{A}\neq \emptyset$. Define $\displaystyle c_0=\inf_{u\in \mathscr{A}} \varphi(u)$, then by Hardy-Sobolev inequality, we see that $c_0>0$. Let $\{u_n\}\subset \mathscr{A}$ be a minimizing sequence, then $\{u_n\}$ is a bounded $(PS)_{c_0}$ sequence of $\varphi$. Note that $\displaystyle 0<c_0\leq c^*<\left(\frac{1}{2}-\frac{1}{2^*(s_1)}\right)\mu_{s_1}(\R_+^N)^{\frac{2^*(s_1)}{2^*(s_1)-2}}$, by \cite[Theorem 3.1]{CeramiZhongZou.2015.DOI10.1007/s00526-015-0844-z.} or Theorem \ref{2014-6-19-th1} again, up to a subsequence, $u_n\rightarrow u_0$, a ground state solution of \eqref{2015-4-24-we1}.
Without loss of generality, we may assume that $u_0\geq 0$.
Now, we let $a(x):=-\lambda^*\frac{1}{|x|^{\sigma_0}}+\frac{|u_0|^{2^*(s_1)-2}}{|x|^{s_1}}$, it is easy to check that $a(x)\in L_{loc}^{\frac{N}{2}}(\Omega)$ and $u_0$ is a weak solution of
$$-\Delta u=a(x)u, u\in H_0^1(\Omega).$$
Then the Br\'{e}zis-Kato theorem \cite{BrezisKato.1978} implies that $u_0\in L_{loc}^{r}(\Omega)$ for all $1\leq r<\infty$. Then $u_0\in W_{loc}^{2,r}(\Omega)$ for all $1\leq r<\infty$. By  the elliptic regularity theory, $u_0\in C^2(\Omega)$.
Finally, by the maximum principle, we obtain that $u_0$ is positive.
\hfill$\Box$

\subsubsection{The cases of  $p=1$ with $|\Omega|=\infty$ or $p>1$ with unbounded domain $\Omega$}
In the following, we are concerned with the existence of  positive solutions to the following problem
\be\lab{UBDP}
\begin{cases}
&-\Delta u+\lambda \frac{u^p}{|x|^{s_1}}=\frac{u^{2^*(s_2)-1}}{|x|^{s_2}}\;\quad \hbox{in}\;\Omega,\\
&u(x)>0\;\hbox{in}\;\Omega,\\
& u(x)=0\;\hbox{on}\;\partial\Omega,
\end{cases}
\ee
 where $\Omega$ is an open unbounded domain in $\R^N$. $N\geq 3, 2^*(s_2):=\frac{2(N-s_2)}{N-2}$.
When $p<2^*(s_1)-1$, there seems to have no results  when  $\Omega$ is an open unbounded domain
but $\Omega\neq \R_+^N$. Here we are aim to give something new. We obtain Theorem \ref{2014-6-20-wth1}.

\vskip 0.02in
When $|\Omega|=\infty$ and $1\leq p<2^*(s_1)-1$, we can not ensure that $\displaystyle \int_\Omega \frac{|u|^{p+1}}{|x|^{s_1}}dx$ is well defined for all $u\in D_{0}^{1,2}(\Omega)$. Hence, Problem (\ref{UBDP}) is posed in the framework of the Sobolev space
$$E:=D_{0}^{1,2}(\Omega)\cap L^{p+1}\left(\Omega, \frac{dx}{|x|^{s_1}}\right).$$  It is easy to see that $E$ is weakly  closed  due to the Fatou's Lemma. We take the norm
$$\|u\|_E:=\left(\int_{\Omega}|\nabla u|^2dx\right)^{\frac{1}{2}}+\left(\int_\Omega \frac{|u|^{p+1}}{|x|^{s_1}}dx\right)^{\frac{1}{p+1}}.$$
In this case, we always assume that $\lambda>0$.

\vskip0.1in

When $|\Omega|<\infty$, $E=D_{0}^{1,2}(\Omega)$, we always assume that $p>1$ since $p=1$ has been studied in subsection 3.2.1. To (\ref{UBDP}) there corresponds the variational functional, defined in $E$, by
$$\varphi(u):=\frac{1}{2}\int_\Omega |\nabla u|^2 dx+\frac{\lambda}{p+1}\int_\Omega \frac{|u|^{p+1}}{|x|^{s_1}}dx-\frac{1}{2^*(s_2)}\int_\Omega \frac{|u|^{2^*(s_2)}}{|x|^{s_2}}dx,$$
which is of class $C^2(E,\R)$.

\bl\lab{2014-6-18-xl1}
Assume that $\Omega\subset \R^N$ is an open unbounded domain and $$
\begin{cases}
\lambda>0, 1\leq p< \min\{2^*(s_2)-1, 2^*(s_1)-1\}\quad &\hbox{if}\;|\Omega|=\infty;\\\\
\hbox{either} \;\lambda< 0, 1<p< 2^*(s_1)-1\;\hbox{or}\\
\lambda> 0, 1<p< \min\{2^*(s_2)-1, 2^*(s_1)-1\}&\hbox{if}\;|\Omega|<\infty.
\end{cases}
$$
Then for any $u\in E\backslash \{0\}$, there exists a  unique $t_u>0$ such that
$$t_uu\in \mathcal{N}:=\{u\neq 0, J(u)=0\},$$
where
$$J(u):=\langle \varphi'(u), u\rangle =\int_\Omega |\nabla u|^2 dx+\lambda\int_\Omega \frac{|u|^{q}}{|x|^{s_1}}dx-\int_\Omega \frac{|u|^{2^*(s_2)}}{|x|^{s_2}}dx.$$
Moreover, $\mathcal{N}$ is close and bounded away from $0$.
\el
\bp
The proof of existence and uniqueness of $t_u$ is standard, we omit the details and  refer to \cite[Lemma 2.1]{ZhongZou.2015arXiv:1504.00730[math.AP]}. And it is trivial that $\mathcal{N}$ is closed. The proof of that $\mathcal{N}$ is bounded away from $0$ is a little different, hence we prefer to give the details of this part.

For any $u\in \mathcal{N}$, we have $J(u)=0$. Then by  the Hardy-Sobolev inequality, for the case of $\lambda>0$, we have
\begin{align*}
\int_\Omega |\nabla u|^2 dx<&\int_\Omega |\nabla u|^2 dx+\lambda\int_\Omega \frac{|u|^{p+1}}{|x|^{s_1}}dx
=\int_\Omega \frac{|u|^{2^*(s_2)}}{|x|^{s_2}}dx\\
\leq&\left(\mu_{s_2}(\R^N)^{-1}\int_\Omega |\nabla u|^2 dx\right)^{\frac{2^*(s_2)}{2}},
\end{align*}
and for the case of $\lambda<0, |\Omega|<\infty, p>1$, by Lemma \ref{2014-6-20-l1}, we also have
\begin{align*}
\int_\Omega |\nabla u|^2 dx=&\int_\Omega \frac{|u|^{2^*(s_2)}}{|x|^{s_2}}dx-\lambda\int_\Omega \frac{|u|^{p+1}}{|x|^{s_1}}dx\\
\leq&\left(\mu_{s_2}(\R^N)^{-1}\int_\Omega |\nabla u|^2 dx\right)^{\frac{2^*(s_2)}{2}}+C|\lambda| \left(\int_\Omega |\nabla u|^2 dx\right)^{\frac{p+1}{2}}.
\end{align*}
Hence, there exists some $\delta_0>0$ such that
$\displaystyle \big(\int_\Omega |\nabla u|^2 dx\big)^{\frac{1}{2}}\geq \delta_0$ and it follows that
\be\lab{2014-6-18-e4}
\|u\|_E>\left(\int_\Omega |\nabla u|^2 dx\right)^{\frac{1}{2}}\geq \delta_0>0.
\ee
Hence, $\mathcal{N}$ is bounded away from $0$.
\ep
\bl\lab{2014-6-18-xl2}
Assume that $\Omega\subset \R^N$ is an open unbounded domain and $$
\begin{cases}
\lambda>0, 1\leq p< \min\{2^*(s_2)-1, 2^*(s_1)-1\}\quad &\hbox{if}\;|\Omega|=\infty;\\\\
\hbox{either}\; \lambda< 0, 1<p< 2^*(s_1)-1\;\hbox{ or }\\
\lambda> 0, 1<p< \min\{2^*(s_2)-1, 2^*(s_1)-1\}&\hbox{if}\;|\Omega|<\infty
\end{cases}
$$
Then any $(PS)_c$ sequence of $\varphi$ is bounded in $E$ and any $\{u_n\}\subset \mathcal{N}$, a $(PS)_c$ sequence of $\varphi\big|_{\mathcal{N}}$, is also a $(PS)_c$ sequence of $\varphi$.
\el
\bp
We firstly consider the case of $|\Omega|=\infty$ and
$\lambda>0, 1\leq p< \min\{2^*(s_2)-1, 2^*(s_1)-1\}$.

Le $\{u_n\}\subset E$ be a $(PS)_c$ sequence of $\varphi(u)$, that is
$$
\begin{cases}
\varphi'(u_n)\rightarrow 0\;\hbox{in}\;E^*,\\
\varphi(u_n)\rightarrow c.
\end{cases}
$$
Denote
$\displaystyle a_n:=\int_\Omega |\nabla u_n|^2dx, b_n:=\int_\Omega \frac{|u|^q}{|x|^{s_1}}dx, c_n:=\int_\Omega \frac{|u_n|^{2^*(s_2)}}{|x|^{s_2}}, q=p+1$.
Then, we have
\be\lab{2014-6-21-we1}
\langle \varphi'(u_n), u_n\rangle=a_n+\lambda b_n-c_n=o(1)\|u_n\|_E
\ee
and
\be\lab{2014-6-21-we2}
\varphi(u_n)=\frac{1}{2}a_n+\frac{\lambda}{q}b_n-\frac{1}{2^*(s_2)}c_n=c+o(1).
\ee
Hence,
\be\lab{2014-6-21-we3}
\left(\frac{1}{2}-\frac{1}{2^*(s_2)}\right)a_n+\left(\frac{1}{q}-\frac{1}{2^*(s_2)}\right)\lambda b_n=c+o(1)\left(1+\|u_n\|_E\right).
\ee
Noting that $\|u_n\|_E=a_{n}^{\frac{1}{2}}+b_{n}^{\frac{1}{q}}$, without loss of generality, going to a subsequence if necessary, we may assume that
$$a_{n}^{\frac{1}{2}}\geq \frac{1}{2}\|u_n\|_E\;\hbox{or}\;b_{n}^{\frac{1}{q}}\geq \frac{1}{2}\|u_n\|_E.$$
If $a_{n}^{\frac{1}{2}}\geq \frac{1}{2}\|u_n\|_E$,  recalling that  $\lambda>0, q<2^*(s_2)$, we have
\begin{align*}
\left(\frac{1}{2}-\frac{1}{2^*(s_2)}\right)\frac{1}{4}\|u_n\|_E^2\leq&\left(\frac{1}{2}-\frac{1}{2^*(s_2)}\right)a_n\\
\leq&\left(\frac{1}{2}-\frac{1}{2^*(s_2)}\right)a_n+\left(\frac{1}{q}-\frac{1}{2^*(s_2)}\right)\lambda b_n\\
=&c+o(1)\left(1+\|u_n\|_E\right),
\end{align*}
and if $b_{n}^{\frac{1}{q}}\geq \frac{1}{2}\|u_n\|_E$, we have
\begin{align*}
\left(\frac{1}{q}-\frac{1}{2^*(s_2)}\right)\lambda \left(\frac{\|u_n\|_E}{2}\right)^q\leq&\left(\frac{1}{q}-\frac{1}{2^*(s_2)}\right)\lambda b_n\\
\leq&\left(\frac{1}{2}-\frac{1}{2^*(s_2)}\right)a_n+\left(\frac{1}{q}-\frac{1}{2^*(s_2)}\right)\lambda b_n\\
=&c+o(1)\left(1+\|u_n\|_E\right),
\end{align*}
which implies that $\{u_n\}$ is bounded in $E$.

\vskip0.2in

Secondly, we consider the case of $|\Omega|<\infty$. We note that if $\lambda>0, q<2^*(s_2)$ , the arguments displayed above are valid. Next, we assume that $\lambda<0, 2<q<2^*(s_1)$. Notice that in this case, the embedding $D_{0}^{1,2}(\Omega)\hookrightarrow L^q(\Omega,\frac{dx}{|x|^{s_1}})$ is compact and the norm $\|u_n\|_E$ is equivalent to the usual one $\|u_n\|=a_{n}^{\frac{1}{2}}$.
Hence,  note that
\begin{align*}
o(1)\|u_n\|_E+2c+o(1)=&\frac{2}{q}\lambda b_n-\frac{2}{2^*(s_2)}c_n-\lambda b_n+c_n\\
=&\left(1-\frac{2}{q}\right)|\lambda|b_n+\left(1-\frac{2}{2^*(s_2)}\right)c_n,
\end{align*}
we have
\begin{align*}
\|u_n\|_E^2\leq& C a_n=C\left(o(1)\|u_n\|_E+c_n+|\lambda|b_n\right)\\
\leq&\bar{C}\left(o(1)\|u_n\|_E+(1-\frac{2}{q})|\lambda|b_n+(1-\frac{1}{2^*(s_2)})c_n\right)\\
\leq&\tilde{C}\left(o(1)\|u_n\|_E+1\right),
\end{align*}
which also implies that $\{u_n\}$ is bounded in $E$.

\vskip0.1in
The remaining proof is similar to \cite[Lemma 2.3]{ZhongZou.2015arXiv:1504.00730[math.AP]}, we omit the details.
\ep

\bl\lab{2015-4-25-l1}
Let $\Omega\subset \R^N$ be a $C^1$ open domain with $0\in\partial\Omega$.
The boundary $\partial\Omega$ is $C^2$ at $0$ with the mean curvature $H(0)<0$. Assume that
$$
\begin{cases}
\lambda>0, 1\leq p< \min\{2^*(s_2)-1, 2^*(s_1)-1\}\quad &\hbox{if}\;|\Omega|=\infty;\\\\
\hbox{either}\; \lambda< 0, 1<p< 2^*(s_1)-1\;\hbox{ or }\\
\lambda> 0, 1<p< \min\{2^*(s_2)-1, 2^*(s_1)-1\}&\hbox{if}\;|\Omega|<\infty.
\end{cases}
$$
Furthermore, if $\lambda>0$, we require that $p<\frac{N-2s_1}{N-2}$ or $p\geq \frac{N-2s_1}{N-2}$ with $|\lambda|$ small enough. Then we have
$\displaystyle 0<c_0<\bar{\Theta}:=\left(\frac{1}{2}-\frac{1}{2^*(s_2)}\right)\mu_{s_2}(\R^N)^{\frac{2^*(s_2)}{2^*(s_2)-2}}$, where
\be c_0:=\inf_{u\in \mathcal{N}}\varphi(u).\ee
\el
\bp
Under the assumptions, by Lemma \ref{2014-6-18-xl1}, the corresponding Nehari manifold is well defined, hence $c_0$ is well defined. By \eqref{2014-6-18-e4}, we obtain that $c_0>0$. When $\Omega$ is bounded, $c_0<\bar{\Theta}$ is given by \cite[Lemma 2.4]{CeramiZhongZou.2015.DOI10.1007/s00526-015-0844-z.}. And we note that the proof is completely valid for the unbounded case, since the construction of testing function only required the local information of $\Omega$ near $x=0$.
\ep

\noindent{\bf Proof of Theorem \ref{2014-6-20-wth1}:}Let $\{u_n\}\subset \mathcal{N}$ be a minimizing sequence, then
$$\begin{cases}
\varphi(u_n)\rightarrow c_0:=\inf_{u\in \mathcal{N}}\varphi(u)\\
\varphi'(u_n)\rightarrow 0\;\hbox{in}\;E^*.
\end{cases}$$
By Lemma \ref{2015-4-25-l1}, we have $c_0\in (0,\bar{\Theta})$ and by Lemma \ref{2014-6-18-xl2}, $\{u_n\}$ is bounded in $E$.
Then by Theorem \ref{2014-6-19-th1}, up to a subsequence, there exists some $u\in E$ such that $u_n\rightarrow u$ strongly in $E$ and $\varphi(u)=c_0$.  Hence, $u\in \mathcal{N}$ is a minimizer. Without loss of generality, we may assume that $u\geq 0$  is a ground state solution of (\ref{UBDP}). Finally, it is standard to show $u_0$ is positive.\hfill$\Box$

\subsubsection{Li-Lin's Open problem on unbounded domain}
\noindent {\bf  Li-Lin's Open Problem} (see \cite[Remark 1.2]{LiLin.2012}) {\it  For the situation
$s_1<s_2$ and $\lambda<0$, the existence of positive solutions to
$$
\begin{cases}
&\Delta u+\lambda \frac{u^{2^*(s_1)-1}}{|x|^{s_1}}+\frac{u^{2^*(s_2)-1}}{|x|^{s_2}}=0\;\quad \hbox{in}\;\Omega,\\
&u(x)>0\;\hbox{in}\;\Omega,\quad u(x)=0\;\hbox{on}\;\partial\Omega.
\end{cases}
$$
 is completely open. Even for the equation
\be\lab{zou=1} \Delta u-u^p+\frac{u^{2^*(s)-1}}{|x|^s}=0 \;\hbox{in}\;    \Omega, \ee
where $0<s<2$ and $2^*(s)-1<p<\frac{N+2}{N-2}$, the existence problem still remains an interesting open question.} \hfill$\Box$

\vskip0.11in

The first partial answer to Li-Lin's open problem is obtained in   \cite{CeramiZhongZou.2015.DOI10.1007/s00526-015-0844-z.},
where the authors consider the following problem:
\be\lab{2014-6-24-e1}
\begin{cases}
&\Delta u+\lambda \frac{u^p}{|x|^{s_1}}+\frac{u^{2^*(s_2)-1}}{|x|^{s_2}}=0\;\quad \hbox{in}\;\Omega,\\
&u(x)>0\;\hbox{in}\;\Omega,\\
& u(x)=0\;\hbox{on}\;\partial\Omega,
\end{cases}
\ee
 where $\Omega$ is an open  bounded domain in $\R^N,\  N\geq 3, \ 2^*(s_2):=\frac{2(N-s_2)}{N-2}$.
 When $\ 2^*(s_2)-1<p\leq 2^*(s_1)-1 \ $, it is a case considered in the above mentioned Li-Lin's open problem. For such cases, we see that the functional still possesses the mountain pass geometric structure. We also note that the method of Nehari manifold fails now and one can not ensure any $(PS)_c$ sequence is bounded any more.  Hence, the Ekeland's principle is not enough for this problem, to overcome this difficulty, one thus needs  to develop more sophisticated critical point theories which assert the existence of a bounded $(PS)$ sequence.  Cerami-Zhong-Zou obtain the following results thanks to the perturbation method and the well-known Struwe's monotonicity trick (see \cite{JeanjeanToland.1998,Struwe.1988}).

\vskip 0.2in
\noindent
{\bf Theorem A1} {\it (cf. \cite[Theorem 1.5]{CeramiZhongZou.2015.DOI10.1007/s00526-015-0844-z.})
 Let $\Omega\subset \R^N$ be a $C^1$  bounded domain such  that $0 \in \partial\Omega. $   Assume  that  $ \partial\Omega$ is of  $C^2$ at $0$ and  that the mean curvature of $\partial\Omega$ at $0$, $H(0),$ is negative.  If  $$\lambda < 0,\quad  \quad   2^*(s_2)-1<p\leq 2^*(s_1)-1, $$  then there exists $\lambda_0<0$   such that (\ref{2014-6-24-e1})  has  a positive solution for all $\lambda_0<\lambda<0$.}
\hfill$\Box$

\vskip 0.2in
\noindent
{\bf Theorem A2} {\it (cf. \cite[Theorem 1.6]{CeramiZhongZou.2015.DOI10.1007/s00526-015-0844-z.})
Let $\Omega\subset \R^N$ be a $C^1$ bounded domain such  that $0 \in \partial\Omega.  $   Assume that $ \partial\Omega$ is of  $C^2$ at $0$ and   the mean curvature of $\partial\Omega$ at $0$, $H(0),$ is  negative.
If  $$\lambda < 0;   2^*(s_2)-1<p<\frac{N-2s_1}{N-2}, $$    then for almost every $\lambda<0$, (\ref{2014-6-24-e1})  has  a positive solution.}
\hfill$\Box$

Next, we are going to  study  the   Li-Lin's open problems on the unbounded domain $\Omega$. The framework of the Sobolev space $E=H_0^1(\Omega)$ if $|\Omega|<\infty$ and $E=D_{0}^{1,2}(\Omega)\cap L^{p+1}(\Omega, \frac{dx}{|x|^{s_1}})$ if $|\Omega|=\infty$.
Usually, the  nontrivial   weak limit $u$ of  a $(PS)$ sequence $u_n$   is a solution we  are searching  for.
The perturbation method state out the weak limit is nontrivial through the geometrical position without the splitting result.

\bt\lab{2014-6-24-th1}
 Let $\Omega\subset \R^N$ be a $C^1$  unbounded domain with $|\Omega|<\infty$ such  that $0 \in \partial\Omega. $   Suppose that $ \partial\Omega$  is of  $C^2$ at $0$ and  the mean curvature of $\partial\Omega$ at $0$, $H(0),$ is  negative.  Let assume $$\lambda < 0,   \hspace{0,4cm} 2^*(s_2)-1<p\leq 2^*(s_1)-1, $$  then there exists $\lambda_0<0$   such that (\ref{2014-6-24-e1})  has  a positive solution for all $\lambda_0<\lambda<0$.
\et
\bp
Replace $H_0^1(\Omega)$ by $E$ and the arguments of \cite[Theorem 1.5]{CeramiZhongZou.2015.DOI10.1007/s00526-015-0844-z.} are completely valid. We omit the details.
\ep

Recalling Lemma \ref{2015-4-25-l1} is still valid for unbounded domain. We see that following results obtain in \cite{CeramiZhongZou.2015.DOI10.1007/s00526-015-0844-z.} hold.
\vskip 0.2in
\noindent
{\bf Lemma A3} (cf.\cite[Lemma 6.1]{CeramiZhongZou.2015.DOI10.1007/s00526-015-0844-z.})
Consider the problem (\ref{2014-6-24-e1}) with $2^*(s_2)-1<p<\frac{N-2s_1}{N-2}$. Then
$\forall\;\lambda<0$, there exists $0<\delta<|\lambda|$ such that for almost every $\eta\in [\lambda-\delta, \lambda+\delta]$, there exists a bounded sequence $\{u_n\}\subset H_0^1(\Omega)\;\big(\{u_n\}\subset E\big)$ with
$$I_\eta(u_n)\rightarrow c_\eta<\Theta=:\big(\frac{1}{2}-\frac{1}{2^*(s_2)}\big)\mu_{s_2}(\R_+^N)^{\frac{2^*(s_2)}{2^*(s_2)-2}}, $$
$c_\eta > 0$ and $ I'_\eta(u_n)\rightarrow 0\;\hbox{in}\;H^{-1}(\Omega)\;\;\big(I'_\eta(u_n)\rightarrow 0\;\hbox{in}\;E^*\big)$ Where $$I_\lambda(u):=\frac{1}{2}\int_\Omega|\nabla u|^2 dx-\frac{\lambda}{p+1}\int_\Omega \frac{u_{+}^{p+1}}{|x|^{s_1}}-\frac{1}{2^*(s_2)}\int_\Omega \frac{u_{+}^{2^*(s_2)}}{|x|^{s_2}}dx$$
\hfill$\Box$

\vskip 0.2in
\noindent
{\bf Corollary A4}(cf.\cite[Corollary 6.1]{CeramiZhongZou.2015.DOI10.1007/s00526-015-0844-z.})
For almost every $\lambda<0$, there exists a bounded sequence $\{u_n\}\subset H_0^1(\Omega)\;\big(\{u_n\}\subset E\big)$ such that
$I_\lambda(u_n)\rightarrow c_\lambda<\Theta$
and
$I'_\lambda(u_n)\rightarrow 0\;\hbox{in}\;H^{-1}(\Omega)\;\;\big(I'_\eta(u_n)\rightarrow 0\;\hbox{in}\;E^*\big).$
\hfill$\Box$

Recall  that Theorem \ref{2014-6-19-th1}  gives  a   splitting result for the case of unbounded domain, we can also obtain the following theorem:
\bt\lab{2014-6-24-th2}
Let $\Omega\subset \R^N$ be a $C^1$ unbounded domain such  that $0 \in \partial\Omega.$   Suppose that  $ \partial\Omega$ is of $C^2$ at $0$ and    the mean curvature of $\partial\Omega$ at $0$, $H(0),$  is negative.  Assume further
$$\lambda < 0, 2^*(s_2)-1<p<\frac{N-2s_1}{N-2}, $$
 then for almost every $\lambda<0$, (\ref{2014-6-24-e1})  has  a positive solution.
\et
\bp
By  Lemma A3 and  Corollary A4,  for almost every $\lambda<0$, there exists a bounded  $(PS)_{c_\lambda}$  sequence with  $0 < c_\lambda<\Theta$.
Then, up to a subsequence, there exists some $U_\lambda\in E$ such that $u_n\rightharpoonup U_\lambda$. We only need to prove that $U_\lambda\neq 0$. Indeed, assume $U_\lambda= 0,$ since $p<2^*(s_1)-1$ we can apply Theorem \ref{2014-6-19-th1} obtaining $k\geq 1$ and

 $$c_\lambda= I_\lambda(u_n)+o(1)=I_\lambda(U_\lambda)+\sum_{i=1}^{k}I(U^k)+o(1)
= \sum_{i=1}^{k}I(U^k)+o(1)\geq\Theta,
$$
a contradiction.
Hence, a positive solution $U_\lambda$ to problem (\ref{2014-6-24-e1}) is found.
\ep

\br
Compare Theorem \ref{2014-6-24-th1}, Theorem \ref{2014-6-24-th2} with the ones for bounded domain case, i.e., Theorem A1 and Theorem A2.
What surprising us is that the representation  are exactly the same even for the range of $p$. Now, let us give an explanation. When we apply the perturbation argument or monotonicity trick, the framework of work space is $E$ and thus the bounded $(PS)_c$ sequence obtained are relate to the norm $\|\cdot\|_E$. Another related explanation see Remark \ref{2015-4-24-r1}.
\er

\s{Compactness results for the system case }
\renewcommand{\theequation}{4.\arabic{equation}}
\renewcommand{\theremark}{4.\arabic{remark}}
\renewcommand{\thedefinition}{4.\arabic{definition}}

\subsection{The case of $\alpha+\beta<2^*(s_2)$ and $|\Omega|<\infty$}
\bl\lab{2014-4-11-wl1}
Suppose that $\Omega\subset \R^N$ with the Lebesgue measure $|\Omega|<\infty$ and $0\leq s_2<2,   1<\alpha,1<\beta, \alpha+\beta<2^*(s_2)$. Let $\{(u_n,v_n)\}\subset \mathscr{D}:=H_0^1(\Omega)\times H_0^1(\Omega)$ be a bounded sequence such that $(u_n,v_n)\rightharpoonup (u,v)$ in $\mathscr{D}^*$, then
\be\lab{2014-4-11-we1}
\frac{|u_n|^{\alpha-2}u_n|v_n|^\beta}{|x|^{s_2}}\rightarrow \frac{|u|^{\alpha-2}u|v|^\beta}{|x|^{s_2}}\;\hbox{in}\;H^{-1}(\Omega)
\ee
and
\be\lab{2014-4-11-we2}
\frac{|u_n|^\alpha|v_n|^{\beta-2}v_n}{|x|^{s_2}}\rightarrow \frac{|u|^\alpha|v|^{\beta-2}v}{|x|^{s_2}}\;\hbox{in}\;H^{-1}(\Omega).
\ee
\el
\bp
We only prove (\ref{2014-4-11-we1}).
Since $|\Omega|<\infty$, there exists some $R>0$ such that $|\Omega\cap \mathbb{B}_{R}^c|<1$ and $\displaystyle \int_{\Omega\cap \mathbb{B}_{R}^c}\frac{1}{|x|^{s_2}}dx<1$ for $s_2\geq 0$.
Note that $\frac{1}{|x|^{s_2}}\in L^1(\Omega\cap \mathbb{B}_{R})$ since $\Omega\cap \mathbb{B}_{R}$ is bounded domain in $\R^N$ and $s_2<N$.  Hence, $\displaystyle  \frac{1}{|x|^{s_2}}\in L^1(\Omega)$. By the absolute continuity of the integral, we obtain that
\be\lab{2014-4-11-we3}
\int_{\Lambda\cap \Omega}\frac{1}{|x|^{s_2}}dx\rightarrow 0\;\hbox{as}\;meas(\Lambda)\rightarrow 0.
\ee
Without loss of generality, we may assume that $u_n\xrightarrow{a.e.} u, v_n\xrightarrow{a.e.} v$ in $\Omega$ up to a subsequence.
Thanks to the Egoroff Theorem,  for any $\varepsilon>0$, there exists some $\Omega_1\subset \Omega$ such that $| \Omega_1|<\varepsilon$ and
\be\lab{2014-4-11-we4}
u_n\rightarrow u, v_n\rightarrow v\;\hbox{uniformly in}\;\Omega\backslash \Omega_1\;\hbox{as}\;n\rightarrow \infty.
\ee
For $\forall\;h\in H_0^1(\Omega)$, we have
\begin{align*}
&\left|\int_\Omega \left(\frac{|u_n|^{\alpha-2}u_n|v_n|^\beta}{|x|^{s_2}}- \frac{|u|^{\alpha-2}u|v|^\beta}{|x|^{s_2}}\right)hdx\right|\\
& \leq \int_{\Omega_1} \left|\frac{|u_n|^{\alpha-2}u_n|v_n|^\beta}{|x|^{s_2}}h\right|dx+\int_{\Omega_1} \left|\frac{|u|^{\alpha-2}u|v|^\beta}{|x|^{s_2}}h\right|dx\\
&\quad +\left|\int_{\Omega\backslash \Omega_1} \left(\frac{|u_n|^{\alpha-2}u_n|v_n|^\beta}{|x|^{s_2}}- \frac{|u|^{\alpha-2}u|v|^\beta}{|x|^{s_2}}\right)hdx\right|\\
&:= I+II+III.
\end{align*}
Recalling that $(u_n,v_n)$ is bounded, by (\ref{2014-4-11-we3}), the  H\"{o}lder inequality and the Hardy-Sobolev inequality, we have
\begin{align*}
I=&
\int_{\Omega_1} \left|\frac{|u_n|^{\alpha-2}u_n|v_n|^\beta}{|x|^{s_2}}h\right|dx\\
=&\int_{\Omega_1}\frac{|h|}{|x|^{\frac{s_2}{2^*(s_2)}}}
\frac{|v_n|^\beta}{|x|^{\frac{\beta s_2}{2^*(s_2)}}}
\frac{|u_n|^{\alpha-1}}{|x|^{\frac{(\alpha-1)s_2}{2^*(s_2)}}}
\frac{1}{|x|^{\frac{(2^*(s_2)-\alpha-\beta)s_2}{2^*(s_2)}}}dx\\
\leq&\left(\int_{\Omega_1}\frac{|h|^{2^*(s_2)}}{|x|^{s_2}}\right)^{\frac{1}{2^*(s_2)}}
\left(\int_{\Omega_1}\frac{|v_n|^{2^*(s_2)}}{|x|^{s_2}}\right)^{\frac{\beta}{2^*(s_2)}}
\left(\int_{\Omega_1}\frac{|u_n|^{2^*(s_2)}}{|x|^{s_2}}\right)^{\frac{\alpha-1}{2^*(s_2)}}\\
&\left(\int_{\Omega_1}\frac{1}{|x|^{s_2}}\right)^{\frac{2^*(s_2)-\alpha-\beta}{2^*(s_2)}}\\
=&o(1)\|h\|\;\hbox{as}\;|\Omega_1|\rightarrow 0.
\end{align*}
Similarly, we also have $II=o(1)\|h\|$ as $|\Omega_1|\rightarrow 0$.
On the other hand, for any fixed $\Omega_1$, by (\ref{2014-4-11-we4}) and the Hardy-Sobolev inequality, it is easy to see that
\be\lab{2014-4-11-we5}
III=o(1)\|h\|\;\hbox{as}\;n\rightarrow \infty.
\ee
Hence, we obtain
$$
\left|\int_\Omega \left(\frac{|u_n|^{\alpha-2}u_n|v_n|^\beta}{|x|^{s_2}}- \frac{|u|^{\alpha-2}u|v|^\beta}{|x|^{s_2}}\right)hdx\right|=o(1)\|h\|,
$$
which means that
$$\frac{|u_n|^{\alpha-2}u_n|v_n|^\beta}{|x|^{s_2}}\rightarrow \frac{|u|^{\alpha-2}u|v|^\beta}{|x|^{s_2}}\;\hbox{in}\;H^{-1}(\Omega).$$
\ep

\subsection{The case of $\alpha+\beta<2^*(s_2)$ and $|\Omega|=\infty$}
For the case of $|\Omega|=\infty$, we always assume that $\lambda^*, \mu^*>0$, then $\|u\|_E, \|v\|_F$ defined by \eqref{normE} and \eqref{normF} are norms.
We denote
\begin{itemize}
\item[{\bf $(A_{\sigma_0}^{*})$}]  either $\sigma_0<s_2$ or $\begin{cases}0<s_2\leq \sigma_0,\\ 2+\frac{4}{\sigma_0}\frac{\sigma_0-s_2}{N-2}< \alpha+\beta\end{cases}$;\\
\item[{\bf $(A_{\eta_0}^{*})$}]  either $\eta_0<s_2$ or $\begin{cases} 0<s_2\leq \eta_0,\\ 2+\frac{4}{\eta_0}\frac{\eta_0-s_2}{N-2}< \alpha+\beta\end{cases}$.
\end{itemize}
\bl\lab{2014-6-6-zl1}
Suppose that $\Omega\subset \R^N$ with   $|\Omega|=\infty$ and $(A_{\sigma_0}^{*}), (A_{\eta_0}^{*})$, $1<\alpha,1<\beta, \alpha+\beta<2^*(s_2)$. Let $\{(u_n,v_n)\}\subset \mathscr{D}:E\times F$ be a bounded sequence such that $(u_n,v_n)\rightharpoonup (u,v)$ in $\mathscr{D}^*$, then
\be\lab{2014-6-6-ze1}
\frac{|u_n|^{\alpha-2}u_n|v_n|^\beta}{|x|^{s_2}}\rightarrow \frac{|u|^{\alpha-2}u|v|^\beta}{|x|^{s_2}}\;\hbox{in}\;H^{-1}(\Omega)
\ee
and
\be\lab{2014-6-6-ze2}
\frac{|u_n|^\alpha|v_n|^{\beta-2}v_n}{|x|^{s_2}}\rightarrow \frac{|u|^\alpha|v|^{\beta-2}v}{|x|^{s_2}}\;\hbox{in}\;H^{-1}(\Omega).
\ee
\el
\bp
We only prove (\ref{2014-6-6-ze1}). Without loss of generality, we assume that $u_n\rightarrow u$ a.e. in $\Omega$.
Under the assumptions $(A_{\sigma_0}^{*}), (A_{\eta_0}^{*})$,
by Proposition \ref{2014-6-5-xl1}, we see that $\frac{|u_n|^{\alpha+\beta}}{|x|^{s_2}}, \frac{|v_n|^{\alpha+\beta}}{|x|^{s_2}}$ are tight sequences.
Then for any $h\in \mathscr{D}$, by  the H\"{o}lder inequality, we have
\begin{align}\lab{2014-6-6-ze3}
&\left|\int_{\mathbb{B}_R^c\cap \Omega} \left(\frac{|u_n|^{\alpha-2}u_n|v_n|^\beta}{|x|^{s_2}}- \frac{|u|^{\alpha-2}u|v|^\beta}{|x|^{s_2}}\right)hdx\right|\nonumber\\
& \leq \int_{\mathbb{B}_R^c\cap \Omega} \left|\frac{|u_n|^{\alpha-2}u_n|v_n|^\beta}{|x|^{s_2}}h\right|dx+\int_{\mathbb{B}_R^c\cap \Omega} \left|\frac{|u|^{\alpha-2}u|v|^\beta}{|x|^{s_2}}h\right|dx\nonumber\\
& \leq \left(\int_{\mathbb{B}_R^c\cap \Omega}\frac{|u_n|^{\alpha+\beta}}{|x|^{s_2}}dx\right)^{\frac{\alpha-1}{\alpha+\beta}}
\left(\int_{\mathbb{B}_R^c\cap \Omega}\frac{|v_n|^{\alpha+\beta}}{|x|^{s_2}}dx\right)^{\frac{\beta}{\alpha+\beta}} \left(\int_{\mathbb{B}_R^c\cap \Omega}\frac{|h|^{\alpha+\beta}}{|x|^{s_2}}dx\right)^{\frac{1}{\alpha+\beta}}\nonumber\\
&\quad +\left(\int_{\mathbb{B}_R^c\cap \Omega}\frac{|u|^{\alpha+\beta}}{|x|^{s_2}}dx\right)^{\frac{\alpha-1}{\alpha+\beta}}
\left(\int_{\mathbb{B}_R^c\cap \Omega}\frac{|v|^{\alpha+\beta}}{|x|^{s_2}}dx\right)^{\frac{\beta}{\alpha+\beta}} \left(\int_{\mathbb{B}_R^c\cap \Omega}\frac{|h|^{\alpha+\beta}}{|x|^{s_2}}dx\right)^{\frac{1}{\alpha+\beta}}\nonumber\\
& = o(1)\|h\|\;\hbox{as}\;R\rightarrow \infty.
\end{align}
Apply the similar arguments as that in Lemma \ref{2014-4-11-wl1}, we can obtain that
\be\lab{2014-6-6-ze4}
\left|\int_{\{|x|\leq R\}\cap \Omega} \left(\frac{|u_n|^{\alpha-2}u_n|v_n|^\beta}{|x|^{s_2}}- \frac{|u|^{\alpha-2}u|v|^\beta}{|x|^{s_2}}\right)hdx\right|=o(1)\|h\|\;\hbox{as}\;n\rightarrow \infty.
\ee
By (\ref{2014-6-6-ze3}) and (\ref{2014-6-6-ze4}), we get (\ref{2014-6-6-ze1}).
\ep


\vskip0.3in

\subsection{The case of $\alpha+\beta=2^*(s_2)$}
\bl\lab{2014-6-2-yl1}
Let $\sigma\geq1$, then for any $\varepsilon>0$, there exists $C(\varepsilon)>0$ such that
\be\lab{2014-6-2-ye1}
\Big||x+y|^{\sigma-2}(x+y)-|x|^{\sigma-2}x\Big|\leq \varepsilon |x|^{\sigma-1}+C(\varepsilon)|y|^{\sigma-1}.
\ee
\el
\bp

It is trivial for $x=0$ or $\sigma=1$. Next, we always assume that $\sigma>1$ and $x\neq 0$.
For the case of $x>0$, let $t=\frac{y}{x}$, then (\ref{2014-6-2-ye1}) is equivalent to the following inequality
\be\lab{2014-6-2-ye2}
\Big||1+t|^{\sigma-2}(1+t)-1\Big|\leq \varepsilon+C(\varepsilon)|t|^{\sigma-1}.
\ee
We define $$f(t):=\Big||1+t|^{\sigma-2}(1+t)-1\Big|,\quad t\in \R,$$
then it is easy to see that $f(t)$ is continuous on $\R$. Notice that $f(0)=0$, we obtain that for any $\varepsilon>0$, there exists some $t_0>0$ such that
\be\lab{2014-6-2-ye3}
f(t)\leq \varepsilon\quad \hbox{if}\;|t|\leq t_0.
\ee
On the other hand, since
$$\limsup_{|t|\rightarrow \infty}\frac{f(t)}{|t|^{\sigma-1}}=1.$$
Thus, by the continuity, there exists some $C(\varepsilon)>0$ such that
\be\lab{2014-6-2-ye4}
f(t)\leq C(\varepsilon) |t|^{\sigma-1}, \quad |t|\geq t_0.
\ee
Then,  (\ref{2014-6-2-ye2}) follows from (\ref{2014-6-2-ye3}) and (\ref{2014-6-2-ye4}).

For the case of $x<0$, we let $t=-\frac{y}{x}$, then (\ref{2014-6-2-ye1}) is equivalent to the following inequality
\be\lab{2014-6-2-ye5}
\Big||t-1|^{\sigma-2}(t-1)+1\Big|\leq \varepsilon +C(\varepsilon)|t|^{\sigma-1}.
\ee
Let $$g(t):=\Big||t-1|^{\sigma-2}(t-1)+1\Big|, \quad t\in \R,$$
the similar arguments also lead to (\ref{2014-6-2-ye5}).
\ep
\bc\lab{2014-6-2-cro1}
Assume that $0\leq s<2, \sigma>1$. Let $\{u_n\}\subset L^{2^*(s)}\left(\R^N,   \frac{dx}{|x|^s}\right)$ be a bounded sequence such that $u_n\rightarrow u$ a.e in $\R^N$, then
\be\lab{2014-6-2-ye6}
\frac{|u_n|^{\sigma-2}u_n-|u_n-u|^{\sigma-2}(u_n-u)-|u|^{\sigma-2}u}{|x|^{\frac{(\sigma-1)s}{2^*(s)}}}\rightarrow 0\;\hbox{in}\;L^{\frac{2^*(s)}{\sigma-1}}(\R^N).
\ee
\ec
\bp
For any fixed $\varepsilon>0$, define
$$f_n^\varepsilon:=\Big(\big||u_n|^{\sigma-2}u_n- |u_n-u|^{\sigma-2}(u_n-u)-|u|^{\sigma-2}u\big|-\varepsilon |u_n-u|^{\sigma-1}\Big)^+,$$
then by Lemma  \ref{2014-6-2-yl1}, we have
\be\lab{2014-6-2-ye7}
f_n^\varepsilon \leq \big(1+C(\varepsilon)\big)|u|^{\sigma-1}.
\ee
Notice that $f_n^\varepsilon\rightarrow 0$ a.e in $\R^N$, then by Lebesgue's dominated convergence theorem,
\be\lab{2014-6-2-ye8}
\frac{f_n^\varepsilon}{|x|^{\frac{(\sigma-1)s}{2^*(s)}}}\rightarrow 0\;\hbox{in}\;L^{\frac{2^*(s)}{\sigma-1}}(\R^N).
\ee
By (\ref{2014-6-2-ye7}) again, we have
\be\lab{2014-6-2-ye9}
\left|\frac{|u_n|^{\sigma-2}u_n- |u_n-u|^{\sigma-2}(u_n-u)-|u|^{\sigma-2}u}{|x|^{\frac{(\sigma-1)s}{2^*(s)}}}\right|\leq \frac{f_n^\varepsilon}{|x|^{\frac{(\sigma-1)s}{2^*(s)}}}+\varepsilon \frac{|u_n-u|^{\sigma-1}}{|x|^{\frac{(\sigma-1)s}{2^*(s)}}}.
\ee
Hence,
\begin{align*}
&\limsup_{n\rightarrow \infty}\int_{\R^N}\left|\frac{|u_n|^{\sigma-2}u_n- |u_n-u|^{\sigma-2}(u_n-u)-|u|^{\sigma-2}u}{|x|^{\frac{(\sigma-1)s}{2^*(s)}}}\right|^{\frac{2^*(s)}{\sigma-1}}\\
& \leq \limsup_{n\rightarrow \infty}\int_{\R^N}\left(\frac{f_n^\varepsilon}{|x|^{\frac{(\sigma-1)s}{2^*(s)}}}+\varepsilon \frac{|u_n-u|^{\sigma-1}}{|x|^{\frac{(\sigma-1)s}{2^*(s)}}}\right)^{\frac{2^*(s)}{\sigma-1}}\\
& \leq 2^{\frac{2^*(s)}{\sigma-1}}\left[\limsup_{n\rightarrow \infty}\int_{\R^N}\left(\frac{f_n^\varepsilon}{|x|^{\frac{(\sigma-1)s}{2^*(s)}}} \right)^{\frac{2^*(s)}{(\sigma-1)}}+\left(\varepsilon \frac{|u_n-u|^{\sigma-1}}{|x|^{\frac{(\sigma-1)s}{2^*(s)}}}\right)^{\frac{2^*(s)}{\sigma-1}}\right]\\
&=2^{\frac{2^*(s)}{\sigma-1}}\limsup_{n\rightarrow \infty}\int_{\R^N}\left(\varepsilon \frac{|u_n-u|^{\sigma-1}}{|x|^{\frac{(\sigma-1)s}{2^*(s)}}}\right)^{\frac{2^*(s)}{\sigma-1}}\\
& \leq C\varepsilon^{\frac{2^*(s)}{\sigma-1}}.
\end{align*}
Now, let $\varepsilon\rightarrow 0$, we obtain (\ref{2014-6-2-ye6}).
\ep
\bl\lab{2014-6-3-xl1}
Assume that $0\leq s_2<2, 1<\alpha, 1<\beta, \alpha+\beta=2^*(s_2)$.
Let $\displaystyle \{(u_n,v_n)\}\subset L^{2^*(s_2)}\left(\R^N,   \frac{dx}{|x|^{s_2}}\right)\times L^{2^*(s_2)}\left(\R^N, \frac{dx}{|x|^{s_2}}\right)$ be a bounded sequence and
$u_n\xrightarrow{a.e.} u, v_n\xrightarrow{a.e.} v$ in $\R^N$, then we have
\be\lab{2014-6-2-e4}
\frac{|u_n|^\alpha-|u_n-u|^\alpha-|u|^\alpha}{|x|^{\frac{\alpha s_2}{2^*(s_2)}}}\rightarrow 0\;\hbox{in}\;L^{\frac{2^*(s_2)}{\alpha}}(\R^N)
\ee
and
\be\lab{2014-6-2-e5}
\frac{|v_n|^\beta-|v_n-v|^\beta-|v|^\beta}{|x|^{\frac{\beta s_2}{2^*(s_2)}}}\rightarrow 0\;\hbox{in}\;L^{\frac{2^*(s_2)}{\beta}}(\R^N).
\ee
\el
\bp
For any fixed $\varepsilon>0$,  we define
$$f_n^\varepsilon:=\Big(\big||u_n|^\alpha-|u_n-u|^\alpha-|u|^\alpha\big|-\varepsilon|u_n-u|^\alpha\Big)^+\leq \big(1+c(\varepsilon)\big)|u|^\alpha$$
and
$$g_n^\varepsilon:=\Big(\big||v_n|^\beta-|v_n-v|^\beta-|v|^\beta\big|-\varepsilon|v_n-v|^\beta\Big)^+\leq \big(1+c(\varepsilon)\big)|v|^\beta,$$
where we  use  the facts of $\alpha>1,\beta>1$ and
$$\Big||x+y|^t-|x|^t\Big|\leq \varepsilon |x|^t+c(\varepsilon)|y|^t, 0< t<\infty.$$
Then by Lebesgue's dominated convergence theorem,
$$\frac{f_n^\varepsilon}{|x|^{\frac{\alpha s_2}{2^*(s_2)}}}\rightarrow 0\;\hbox{in}\;L^{\frac{2^*(s_2)}{\alpha}}(\R^N)\;\hbox{and}\; \frac{g_n^\varepsilon}{|x|^{\frac{\beta s_2}{2^*(s_2)}}}\rightarrow 0\;\hbox{in}\;L^{\frac{2^*(s_2)}{\beta}}(\R^N).$$
Since
$$\big||u_n|^\alpha-|u_n-u|^\alpha-|u|^\alpha\big|\leq f_n^\varepsilon +\varepsilon|u_n-u|^\alpha$$
and
$$\big||v_n|^\beta-|v_n-v|^\beta-|v|^\beta\big|\leq g_n^\varepsilon+\varepsilon|v_n-v|^\beta,$$
we obtain that
\begin{align}\lab{2014-6-2-e2}
&\limsup_{n\rightarrow \infty}\int_{\R^N}\left|\frac{|u_n|^\alpha-|u_n-u|^\alpha-|u|^\alpha}{|x|^{\frac{\alpha s_2}{2^*(s_2)}}}\right|^{\frac{2^*(s_2)}{\alpha}}\nonumber\\
& \leq \limsup_{n\rightarrow \infty}\int_{\R^N} \left|\frac{f_n^\varepsilon +\varepsilon|u_n-u|^\alpha}{|x|^{\frac{\alpha s_2}{2^*(s_2)}}}\right|^{\frac{2^*(s_2)}{\alpha}}\nonumber\\
& \leq  2^{\frac{2^*(s_2)}{\alpha}}\left[\limsup_{n\rightarrow\infty}\int_{\R^N} \left|\frac{f_n^\varepsilon}{|x|^{\frac{\alpha s_2}{2^*(s_2)}}}\right|^\frac{2^*(s_2)}{\alpha}+\varepsilon^{\frac{2^*(s_2)}{\alpha}}\frac{|u_n-u|^{2^*(s_2)}}{|x|^{s_2}}\right]\nonumber\\
& \leq  C\varepsilon.
\end{align}
Similarly,
\be\lab{2014-6-2-e3}
\limsup_{n\rightarrow \infty}\int_{\R^N}\left|\frac{|v_n|^\beta-|v_n-v|^\beta-|v|^\beta}{|x|^{\frac{\beta s_2}{2^*(s_2)}}}\right|^{\frac{2^*(s_2)}{\beta}}\leq C\varepsilon,
\ee
where $C$ is independent of $n$.
Let $\varepsilon\rightarrow 0$, we obtain (\ref{2014-6-2-e4}) and (\ref{2014-6-2-e5}).
\ep
\bl\lab{2014-6-3-l1}
Assume that $0\leq s_2<2, 1<\alpha, 1<\beta, \alpha+\beta=2^*(s_2)$.
Let $\{(u_n,v_n)\}\subset L^{2^*(s_2)}(\R^N, \frac{dx}{|x|^{s_2}})\times L^{2^*(s_2)}(\R^N, \frac{dx}{|x|^{s_2}})$ be a bounded sequence and
$u_n\xrightarrow{a.e.} u, v_n\xrightarrow{a.e.} v$ in $\R^N$, then
\be\lab{2014-6-3-e1}
\frac{|u_n|^{\alpha-2}u_n|v_n|^\beta}{|x|^{s_2}}-\frac{|u_n-u|^{\alpha-2}(u_n-u)|v_n-v|^\beta}{|x|^{s_2}}-\frac{|u|^{\alpha-2}u|v|^\beta}{|x|^{s_2}}\rightarrow 0\;\hbox{in}\;H^{-1}(\R^N)
\ee
and
\be\lab{2014-6-3-e2}
\frac{|u_n|^{\alpha}|v_n|^{\beta-2}v_n}{|x|^{s_2}}-\frac{|u_n-u|^{\alpha}|v_n-v|^{\beta-2}(v_n-v)}{|x|^{s_2}}-\frac{|u|^{\alpha}|v|^{\beta-2}v}{|x|^{s_2}}\rightarrow 0\;\hbox{in}\;H^{-1}(\R^N)
\ee
\el
\bp
We only prove (\ref{2014-6-3-e1}).
Let $h\in D_{0}^{1,2}(\R^N)$, it is sufficient to prove that
\be\lab{2014-6-3-e3}
\int_{\R^N}\left(\frac{|u_n|^{\alpha-2}u_n|v_n|^\beta}{|x|^{s_2}}-\frac{|u_n-u|^{\alpha-2}(u_n-u)|v_n-v|^\beta}{|x|^{s_2}}-\frac{|u|^{\alpha-2}u|v|^\beta}{|x|^{s_2}}\right)h=o(1)\|h\|.
\ee
We note that for any subset $\Omega$ of $\R^N$,
\begin{align*}
&\int_\Omega \frac{\Big(|u_n|^{\alpha-2}u_n|v_n|^\beta-|u_n-u|^{\alpha-2}(u_n-u)|v_n-v|^\beta-|u|^{\alpha-2}u|v|^\beta\Big)h}{|x|^{s_2}}\\
& = \int_\Omega \frac{\big[|u_n|^{\alpha-2}u_n-|u_n-u|^{\alpha-2}(u_n-u)-|u|^{\alpha-2}u\big]|v_n|^\beta h}{|x|^{s_2}}\\
&\quad +\int_\Omega \frac{|u_n-u|^{\alpha-2}(u_n-u)\big[|v_n|^\beta-|v_n-v|^\beta-|v|^\beta\big]h}{|x|^{s_2}}\\
&\quad +\int_\Omega \frac{|u|^{\alpha-2}u\big[|v_n|^\beta-|v_n-v|^\beta-|v|^\beta\big]h}{|x|^{s_2}}\\
&\quad +\int_\Omega \frac{|u_n-u|^{\alpha-2}(u_n-u)|v|^\beta h}{|x|^{s_2}}\\
&\quad +\int_\Omega \frac{|u|^{\alpha-2}u|v_n-v|^\beta h}{|x|^{s_2}}\\
& :=I(\Omega)+II(\Omega)+III(\Omega)+IV(\Omega)+V(\Omega).
\end{align*}
Take $\Omega=\R^N$, since $v_n$ is bounded in $\displaystyle L^{2^*(s_2)}\left(\R^N, \frac{dx}{|x|^{s_2}}\right)$ and $\frac{\alpha-1}{2^*(s_2)}+\frac{\beta}{2^*(s_2)}+\frac{1}{2^*(s_2)}=1$, by  the H\"{o}lder inequality, the Hardy-Sobolev inequality and Corollary \ref{2014-6-2-cro1}, we can have
\be\lab{2014-6-3-e4}
I(\R^N)=\int_{\R^N} \frac{\big[|u_n|^{\alpha-2}u_n-|u_n-u|^{\alpha-2}(u_n-u)-|u|^{\alpha-2}u\big]|v_n|^\beta h}{|x|^{s_2}}=o(1)\|h\|.
\ee
Similarly, by (\ref{2014-6-2-e5}), we can obtain that
\be\lab{2014-6-3-e5}
II(\R^N)=o(1)\|h\|, III(\R^N)=o(1)\|h\|.
\ee
Hence, we only need to prove that
\be\lab{2014-6-3-e6}
IV(\R^N)=o(1)\|h\|\;\hbox{and}\;V(\R^N)=o(1)\|h\|.
\ee
By the H\"{o}lder inequality and  the Hardy-Sobolev inequality again,  we have
\begin{align}\lab{2014-6-3-e7}
&\left|IV(\mathbb{B}_R^c)\right|=\left|\int_{\mathbb{B}_R^c}\frac{|u_n-u|^{\alpha-2}(u_n-u)|v|^\beta h}{|x|^{s_2}}\right|\nonumber\\
\leq&\left(\int_{\mathbb{B}_R^c}\frac{|u_n-u|^{2^*(s_2)}}{|x|^{s_2}}\right)^{\frac{\alpha-1}{2^*(s_2)}}
\left(\int_{\mathbb{B}_R^c}\frac{|h|^{2^*(s_2)}}{|x|^{s_2}}\right)^{\frac{1}{2^*(s_2)}}
\left(\int_{\mathbb{B}_R^c}\frac{|v|^{2^*(s_2)}}{|x|^{s_2}}\right)^{\frac{\beta}{2^*(s_2)}}\nonumber\\
\leq&C\|h\|\left(\int_{\mathbb{B}_R^c}\frac{|v|^{2^*(s_2)}}{|x|^{s_2}}\right)^{\frac{\beta}{2^*(s_2)}}.
\end{align}
Similarly, we can prove that
\be\lab{2014-6-3-e8}
\big|V(\mathbb{B}_R^c)\big|\leq C\|h\|\left(\int_{\mathbb{B}_R^c}\frac{|u|^{2^*(s_2)}}{|x|^{s_2}}\right)^{\frac{\alpha-1}{2^*(s_2)}}.
\ee
Hence, by the absolute continuity of the integral, we have
\be\lab{2014-6-3-e9}
IV(\mathbb{B}_R^c)=o(1)\|h\|, V(\mathbb{B}_R^c)=o(1)\|h\|\;\hbox{as}\;R\rightarrow \infty.
\ee
For some fixed $R$ large enough, apply the Egoroff Theorem in the bounded domain $\mathbb{B}_R$, for any $\delta>0$, there exists $\Omega_1\subset \mathbb{B}_R$ such that $|\Omega|<\delta$ and $u_n\rightarrow u, v_n\rightarrow v$ uniformly in $\mathbb{B}_R\backslash \Omega_1$.
The similar arguments above shows that
\be\lab{2014-6-3-e10}
IV(\Omega_1)=o(1)\|h\|, V(\Omega_1)=o(1)\|h\|\;\hbox{as}\;|\Omega_1|\rightarrow 0.
\ee
On the other hand, the uniform convergence leads to that
\be\lab{2014-6-3-e11}
IV(\mathbb{B}_R\backslash \Omega_1)=o(1)\|h\|, V(\mathbb{B}_R\backslash \Omega_1)=o(1)\|h\|\;\hbox{as}\;n\rightarrow \infty.
\ee
Then (\ref{2014-6-3-e6}) follows from (\ref{2014-6-3-e9})-(\ref{2014-6-3-e11}).
\ep

\bl\lab{2014-6-1-wl1}
Assume that  $ 0\leq s_2<2, 1<\alpha, 1<\beta, \alpha+\beta=2^*(s_2)$. Let $\{(u_n,v_n)\}$ be a bounded sequence of $\displaystyle L^{2^*(s_2)}\left(\R^N, \frac{dx}{|x|^{s_2}}\right)\times L^{2^*(s_2)}\left(\R^N, \frac{dx}{|x|^{s_2}}\right)$ such that $(u_n,v_n)\xrightarrow{a.e.} (u,v)$  in $\R^N$. Then
\be\lab{2014-6-3-xe0}
\lim_{n\rightarrow \infty}\int_{\R^N}\left( \frac{|u_n|^\alpha|v_n|^\beta}{|x|^{s_2}}-\frac{|u_n-u|^\alpha|v_n-v|^\beta}{|x|^{s_2}}\right)dx=\int_{\R^N} \frac{|u|^\alpha|v|^\beta}{|x|^{s_2}}dx.
\ee
\el
\bp
Since $v_n$ is bounded in $\displaystyle L^{2^*(s_2)}\left(\R^N, \frac{dx}{|x|^{s_2}}\right)$, by the H\"{o}lder inequality and \eqref{2014-6-2-e4}, we have
\be\lab{2014-6-3-xe1}
\int_{\R^N}\frac{\left[|u_n|^\alpha-|u_n-u|^\alpha-|u|^\alpha\right]|v_n|^\beta}{|x|^{s_2}}=o(1)\;\hbox{as}\;n\rightarrow \infty.
\ee
Similarly, since by (\ref{2014-6-2-e5}), we obtain that
\be\lab{2014-6-3-xe2}
\int_{\R^N}\frac{|u_n-u|^\alpha\big[|v_n|^\beta-|v_n-v|^\beta-|v|^\beta\big]}{|x|^{s_2}}=o(1)\;\hbox{as}\;n\rightarrow \infty
\ee
and
\be\lab{2014-6-3-xe3}
\int_{\R^N}\frac{|u|^\alpha\big[|v_n|^\beta-|v_n-v|^\beta-|v|^\beta\big]}{|x|^{s_2}}=o(1)\;\hbox{as}\;n\rightarrow \infty.
\ee
By the  H\"{o}lder inequality again and the absolute continuity of the integral, we have
\be\lab{2014-6-3-xe4}
\int_{\mathbb{B}_R^c}\frac{|u_n-u|^\alpha |v|^\beta}{|x|^{s_2}}=o(1), \int_{\mathbb{B}_R^c}\frac{|u|^\alpha |v_n-v|^\beta}{|x|^{s_2}}=o(1)\;\hbox{as}\;R\rightarrow \infty.
\ee
Apply the Egoroff Theorem in the bounded domain $\mathbb{B}_R$, for any $\delta>0$, there exists some $\Omega_1\subset \mathbb{B}_R$ such that
$|\Omega_1|<\delta$ and $u_n\rightarrow u, v_n\rightarrow v$ uniformly in $\mathbb{B}_R\backslash \Omega_1$. The similar arguments above shows that
\be\lab{2014-6-3-xe5}
\int_{\Omega_1}\frac{|u_n-u|^\alpha |v|^\beta}{|x|^{s_2}}=o(1), \int_{\Omega_1}\frac{|u|^\alpha |v_n-v|^\beta}{|x|^{s_2}}=o(1)\;\hbox{as}\;|\Omega_1|\rightarrow \infty.
\ee
The uniformly convergence leads to that
\be\lab{2014-6-3-xe6}
\int_{\mathbb{B}_R\backslash\Omega_1}\frac{|u_n-u|^\alpha |v|^\beta}{|x|^{s_2}}=o(1), \int_{\mathbb{B}_R\backslash\Omega_1}\frac{|u|^\alpha |v_n-v|^\beta}{|x|^{s_2}}=o(1)\;\hbox{as}\;n\rightarrow \infty.
\ee
It follows  from (\ref{2014-6-3-xe4})-(\ref{2014-6-3-xe6}) that
\be\lab{2014-6-3-xe7}
\int_{\R^N}\frac{|u_n-u|^\alpha |v|^\beta}{|x|^{s_2}}=o(1), \int_{\R^N}\frac{|u|^\alpha |v_n-v|^\beta}{|x|^{s_2}}=o(1)\;\hbox{as}\;n\rightarrow \infty.
\ee
Then by (\ref{2014-6-3-xe1})-(\ref{2014-6-3-xe3}) and (\ref{2014-6-3-xe7}), we obtain
(\ref{2014-6-3-xe0}). The proof is complete.
\ep


\bc\lab{2014-6-3-xcro1}
Assume that $0\leq s_2<2,   1<\alpha,1<\beta, \alpha+\beta=2^*(s_2)$. Let $\displaystyle \{(u_n,v_n)\}\subset L^{2^*(s_2)}\left(\R^N, \frac{dx}{|x|^{s_2}}\right)\times L^{2^*(s_2)}\left(\R^N, \frac{dx}{|x|^{s_2}}\right)$ be a bounded sequence such that $(u_n,v_n)\xrightarrow{a.e.} (u,v)$  in $\R^N$.
Furthermore, if $u_n\rightarrow u$ or $v_n\rightarrow v$ strongly in $\displaystyle L^{2^*(s_2)}\left(\R^N, \frac{dx}{|x|^{s_2}}\right)$, then
\be\lab{2014-6-3-we1}
\frac{|u_n|^{\alpha-2}u_n|v_n|^\beta}{|x|^{s_2}}\rightarrow \frac{|u|^{\alpha-2}u|v|^\beta}{|x|^{s_2}}\;\hbox{in}\;H^{-1}(\R^N)
\ee
and
\be\lab{2014-6-3-we2}
\frac{|u_n|^\alpha|v_n|^{\beta-2}v_n}{|x|^{s_2}}\rightarrow \frac{|u|^\alpha|v|^{\beta-2}v}{|x|^{s_2}}\;\hbox{in}\;H^{-1}(\R^N).
\ee
\ec
\bp
They are straightforward  conclusions from Lemma \ref{2014-6-3-l1}.
\ep

\bc\lab{2014-4-11-Cro1}
Suppose that $\Omega\subset \R^N$ is an open set and $0\leq s_2<2,   1<\alpha,1<\beta, \alpha+\beta\leq2^*(s_2)$. Let $\{(u_n,v_n)\}\subset \mathscr{D}$ be a bounded sequence such that $(u_n,v_n)\xrightarrow{a.e.} (u,v)$ in $\Omega$. Furthermore, for the case of $\alpha+\beta<2^*(s_2)$, we assume that  either $|\Omega|<\infty$ or $|\Omega|=\infty$ with $(A_{\sigma_0}^{*}), (A_{\eta_0}^{*})$.   Furthermore,  for the case of $\alpha+\beta=2^*(s_2)$, we require  that $u_n\rightarrow u$ or $v_n\rightarrow v$ strongly in $\displaystyle L^{2^*(s_2)}\left(\Omega, \frac{dx}{|x|^{s_2}}\right)$.
Then if $u=0$ or $v=0$, we have
\be\lab{2014-4-11-we6}
\frac{|u_n|^{\alpha-2}u_n|v_n|^\beta}{|x|^{s_2}}\rightarrow 0\;\hbox{in}\;H^{-1}(\Omega)
\ee
and
\be\lab{2014-4-11-we7}
\frac{|u_n|^\alpha|v_n|^{\beta-2}v_n}{|x|^{s_2}}\rightarrow 0\;\hbox{in}\;H^{-1}(\Omega).
\ee
\ec
\bp
It is a direct conclusion from Lemma \ref{2014-4-11-wl1}, Lemma \ref{2014-6-6-zl1} or Corollary \ref{2014-6-3-xcro1}.
\ep


\subsection{Splitting results}
Next, we will establish  the  splitting  theorem for the system case.
Assume that $\Omega\subset \R^N$ is an open domain with $0\in \partial\Omega$ such that $\partial\Omega$ is $C^2$ at $0$ and the mean curvature of $\partial\Omega$ at $0$ is negative. $0\leq s_2<2, 1<\alpha,1<\beta, \alpha+\beta\leq 2^*(s_2)$ and one of the following holds:
\begin{itemize}
\item[$(i)$]$\Omega$ is bounded, $0\neq \lambda^*>-\lambda_{1,\sigma_0}(\Omega), 0\neq \mu^*>-\lambda_{1,\eta_0}(\Omega)$;
\item[$(ii)$]$\Omega$ is unbounded but $|\Omega|<\infty$, $0\neq \lambda^*>-\lambda_{1,\sigma_0}(\Omega), 0\neq \mu^*>-\lambda_{1,\eta_0}(\Omega)$;
\item[$(iii)$]$|\Omega|=\infty$, $\lambda^*>0,\mu^*>0$ and conditions $(A_{\sigma_0}^*)$ and $(A_{\eta_0}^{*})$ are satisfied.
\end{itemize}

\br\lab{2015-bur1}
Since the domain we considered can be bounded or unbounded, and for the unbounded case, the measure is also allowed to be finite or infinite. The conditions we required will change a little in different situations. For example, when $|\Omega|<\infty$, we consider $\lambda^*>-\lambda_{1,\sigma_0}(\Omega),  \mu^*>-\lambda_{1,\eta_0}(\Omega)$. While $|\Omega|=\infty$, we only consider $\lambda^*>0,\mu^*>0$, moreover, we require the assumptions $(A_{\sigma_0}^*)$ and $(A_{\eta_0}^{*})$ if $\alpha+\beta<2^*(s_2)$. For the convenience, we prefer to denote our work space by the same denotation $\mathscr{D}:=E\times F$. When $|\Omega|<\infty$, it is easy to see that $E=F=H_0^1(\Omega)$.
\er

We denote the least energy corresponding to \eqref{2015-4-27-le1} by $c_{\sigma_0,\lambda^*,\lambda}$, then
\begin{align*}
c_{\sigma_0,\lambda^*,\lambda}=&\big(\frac{1}{2}-\frac{1}{2^*(s_1)}\big)\int_\Omega|\nabla u|^2+\lambda^*\frac{|u|^2}{|x|^{\sigma_0}}dx\\
=&\lambda^{\frac{-2}{2^*(s_1)-2}}\left(\frac{1}{2}-\frac{1}{2^*(s_1)}\right)\int_\Omega|\nabla \tilde{u}|^2+\lambda^*\frac{|\tilde{u}^2|}{|x|^{\sigma_0}}dx\\
=&\lambda^{\frac{-2}{2^*(s_1)-2}}c_{\sigma_0,\lambda^*, 1}.
\end{align*}
Hence, $c_{\sigma_0, \lambda^*,\lambda}$ is continuous and strictly decreasing relate to $\lambda\in (0,+\infty)$.  By  \cite[Theorem 4.1]{CeramiZhongZou.2015.DOI10.1007/s00526-015-0844-z.}, any ground state solution of \eqref{2015-4-27-le1} is a mountain pass solution.
\br\lab{2015-4-25-xr1}
When $\lambda^*>0$, the assumption $H(0)<0$ can guarantee
$$ \left(\frac{1}{2}-\frac{1}{2^*(s_1)}\right) \mu_{s_1}(\Omega)^{\frac{N-s_1}{2-s_1}}<c_{\sigma_0,\lambda^*,1}<\left(\frac{1}{2}-\frac{1}{2^*(s_1)}\right) \mu_{s_1}(\R_+^N)^{\frac{N-s_1}{2-s_1}}.$$ And when $|\Omega|<\infty$, our assumption is $\lambda^*>-\lambda_{1,\sigma_0}(\Omega)$, i.e., $\lambda^*$ is allowed to be negative. When $\lambda^*<0$, it is easy to prove that $\displaystyle c_{\sigma_0,\lambda^*,1}< \left(\frac{1}{2}-\frac{1}{2^*(s_1)}\right) \mu_{s_1}(\Omega)^{\frac{N-s_1}{2-s_1}}$. The details of this kind of estimation we refer to \cite{CeramiZhongZou.2015.DOI10.1007/s00526-015-0844-z.}.
\er

For any $u\in E, v\in F$, we define
$$\Psi_\lambda(u):=\frac{1}{2}\int_\Omega \left(|\nabla u|^2+\lambda^*\frac{|u|^2}{|x|^{\sigma_0}}\right)dx-\frac{\lambda}{2^*(s_1)}\int_\Omega \frac{|u|^{2^*(s_1)}}{|x|^{s_1}}dx$$
and
$$\Upsilon_\mu(v):=\frac{1}{2}\int_\Omega \left(|\nabla v|^2+\mu^*\frac{|v|^2}{|x|^{\eta_0}}\right)dx-\frac{\mu}{2^*(s_1)}\int_\Omega \frac{|v|^{2^*(s_1)}}{|x|^{s_1}}dx.$$
For any $(u,v)\in \mathscr{D}$,
\be\lab{2014-4-11-we12}
\Phi(u,v):=\Psi_{\lambda}(u)+\Upsilon_{\mu}(v)-\kappa\int_\Omega \frac{|u|^\alpha|v|^\beta}{|x|^{s_2}}dx.
\ee

\bl\lab{2014-4-11-wl2}
Suppose that $\Omega\subset \R^N$ is an open domain with $0\in \partial\Omega$ such that $\partial\Omega$ is $C^2$ at $0$ and the mean curvature of $\partial\Omega$ at $0$ is negative. $0\leq s_2<2,   1<\alpha,1<\beta, \alpha+\beta\leq2^*(s_2)$ and one of the following holds:
 \begin{itemize}
\item[$(1)$]$\Omega$ is bounded, $0\neq \lambda^*>-\lambda_{1,\sigma_0}(\Omega), 0\neq \mu^*>-\lambda_{1,\eta_0}(\Omega)$;
\item[$(2)$]$\Omega$ is unbounded but $|\Omega|<\infty$, $0\neq \lambda^*>-\lambda_{1,\sigma_0}(\Omega), 0\neq \mu^*>-\lambda_{1,\eta_0}(\Omega)$;
\item[$(3)$]$|\Omega|=\infty$, $\lambda^*>0,\mu^*>0$, especially when $\alpha+\beta<2^*(s_2)$ we require furthermore $(A_{\sigma_0}^*)$ and $(A_{\eta_0}^{*})$.
\end{itemize}

 Let $\{(u_n,v_n)\}\subset \mathscr{D}$ be a bounded $(PS)_m$ sequence such that $(u_n,v_n)\rightharpoonup (u,v)$ in $\mathscr{D}^*$. Furthermore, if $\alpha+\beta=2^*(s_2)$, we assume  further that $u_n\rightarrow u$ or $v_n\rightarrow v$ strongly in $\displaystyle L^{2^*(s_2)}\left(\Omega, \frac{dx}{|x|^{s_2}}\right)$.
Then
going to a subsequence if necessary, we can define
$$\tilde{m}_1:=\lim_{n\rightarrow \infty}\Psi_{\lambda}(u_n)$$
and
$$\tilde{m}_2:=\lim_{n\rightarrow \infty}\Upsilon_{\mu}(v_n).$$
Moreover, if $u=0$ or $v=0$, then the following are satisfied:
\begin{itemize}
\item[$(i)$]$\displaystyle m=\tilde{m}_1+\tilde{m}_2.$\\
\item[$(ii)$]$\{u_n\}$ is a $(PS)_{\tilde{m}_1}$ sequence for $\Psi_{\lambda}$ and $\{v_n\}$ is a $(PS)_{\tilde{m}_2}$ sequence for $\Upsilon_{\mu}$.\\
\item[$(iii)$] either $\tilde{m}_1=0$ or $\tilde{m}_1\geq c_{\sigma_0,\lambda^*,\lambda}$;\quad
either  $\tilde{m}_2=0$ or $\tilde{m}_2\geq c_{\eta_0,\mu^*,\mu}$.
\end{itemize}
\el
\bp
Since $\{(u_n,v_n)\}$ is bounded in $\mathscr{D}$, it is easy to see that $\Psi_{\lambda}(u_n)$ is bounded, so does $\Upsilon_{\mu}(v_n)$. Hence, up to a subsequence, $\tilde{m}_1$ and $\tilde{m}_2$ are well defined.
Since $\Phi'(u_n,v_n)\rightarrow 0$ in $\mathscr{D}^*$, we have
$$\Psi'_\lambda(u_n)-\kappa \alpha\frac{|u_n|^{\alpha-2}u_n|v_n|^\beta}{|x|^{s_2}}\rightarrow 0\;\hbox{in}\;E^*$$
and
$$\Upsilon'_\mu(v_n)-\kappa \beta\frac{|u_n|^{\alpha}|v_n|^{\beta-2}v_n}{|x|^{s_2}}\rightarrow 0\;\hbox{in}\;F^*.$$
Hence, if $u=0$ or $v=0$, by Corollary  \ref{2014-4-11-Cro1}, we have
$$\lim_{n\rightarrow \infty}\Psi'_{\lambda}(u_n)=\kappa \alpha\lim_{n\rightarrow \infty}\frac{|u_n|^{\alpha-2}u_n|v_n|^\beta}{|x|^{s_2}}=0$$
and
$$\lim_{n\rightarrow \infty}\Upsilon'_{\mu}(v_n)=\kappa\beta\lim_{n\rightarrow \infty}\frac{|u_n|^{\alpha}|v_n|^{\beta-2}v_n}{|x|^{s_2}}=0.$$
Hence,  $(ii)$ is proved.
We also have
$$\int_\Omega \frac{|u_n|^\alpha|v_n|^\beta}{|x|^{s_2}}dx=o(1)\|v_n\|\rightarrow 0$$
due to the boundedness of $\{v_n\}$.
Then it follows that
\begin{align*}
m=&\lim_{n\rightarrow \infty}\Phi(u_n,v_n)\\
=&\lim_{n\rightarrow \infty}\Psi_{\lambda}(u_n)+\lim_{n\rightarrow \infty}\Upsilon_{\mu}(v_n)-\kappa\lim_{n\rightarrow \infty}\int_\Omega \frac{|u_n|^\alpha|v_n|^\beta}{|x|^{s_2}}dx\\
=&\tilde{m}_1+\tilde{m}_2.
\end{align*}
By $(ii)$, it is easy to see $\tilde{m}_1\geq 0, \tilde{m}_2\geq 0$ since $2^*(s_1)>2, \lambda>0,\mu>0$. If $\tilde{m}_1\neq 0$, then $u_n\not\rightarrow 0$. By the definition of $c_{\sigma_0,\lambda^*,\lambda}$, we see that
$$\tilde{m}_1\geq c_{\sigma_0, \lambda^*,\lambda}.$$
Similarly, if $\tilde{m}_2\neq 0$, we have $$\tilde{m}_2\geq c_{\eta_0,\mu^*,\mu}.$$
\ep

\bc\lab{2014-4-11-Cro2}
Under the assumptions of Lemma \ref{2014-4-11-wl2}, let $\{(u_n,v_n)\}$ be a bounded $(PS)_m$ sequence with
$$0<m<\min\left\{c_{\sigma_0, \lambda^*,\lambda}, c_{\eta_0,\mu^*,\mu}\right\}.$$
Then $\Phi$ possesses a nontrivial critical point $(u,v)$,i.e., $u\neq 0, v\neq 0$.
\ec
\bp
It is an easy conclusion directly from Lemma \ref{2014-4-11-wl2}.
\ep

Next, we define
\be\lab{2015-4-25-we1}
I_\lambda(u)=\frac{1}{2}\int_{\R_+^N}|\nabla u|^2dx-\frac{\lambda}{2^*(s_1)}\int_{\R_+^N}\frac{|u|^{2^*(s_1)}}{|x|^{s_1}}dx
\ee
and denote the corresponding ground state value by $m_\lambda$.
Then, a direct calculation leads to that
$$m_1=\left(\frac{1}{2}-\frac{1}{2^*(s_1)}\right)\mu_{s_1}(\R_+^N)^{\frac{N-s_1}{2-s_1}}$$
and
$$m_\lambda=\lambda^{\frac{-2}{2^*(s_1)-2}}m_1.$$

\bc\lab{2014-4-12-Cro1}
Under the assumptions of Lemma \ref{2014-4-11-wl2},
let  $\{(u_n,v_n)\}\subset \mathscr{D}$ be a bounded sequence such that $(u_n,v_n)\rightharpoonup (u,v)$. Furthermore, we assume that $u_n\rightarrow u$ or $v_n\rightarrow v$ strongly in $H_0^1(\Omega)$ if $\alpha+\beta=2^*(s_2)$.
We have the following results:
\begin{itemize}
\item[$(i)$]  if $u=0$,  then either
$\tilde{m}_1=0$ or $\tilde{m}_1\geq m_{\lambda}$,
\item[$(ii)$]  if $v=0$,   then either
$\tilde{m}_2=0$ or $\tilde{m}_2\geq m_{\mu}$,
\end{itemize}
where $\tilde{m}_1, \tilde{m}_2$ are defined in Lemma \ref{2014-4-11-wl2}.
\ec
\bp
Recalling that $\{u_n\}$ is a $(PS)_{\tilde{m}_1}$ sequence for $\Psi_{\lambda}$ and $u_n\rightharpoonup u$ in $E=H_0^1(\Omega)$. By \cite[Theorem 3.1]{CeramiZhongZou.2015.DOI10.1007/s00526-015-0844-z.} or Theorem \ref{2014-6-19-th1}, we obtain that either $u_n\rightarrow u$ in $E$ or there exists $k$ functions $U^1,\cdots,U^k$ such that
$$\Psi_{\lambda}(u_n)=\Psi_{\lambda}(u)+\sum_{j=1}^{k}I_{\lambda}(U^j)+o(1),$$
where $u$ is a critical point of $\Psi_{\lambda}$ and $U^j, j=1,\dots, k$ are critical points of $I_{\lambda}$.
Hence, we have
$$\Psi_{\lambda}(u_n)=\Psi_{\lambda}(u)+o(1)$$
or
$$\Psi_{\lambda}(u_n)\geq \Psi_{\lambda}(u)+I_{\lambda}(U^1)+o(1).$$
Especially, when $u=0$, we have
$$\Psi_{\lambda}(u_n)=o(1)$$
or
$$\Psi_{\lambda}(u_n)\geq I_{\lambda}(U^1)+o(1).$$
Hence,
$$\tilde{m}_1=0\;\hbox{or}\;\tilde{m}_1\geq m_{\lambda}.$$
Similarly, we can prove that if $v=0$, we have
$$\tilde{m}_2=0\;\hbox{or}\;\tilde{m}_2\geq m_{\mu}.$$
\ep
\br\lab{2015-4-25-xbur1}
By Remark \ref{2015-4-25-xr1}, we see that Corollary \ref{2014-4-12-Cro1} is an improvement of Lemma \ref{2014-4-11-wl2}. And naturally we have the following result as an improvement of Corollary \ref{2014-4-11-Cro2}.
\er

\bc\lab{2015-4-25-cro1}
Under the assumptions of Lemma  \ref{2014-4-11-wl2}, let $\displaystyle\left\{(u_n,v_n)\right\}$ be a bounded $(PS)_m$ sequence for $\Phi$ and
\be
0<m<\min\{m_\lambda,m_\mu\}.
\ee
Then $\Phi$ possesses a nontrivial critical point $(u,v)$,i.e., $u\neq 0, v\neq 0$.
\ec
\bp
It is trivial by Corollary \ref{2014-4-12-Cro1}.
\ep

\br\lab{2015-4-26-xbur1}
Noting that for the critical couple case, i.e., $\alpha+\beta=2^*(s_2)$, we assume in Lemma  \ref{2014-4-11-wl2} that $u_n\rightarrow u$ or $v_n\rightarrow v$ strongly in $\displaystyle L^{2^*(s_2)}\left(\Omega, \frac{dx}{|x|^{s_2}}\right)$. This assumption is very strong, hence the Corollary \ref{2015-4-25-cro1} above is not enough for the case of $\alpha+\beta=2^*(s_2)$ in applications. We need a better description of the behavior of the $(PS)$ sequences of the functional $\Phi$.
\er

Here comes our main splitting results. We will discuss it by two cases of $\alpha+\beta<2^*(s_2)$ and $\alpha+\beta=2^*(s_2)$. The idea of proof is similar to the single equation case, however it is much more complicated for the system case. Hence, we will give the details.
\bt\lab{2014-4-13-th1}(Splitting for $\alpha+\beta<2^*(s_2)$)
Assume that $\Omega\subset \R^N$ is an open domain with $0\in \partial\Omega$, $\partial\Omega$ is $C^2$ at $0$, $0\leq s_2<2,    1<\alpha,1<\beta, \alpha+\beta<2^*(s_2)$ and one of the following holds:
 \begin{itemize}
\item[$(1)$]$\Omega$ is bounded, $0\neq \lambda^*>-\lambda_{1,\sigma_0}(\Omega), 0\neq \mu^*>-\lambda_{1,\eta_0}(\Omega)$;
\item[$(2)$]$\Omega$ is unbounded but $|\Omega|<\infty$, $0\neq \lambda^*>-\lambda_{1,\sigma_0}(\Omega), 0\neq \mu^*>-\lambda_{1,\eta_0}(\Omega)$;
\item[$(3)$]$|\Omega|=\infty$, $\lambda^*>0,\mu^*>0$ and conditions $(A_{\sigma_0}^*)$ and $(A_{\eta_0}^{*})$ are satisfied.
\end{itemize}
 Let $\{(u_n,v_n)\}$ be a bounded $(PS)_m$ sequence. Then there exists  a critical point  $(u,v)$ of $\Phi$ and two numbers $k, l\in \NN\cup \{0\}$ and functions $U^1,\cdots, U^k, V^1,\cdots,V^l$ and sequences of radius $r_n^j,\rho_n^i, 1\leq j\leq k, 1\leq i\leq l$ such that the following properties are satisfied going to a subsequence if necessary:
\begin{itemize}
\item[(a)] either $(u_n,v_n)\rightarrow (u,v)$ in $\mathscr{D}$ or\\
\item[(b1)]$v_n\rightarrow v$ in $F$ and $U^j\in D_{0}^{1,2}(\R_+^N)\subset D_{0}^{1,2}(\R^N)$ are nontrivial solutions of
    \be\lab{BP1}
\begin{cases}
\Delta u+\lambda\frac{|u|^{2^*(s_1)-2}u}{|x|^{s_1}}=0\quad &\hbox{in}\;\R_+^N,\\
u=0\;\hbox{on}\;\partial\R_+^N.
\end{cases}
\ee
\item[(b2)]$r_n^j\rightarrow 0$ as $n\rightarrow \infty$;\\
\item[(b3)] $\left\|u_n-u-\sum_{j=1}^{k} (r_n^j)^{\frac{2-N}{2}}U^j(\frac{\cdot}{r_n^j})\right\|\rightarrow 0$, where $\|\cdot\|$ is the norm in $D^{1,2}(\R^N)$;\\
\item[(b4)]$\|u_n\|^2\rightarrow \|u\|^2+\sum_{j=1}^{k}\|U^j\|^2$;\\
\item[(b5)] $\Phi(u_n,v_n)-\Phi(u,v)\rightarrow \sum_{j=1}^{k}I_\lambda(U^j)$ or\\
\item[(c1)]$u_n\rightarrow u$ in $E$ and $V^i\in D_{0}^{1,2}(\R_+^N)\subset D_{0}^{1,2}(\R^N)$ are nontrivial solutions of
    \be\lab{BP2}
\begin{cases}
\Delta v+\mu\frac{|v|^{2^*(s_1)-2}v}{|x|^{s_1}}=0\quad &\hbox{in}\;\R_+^N,\\
v=0\;\hbox{on}\;\partial\R_+^N.
\end{cases}
\ee
\item[(c2)]$\rho_n^i\rightarrow 0$ as $n\rightarrow \infty$;\\
\item[(c3)] $\left\|v_n-v-\sum_{i=1}^{l} (\rho_n^i)^{\frac{2-N}{2}}V^i(\frac{\cdot}{\rho_n^i})\right\|\rightarrow 0$, where $\|\cdot\|$ is the norm in $D^{1,2}(\R^N)$;\\
\item[(c4)]$\|v_n\|^2\rightarrow \|v\|^2+\sum_{i=1}^{l}\|V^i\|^2$;\\
\item[(c5)] $\Phi(u_n,v_n)-\Phi(u,v)\rightarrow \sum_{i=1}^{l}I_\mu(V^i)$ or\\
\item[(d1)]$U^j$ are nontrivial solutions of (\ref{BP1}) and $V^i$ are nontrivial solutions of (\ref{BP2});\\
\item[(d2)]$r_n^j\rightarrow 0, \rho_n^i\rightarrow 0$ as $n\rightarrow \infty$;\\
\item[(d3)] $\left\|u_n-u-\sum_{j=1}^{k} (r_n^j)^{\frac{2-N}{2}}U^j(\frac{\cdot}{r_n^j})\right\|\rightarrow 0$ and $\left\|v_n-v-\sum_{i=1}^{l} (\rho_n^i)^{\frac{2-N}{2}}V^i(\frac{\cdot}{\rho_n^i})\right\|\rightarrow 0$, where $\|\cdot\|$ is the norm in $D^{1,2}(\R^N)$;\\
\item[(d4)]$\|u_n\|^2\rightarrow \|u\|^2+\sum_{j=1}^{k}\|U^j\|^2$ and $\|v_n\|^2\rightarrow \|v\|^2+\sum_{i=1}^{l}\|V^i\|^2$;\\
\item[(d5)] $\Phi(u_n,v_n)-\Phi(u,v)\rightarrow \sum_{j=1}^{k}I_\lambda(U^j)+ \sum_{i=1}^{l}I_\mu(V^j)$.
\end{itemize}
\et
\bp
Since $\{(u_n,v_n)\}$ is bounded in $\mathscr{D}$, up to a subsequence, there exists some $(u,v)$ such that $(u_n,v_n)\rightharpoonup (u,v)$ in $\mathscr{D}$. Since $\{(u_n,v_n)\}$ is a   $(PS)_m$ sequence of $\Phi$, it is easy to see that $\Phi'(u,v)=0$. Then one of the following holds:
\begin{itemize}
\item[(a)]$(u_n,v_n)\rightarrow (u,v)$ in $\mathscr{D}$.\\
\item[(b)]$u_n\not\rightarrow u$ in $E$ but  $v_n\rightarrow v$ in $F$.\\
\item[(c)]$u_n\rightarrow u$ in $E$ but $v_n\not\rightarrow v$ in $F$.\\
\item[(d)]$u_n\not\rightarrow u$ in $E$ and  $v_n\not\rightarrow v$ in $F$.
\end{itemize}
We only prove the case of $(b1)-(b5)$ which is corresponding to the case of (b) above.
Since $\Phi'(u_n,v_n)\rightarrow 0$, we have
\be\lab{2014-4-13-we1}
\lim_{n\rightarrow \infty} \Psi'_\lambda(u_n)-\kappa\alpha\frac{|u_n|^{\alpha-2}u_n|v_n|^\beta}{|x|^{s_2}}=0.
\ee
Use Lemma \ref{2014-4-11-wl1} if $|\Omega|<\infty$ and apply  Lemma \ref{2014-6-6-zl1} if $|\Omega|=\infty$, we have
\be\lab{2014-4-13-we2}
\lim_{n\rightarrow \infty} \frac{|u_n|^{\alpha-2}u_n|v_n|^\beta}{|x|^{s_2}}=\frac{|u|^{\alpha-2}u|v|^\beta}{|x|^{s_2}}.
\ee
Noting that  $\Phi'(u,v)=0$, we have
\be\lab{2014-4-13-we3}
\Psi'_\lambda(u)=\kappa \alpha\frac{|u|^{\alpha-1}u|v|^\beta}{|x|^{s_2}}
\ee
Combined with (\ref{2014-4-13-we1}), (\ref{2014-4-13-we2}) and  (\ref{2014-4-13-we3}), we have
\be\lab{2014-4-13-we4}
\lim_{n\rightarrow \infty} \Psi'_\lambda(u_n)=\Psi'_\lambda(u).
\ee
By 3) of Lemma \ref{2014-3-19-l1}, we have
$$\frac{|u_n|^{2^*(s_1)-2}u_n}{|x|^{s_1}}-\frac{|u_n-u|^{2^*(s_1)-2}(u_n-u)}{|x|^{s_1}}\rightarrow \frac{|u|^{2^*(s_1)-2}u}{|x|^{s_1}}.$$
Hence,
\begin{align*}
&\lim_{n\rightarrow \infty}\Psi'_\lambda(u_n-u)\\
&= \lim_{n\rightarrow \infty}\Big(-\Delta (u_n-u)+\lambda^*\frac{(u_n-u)}{|x|^{\sigma_0}}-\lambda \frac{|u_n-u|^{2^*(s_1)-2}(u_n-u)}{|x|^{s_1}}\Big)\\
& =\lim_{n\rightarrow \infty}\Big[-\Delta (u_n-u)+\lambda^*\frac{(u_n-u)}{|x|^{\sigma_0}}-\lambda \Big(\frac{|u_n|^{2^*(s_1)-2}u_n}{|x|^{s_1}}-\frac{|u|^{2^*(s_1)-2}u}{|x|^{s_1}}\Big)\Big]\\
& =\lim_{n\rightarrow \infty} \Psi'_\lambda(u_n)-\Psi'_\lambda(u)\\
& = 0\quad (\hbox{by (\ref{2014-4-13-we4})}).
\end{align*}
Let $w_n=u_n-u$ and $c=m-\Phi(u,v)$, then $\{w_n\}$ is a $(PS)_c$ sequence of $\Psi_\lambda$. Notice that $w_n\rightharpoonup 0$ but $w_n\not\rightarrow 0$ in $E$. Hence, apply  \cite[Theorem 3.1]{CeramiZhongZou.2015.DOI10.1007/s00526-015-0844-z.} if $\Omega$ is bounded and  use Theorem \ref{2014-6-19-th1} if $\Omega$ is unbounded, we obtain that the case of $(b1)-(b5)$ happens.

By the similar arguments, if $u_n\rightarrow u$ in $E$ but $v_n\not\rightarrow v$ in $F$, we obtain (c1)-(c5); if $u_n\not\rightarrow u$ in $E$ and  $v_n\not\rightarrow v$ in $F$, (d1)-(d5) hold.
\ep

Next, we are going to establish the similar splitting result for the system case with critical couple terms. Denote
$$A(u,v):=I_\lambda(u)+I_\mu(v)-\kappa\int_{\R_+^N} \frac{|u|^\alpha|v|^\beta}{|x|^{s_2}}dx.$$

\br\lab{2015-4-25-xbur3}
Noting that $\mathscr{D}:=E\times F$, next we will also adopt the notations $\mathscr{P}:=D_{0}^{1,2}(\R^N)\times D_{0}^{1,2}(\R^N)$ and $\displaystyle \|(\phi,\psi)\|^{2}:=\int_\Omega \left(|\nabla \phi|^2+|\nabla \psi|^2\right)dx$ while
$\displaystyle \|(\phi,\psi)\|_{\mathscr{P}}^{2}:=\int_{\R^N} \left(|\nabla \phi|^2+|\nabla \psi|^2\right)dx$. For the convenience, we prefer to use the common notations $A(u,v)$ and $A'(u,v)$ for different domains. However, we will distinguish them by the words. For example, if for any $(h,\hbar)\in \mathscr{D}$, we have $\displaystyle\left\langle A'(u,v), (h,\hbar)\right\rangle=0$, then we will say that $A'(u,v)=(0,0)$ in $\mathscr{D}^*$; if for any $(h,\hbar)\in \mathscr{P}$, we have $\displaystyle\left\langle A'(u,v), (h,\hbar)\right\rangle=0$,then we will say that $A'(u,v)=(0,0)$ in $\mathscr{P}^*$. And if we only say that $A'(u,v)=(0,0)$, that means for any $(h,\hbar)\in D_{0}^{1,2}(\R_+^N)\times D_{0}^{1,2}(\R_+^N)$, there holds $\displaystyle\left\langle A'(u,v), (h,\hbar)\right\rangle=0$.
\er

\bl\lab{2014-6-10-l1}
Assume that $\Omega\subset \R^N$ is an open domain with $0\in \partial\Omega$, $\partial\Omega$ is $C^2$ at $0$, $0\leq s_2<2,   1<\alpha,1<\beta, \alpha+\beta=2^*(s_2)$ and one of the following holds:
 \begin{itemize}
\item[$(1)$]$\Omega$ is bounded, $0\neq \lambda^*>-\lambda_{1,\sigma_0}(\Omega), 0\neq \mu^*>-\lambda_{1,\eta_0}(\Omega)$;
\item[$(2)$]$\Omega$ is unbounded but $|\Omega|<\infty$, $0\neq \lambda^*>-\lambda_{1,\sigma_0}(\Omega), 0\neq \mu^*>-\lambda_{1,\eta_0}(\Omega)$;
\item[$(3)$]$|\Omega|=\infty$, $\lambda^*>0,\mu^*>0$.
\end{itemize}
 Let $\{(u_n,v_n)\}\subset \mathscr{D}$ such that $A(u_n,v_n)\rightarrow c$ and $A'(u_n,v_n)\rightarrow 0$ in $\mathscr{D}^*$. For $\{r_n\}\subset (0,+\infty)$ with $r_n\rightarrow 0$, assume that the rescaled sequences
$$\phi_n(x):=r_{n}^{\frac{N-2}{2}}u_n(r_nx),\;\psi_n(x):=r_{n}^{\frac{N-2}{2}}v_n(r_nx)$$
are such that
$$(\phi_n, \psi_n)\rightharpoonup (\phi, \psi)\;\hbox{in}\;\mathscr{P}:=D_{0}^{1,2}(\R^N)\times D_{0}^{1,2}(\R^N)$$
and
$$(\phi_n, \psi_n)\xrightarrow{a.e} (\phi,\psi)\;\hbox{on}\;\R^N.$$
Then $A'(\phi,\psi)=0$ and the sequences
$$\chi_n(x):=u_n(x)-r_{n}^{\frac{2-N}{2}}\phi(\frac{x}{r_n}),\;\varphi_n(x)
:=v_n(x)-r_{n}^{\frac{2-N}{2}}\psi(\frac{x}{r_n})$$
satisfy
$$A(\chi_n, \varphi_n)\rightarrow c-A(u,v),\;A'(\chi_n,\varphi_n)\rightarrow 0\;\hbox{in}\;\mathscr{D}^*$$
and
$$\|(\chi_n, \varphi_n)\|_{\mathscr{P}}^{2}=\|(u_n, v_n)\|_{\mathscr{P}}^{2}-\|(\phi, \psi)\|_{\mathscr{P}}^{2}+o(1),$$
where
$$\|(\phi, \psi)\|_{\mathscr{P}}^{2}:=\int_{\R^N}\big(|\nabla \phi|^2+|\nabla \psi|^2\big)dx.$$
\el
\bp
Without loss of generality, we assume  that $\partial R_+^N=\{x_N=0\}$ is tangent to $\partial\Omega$ at $0$, and that $-e_N=(0,\cdots, -1)$ is the outward normal to $\partial\Omega$ at that point. It is standard to prove that $supp\;\phi, supp\;\psi\subset \R_+^N$ (see \cite[Lemma 3.4]{GhoussoubKang.2004} or \cite[Lemma 3.3]{CeramiZhongZou.2015.DOI10.1007/s00526-015-0844-z.}).
Hence, we may assume that $(\phi_n, \psi_n)\rightharpoonup (\phi,\psi)$ in $D_{0}^{1,2}(\R^N)\times D_{0}^{1,2}(\R^N)$.
It is easy to see that $A(u,v)$ is invariant under dilation and we have
\begin{align*}
\|(\chi_n, \varphi_n)\|_{\mathscr{P}}^{2}=&\|r_{n}^{\frac{N-2}{2}}\big(\chi_n(r_nx), \varphi_n(r_nx)\big)\|_{\mathscr{P}}^{2}\\
=&\|(\phi_n-\phi, \psi_n-\psi)\|_{\mathscr{P}}^{2}\\
=&\|(\phi_n, \psi_n)\|_{\mathscr{P}}^{2}-\|(\phi,\psi)\|_{\mathscr{P}}^{2}+o(1)\\
=&\|(u_n,v_n)\|_{\mathscr{P}}^{2}-\|(\phi,\psi)\|_{\mathscr{P}}^{2}+o(1).
\end{align*}
By Lemma  \ref{2014-6-1-wl1}  and Lemma \ref{2014-3-19-l1}, we have
$$A(\phi_n-\phi, \psi_n-\psi)=A(\phi_n, \psi_n)-A(\phi,\psi)+o(1).$$
Hence
\begin{align*}
A(\chi_n,\varphi_n)=&A(r_{n}^{\frac{N-2}{2}}\chi_n(r_nx), r_{n}^{\frac{N-2}{2}}\varphi_n(r_nx))\\
=&A(\phi_n-\phi, \psi_n-\psi)\\
=&A(\phi_n, \psi_n)-A(\phi,\psi)+o(1)\\
=&A(u_n, v_n)-A(\phi,\psi)+o(1)\\
=&c-A(\phi,\psi)+o(1).
\end{align*}
For any $h, \hbar\in C_0^\infty(\R_+^N)$, let $h_n(x):=(r_n)^{\frac{2-N}{2}}h(\frac{x}{r_n}), \hbar_n:=(r_n)^{\frac{2-N}{2}}\hbar(\frac{x}{r_n})$, then we have that $h_n\times \hbar_n\in \mathscr{D}$ for $n$ large enough due to the assumption that $r_n\rightarrow 0$. Thus
\begin{align*}
\langle A'(\phi,\psi), (h,\hbar)\rangle=&\langle A'(\phi_n,\psi_n), (h,\hbar)\rangle+o(1)\\
=&\langle A'(u_n, v_n), (h_n, \hbar_n)\rangle +o(1)\\
=&o(1)\|(h_n, \hbar_n)\|^{2}+o(1)\\
=&o(1)\|(h, \hbar)\|_{\mathscr{P}}^{2}+o(1),
\end{align*}
which implies that $A'(\phi,\psi)=0$.

Next, for any $(h, \hbar)\in \mathscr{D}$, let $\big(h_n(x), \hbar_n(x)\big):=r_{n}^{\frac{N-2}{2}}\big(h(r_nx), \hbar(r_nx)\big)$, then we can see that for $n$ large enough, $supp\;h_n, supp\;\hbar_n\subset \R_+^N$. By Lemma \ref{2014-6-3-l1} and Lemma  \ref{2014-3-19-l1}, we have
\be\lab{2014-6-10-xe1}
A'(\phi_n,\psi_n)-A'(\phi_n-\phi, \psi_n-\psi)-A'(\phi,\psi)\rightarrow 0\;\hbox{in}\;H^{-1}(\R^N).
\ee
Hence, for $(h, \hbar)\in \mathscr{D}$, we have
\begin{align*}
& \langle A'(\chi_n, \varphi_n), (h,\hbar)\rangle\\
&=\langle A'(r_{n}^{\frac{N-2}{2}}\chi_n(r_nx), r_{n}^{\frac{N-2}{2}}\varphi_n(r_nx)) (h_n,\hbar_n)\rangle\\
&=\langle A'(\phi_n-\phi, \psi_n-\psi), (h_n,\hbar_n)\rangle\\
&=\langle A'(\phi_n, \psi_n), (h_n,\hbar_n)\rangle -\langle A'(\phi, \psi), (h_n,\hbar_n)\rangle \\
&\quad +o(1)\|(h_n, \hbar_n)\|_{\mathscr{P}}^{2}\; \; \hbox{  (by (\ref{2014-6-10-xe1}) )}\\
&=\langle A'(\phi_n, \psi_n), (h_n,\hbar_n)\rangle +o(1)\|(h_n, \hbar_n)\|_{\mathscr{P}}^{2}\;\; \hbox{(since $A'(\phi,\psi)=0$)}\\
&=\langle A'(u_n, v_n), (h,\hbar)\rangle +o(1)\|(h_n, \hbar_n)\|_{\mathscr{P}}^{2}\\
&=o(1)\|(h, \hbar)\|_{\mathscr{D}}^{2}\;\hbox{since $\|(h_n, \hbar_n)\|_{\mathscr{P}}^{2}=\|(h, \hbar)\|_{\mathscr{D}}^{2}$}.
\end{align*}
\ep

Now, it is enough for us to prove the splitting result for critical couple case.

\bt\lab{2014-6-16-th2}(Splitting Theorem for the critical case $\alpha+\beta=2^*(s_2)$)
Assume that $\Omega\subset \R^N$ is an open domain with $0\in \partial\Omega$, $\partial\Omega$ is $C^2$ at $0$, $0< s_2<2,  1<\alpha,1<\beta, \alpha+\beta=2^*(s_2)$ and one of the following holds:
 \begin{itemize}
\item[$(1)$]$\Omega$ is bounded, $0\neq \lambda^*>-\lambda_{1,\sigma_0}(\Omega), 0\neq \mu^*>-\lambda_{1,\eta_0}(\Omega)$;
\item[$(2)$]$\Omega$ is unbounded but $|\Omega|<\infty$, $0\neq \lambda^*>-\lambda_{1,\sigma_0}(\Omega), 0\neq \mu^*>-\lambda_{1,\eta_0}(\Omega)$;
\item[$(3)$]$|\Omega|=\infty$, $\lambda^*>0,\mu^*>0$.
\end{itemize}
 Let $\{(u_n,v_n)\}\subset \mathscr{D}$ be a bounded $(PS)_c$ sequence for the functional $\Phi$, i.e., $\Phi(u_n,v_n)\rightarrow c$ and $\Phi'(u_n,v_n)\rightarrow 0$ strongly in $\mathscr{D}^*$ as $n\rightarrow \infty$. Then there exists a critical point $(U^0, V^0)$ of $\Phi$, a number $k\in \NN$, $k$ pairs of functions $(U^1, V^1); \cdots ;(U^k, V^k)$ and $k$ sequences of radius
$(r^j_n)_n; r^j_n>0; 1 \leq j \leq k$ such that, up to a subsequence if necessary, the following
properties are satisfied. Either
\begin{itemize}
\item[(a)] $(u_n, v_n)\rightarrow (U^0, V^0)$ in $\mathscr{D}$ or
\item[(b)]  all the following  items are true:
\begin{itemize}
\item[(b1)]$(0,0)\neq (U^j,V^j)\in D_{0}^{1,2}(\R_+^N)\times D_{0}^{1,2}(\R_+^N)\subset D_{0}^{1,2}(\R^N) \times D_{0}^{1,2}(\R^N)$ are critical points of $A$;
\item[(b2)]$r_n^j\rightarrow 0$ as $n\rightarrow \infty$;
\item[(b3)] $\big\|\big(u_n-U^0-\sum_{j=1}^{k} (r_n^j)^{\frac{2-N}{2}}U^j(\frac{\cdot}{r_n^j}),v_n-V^0-\sum_{j=1}^{k} (r_n^j)^{\frac{2-N}{2}}V^j(\frac{\cdot}{r_n^j})\big)\big\|_{\mathscr{P}}\rightarrow 0$,
     where $\|\cdot\|_{\mathscr{P}}$ is defined in Lemma \ref{2014-6-10-l1};\\
\item[(b4)]$\|(u_n, v_n)\|_{\mathscr{P}}^{2}\rightarrow \sum_{j=0}^{k}\|(U^j, V^j)\|_{\mathscr{P}}^{2}$;\\
\item[(b5)] $\Phi(u_n, v_n)\rightarrow \Phi(U^0, V^0)+\sum_{j=1}^{k}A(U^j, V^j)$ with $$A(U^j, V^j)\geq c_\infty>0, $$
    where $c_\infty:=\inf\{A(u,v): (u,v)\neq (0,0), A'(u,v)=0\}$.
\end{itemize}
\end{itemize}
\et

\bp
Let $\{(u_n,v_n)\}\subset \mathscr{D}$ be a bounded $(PS)_c$ sequence of $\Phi(u,v)$, that is,
\be\lab{2014-6-10-pe1}
\begin{cases}
\Phi(u_n,v_n)\rightarrow c\;\hbox{in}\;\R,\\
\Phi'(u_n,v_n)\rightarrow 0\;\hbox{in}\;\mathscr{D}^*.
\end{cases}
\ee
Up to a subsequence, there exists some $(U^0, V^0)\in \mathscr{D}$ such that
$$
\begin{cases}
(u_n,v_n)\rightharpoonup (U^0, V^0)\;\hbox{in}\;\mathscr{D},\\
(u_n, v_n)\xrightarrow{a.e.} (U^0, V^0)\;\hbox{on}\;\R^N.
\end{cases}
$$
Then it is easy to see that $\Phi'(U^0, V^0)=0$. Denote $u_n^1:=u_n-U^0, v_n^1:=v_n-V^0$, then by Lemma \ref{2014-6-3-l1}  and Lemma \ref{2014-3-19-l1}, we have
\be\lab{2014-6-10-pe2}
\begin{cases}
\|(u_n^1, v_n^1)\|_{\mathscr{P}}^{2}=\|(u_n,v_n)\|_{\mathscr{P}}^{2}-\|(u,v)\|_{\mathscr{P}}^{2}+o(1),\\
A'(u_n^1, v_n^1)\rightarrow 0\;\hbox{in}\;\mathscr{D}^*,\\
A(u_n^1, v_n^1)\rightarrow c-\Phi(U^0, V^0).
\end{cases}
\ee
If $(u_n^1, v_n^1)\rightarrow 0$ in $\mathscr{D}$, we are done. If $(u_n^1, v_n^1)\not\rightarrow 0$ in $\mathscr{D}$, we claim that
\be\lab{2014-6-16-xe1}
\liminf_{n\rightarrow \infty} \Big(\int_\Omega \frac{|u_n^1|^{2^*(s_1)}}{|x|^{s_1}}dx+\int_\Omega \frac{|v_n^1|^{2^*(s_1)}}{|x|^{s_1}}dx\Big)>0.
\ee
Otherwise, the facts  that  $\{(u_n^1, v_n^1)\}\subset \mathscr{D}$ is bounded and  that $A'(u_n^1, v_n^1)\rightarrow 0$ in $\mathscr{D}^*$ would bring to a contradiction:
\begin{align*}
&\langle A'(u_n^1, v_n^1), (u_n^1,v_n^1)\rangle\\
 & =\int_\Omega \big(|\nabla u_n^1|^2+|\nabla v_n^1|^2\big) dx-\Big(\lambda\int_\Omega \frac{|u_n^1|^{2^*(s_1)}}{|x|^{s_1}}dx+\mu\int_\Omega \frac{|v_n^1|^{2^*(s_1)}}{|x|^{s_1}}dx\Big)\\
 &\quad -\kappa(\alpha+\beta)\int_\Omega \frac{|u_n^1|^\alpha |v_n^1|^\beta}{|x|^{s_2}}dx\\
 & \rightarrow  0.
\end{align*}
On the other hand, by \cite[Corollary 2.1 or
Corollary 2.2]{ZhongZou.2015arXiv:1504.01005v1[math.AP]4Apr2015}, we have
$$\liminf\int_\Omega \frac{|u_n^1|^{2^*(s_2)}}{|x|^{s_2}}dx=0$$
and
$$\liminf\int_\Omega \frac{|v_n^1|^{2^*(s_2)}}{|x|^{s_2}}dx=0$$
and then it follows from the H\"{o}lder inequality, we have
$$\big|\kappa(\alpha+\beta)\int_\Omega \frac{|u_n^1|^\alpha |v_n^1|^\beta}{|x|^{s_2}}dx\big|\rightarrow 0.$$
Hence, $(u_n^1,v_n^1)\rightarrow 0$ in $\mathscr{D}$. Thereby, (\ref{2014-6-16-xe1}) holds true. Define an analogue of Levy's concentration function
$$Q_n(r):=\int_{\mathbb{B}_r}\Big(\frac{|u_n^1|^{2^*(s_1)}}{|x|^{s_1}}+ \frac{|v_n^1|^{2^*(s_1)}}{|x|^{s_1}} \Big)dx.$$
Since $Q_n(0)=0$ and, due to (\ref{2014-6-16-xe1}), $Q_n(\infty)>0$ , we can choose  $\delta_1>0$ small enough such that
$$\mu_{s_1}(\R^N)^{-1}\delta_1^{\frac{2-s_1}{N-s_1}}<\frac{1}{4}$$
and find a positive sequence $\{r_n^1\}$ such that $Q_n(r_n^1)=\delta_1$.
Set $$\left(\tilde{u}_n^1(x),\tilde{v}_n^1(x)\right):=\left((r_n^1)^{\frac{N-2}{2}}u_n^1(r_n^1x), (r_n^1)^{\frac{N-2}{2}}v_n^1(r_n^1x)\right).$$
 Since $\|(u_n^1,v_n^1)\|_{\mathscr{P}}=\|\big(\tilde{u}_n^1, \tilde{v}_n^1\big)\|_{\mathscr{P}}$ is bounded, we may assume $(u_n^1,v_n^1)\rightharpoonup (U^1, V^1)$ in $D_{0}^{1,2}(\R^N)\times D_{0}^{1,2}(\R^N)$, $(u_n^1,v_n^1)\xrightarrow{a.e.} (U^1, V^1)$ in $\R^N$ and
\be\lab{2014-3-21-xe1}
\delta_1=\int_{\mathbb{B}_1}\Big(\frac{|\tilde{u}_n^1|^{2^*(s_1)}}{|x|^{s_1}}+ \frac{|\tilde{v}_n^1|^{2^*(s_1)}}{|x|^{s_1}}\Big)dx.
\ee
Since $0<s_1,s_2<2$,  by \cite[Corollary 2.1 or
Corollary 2.2]{ZhongZou.2015arXiv:1504.01005v1[math.AP]4Apr2015} again, we can also choose $\delta_1$ small enough such that
\be\lab{2014-6-16-we1}
\frac{|\kappa|(\alpha^2+\beta^2)}{2^*(s_2)}\mu_{s_2}(\R^N)^{-1}\left(\int_{\mathbb{B}_1}\frac{|\tilde{u}_n^1|^{2^*(s_2)}}{|x|^{s_2}}\right)^{\frac{2^*(s_2)-2}{2^*(s_2)}}<\frac{1}{4},
\ee
\be\lab{2014-6-16-we2}
\frac{|\kappa|(\alpha^2+\beta^2)}{2^*(s_2)}\mu_{s_2}(\R^N)^{-1}\left(\int_{\mathbb{B}_1}\frac{|\tilde{v}_n^1|^{2^*(s_2)}}{|x|^{s_2}}\right)^{\frac{2^*(s_2)-2}{2^*(s_2)}}<\frac{1}{4},
\ee
We claim that $(U^1, V^1)\not\equiv (0,0)$.

 Set $\Omega_n=\frac{1}{r_n^1}\Omega$ and let $f_n\in D_{0}^{1,2}(\Omega)$ be such that   $\displaystyle \langle A'(u_n^1, v_n^1), (h,0)\rangle=\int_\Omega \nabla f_n\cdot \nabla h$
for any $h\in D_{0}^{1,2}(\Omega)$. Then $g_n(x):=(r_n^1)^{\frac{N-2}{2}}f_n(r_n^1x)$ satisfies
$\displaystyle \int_{\Omega_n}|\nabla g_n|^2 =\int_\Omega |\nabla f_n|^2$ and $\displaystyle \langle A'(\tilde{u}_n^1, \tilde{v}_n^1), (h,0)\rangle=\int_{\Omega_n}\nabla g_n\cdot \nabla h$ for any $h\in D_{0}^{1,2}(\Omega_n)$.

 If $U^1\equiv 0$ would be true, then we could deduce that,
for any $h\in C_0^\infty(\R^N)$ with $supp\;h\subset \mathbb{B}_1,$
\begin{align*}
&\int_{\mathbb{B}_1}|\nabla (h\tilde{u}_n^1)|^2=\int_{\mathbb{B}_1}\nabla \tilde{u}_n^1\cdot \nabla (h^2\tilde{u}_n^1)+o(1)\\
&=\int_{\mathbb{B}_1}\frac{h^2|\tilde{u}_n^1|^{2^*(s_2)}}{|x|^{s_2}}+\int_{\mathbb{B}_1}\nabla g_n\cdot \nabla(h^2\tilde{u}_n^1)+
\kappa \alpha \int_{\mathbb{B}_1}\frac{|\tilde{u}_n^1|^\alpha |\tilde{v}_n^1|^\beta h^2}{|x|^{s_2}}+o(1)\\
&\leq  \int_{\mathbb{B}_1}\frac{h^2|\tilde{u}_n^1|^{2^*(s_2)}}{|x|^{s_2}}+o(1)\\
&\quad +|\kappa|\alpha\frac{\alpha}{2^*(s_2)}\int_{\mathbb{B}_1}\frac{|\tilde{u}_n^1|^{2^*(s_2)}h^2}{|x|^{s_2}}
+|\kappa|\alpha\frac{\beta}{2^*(s_2)}\int_{\mathbb{B}_1}\frac{|\tilde{v}_n^1|^{2^*(s_2)}h^2}{|x|^{s_2}}\\
&\leq  \mu_{s_1}(\R^N)^{-1}\Big(\int_{supp\; h}\frac{|\tilde{u}_n^1|^{2^*(s_1)}}{|x|^{s_1}}\Big)^{\frac{2^*(s_1)-2}{2^*(s_1)}}
\int_{\mathbb{B}_1} |\nabla (h\tilde{u}_n^1)|^2+o(1)\\
&\quad +\frac{|\kappa|\alpha^2}{2^*(s_2)}\mu_{s_2}(\R^N)^{-1}\Big(\int_{supp\; h}\frac{|\tilde{u}_n^1|^{2^*(s_2)}}{|x|^{s_2}}\Big)^{\frac{2^*(s_2)-2}{2^*(s_2)}}
\int_{\mathbb{B}_1} |\nabla (h\tilde{u}_n^1)|^2\\
&\quad +\frac{|\kappa|\alpha\beta}{2^*(s_2)}\mu_{s_2}(\R^N)^{-1}\Big(\int_{supp\; h}\frac{|\tilde{v}_n^1|^{2^*(s_2)}}{|x|^{s_2}}\Big)^{\frac{2^*(s_2)-2}{2^*(s_2)}}
\int_{\mathbb{B}_1} |\nabla (h\tilde{u}_n^1)|^2\\
&<\frac{3}{4}\int_{\mathbb{B}_1} |\nabla (h\tilde{u}_n^1)|^2+o(1),
\end{align*}
here we  use  the Young inequality and the following inequality:
$$\int_{\R^N} \frac{h^2|u|^{2^*(s)}}{|x|^s}\leq \mu_s(\R^N)^{-1}\big(\int_{supp\;h}\frac{|u|^{2^*(s)}}{|x|^s}\big)^{\frac{2^*(s)-2}{2^*(s)}}\int_{\R^N} |\nabla (hu)|^2,$$
which  holds for all $u\in D_{0}^{1,2}(\R^N)$ and $h\in C_0^\infty(\R^N)$ (see \cite[Lemma 3.5.]{GhoussoubKang.2004}).

\vskip0.123in

Hence, $\nabla \tilde{u}_n^1\rightarrow 0$ in $\displaystyle L^{2}_{loc}\left(\mathbb{B}_1\right)$. Similarly, if $V^1\equiv 0$, we have
$\nabla \tilde{v}_n^1\rightarrow 0$ in $\displaystyle L^{2}_{loc}\left(\mathbb{B}_1\right)$.
Thus, by the Hardy-Sobolev inequality,  $\tilde{u}_n^1\rightarrow 0, \tilde{v}_n^1\rightarrow 0$ in $\displaystyle L^{2^*(s_1)}_{loc}\left(\mathbb{B}_1, \frac{dx}{|x|^{s_1}}\right)$, which  follows  a contradiction with  (\ref{2014-3-21-xe1}).  Thus, the claim is proved.  Now, let us prove that $r_n^1\rightarrow 0$. Otherwise, when $\Omega$ is bounded, we see that $\{r_n^1\}$ is a bounded sequence. And if $\Omega$ is unbounded,  under our assumptions, similar to Theorem \ref{2014-6-19-th1}, we can still guarantee the boundedness of $\{r_n^1\}$ when $\delta_1$ is small enough. Up to a subsequence, we may assume that $r_n^1\rightarrow r_\infty^1>0$, then the fact that $(u_n^1, v_n^1)\rightharpoonup (0,0)$ in $\mathscr{D}$ would imply  $$\left(\tilde{u}_n^1(x), \tilde{v}_n^1(x)\right):=\left((r_n^1)^{\frac{N-2}{2}}u_n^1(r_n^1x), (r_n^1)^{\frac{N-2}{2}}v_n^1(r_n^1x)\right)\rightharpoonup (0,0)   \;\;\hbox{  in  }  D_{0}^{1,2}(\R^N)\times D_{0}^{1,2}(\R^N),$$  in contradiction with  $(U^1, V^1)\not\equiv (0,0)$. Therefore,  $r_n^1\rightarrow 0$.

\vskip0.123in

Next, we prove that $supp\;U^1\subset \R_+^N, supp\;V^1\subset \R_+^N$. Without loss of generality, we assume  that   $\partial R_+^N=\{x_N=0\}$ is tangent to $\partial\Omega$ at $0$, and that $-e_N=(0,\cdots, -1)$ is the outward normal to $\partial\Omega$ at that point. For any compact $K\subset \R_-^N$, we have for $n$ large enough, that $\frac{\Omega}{r_n^1}\cap K=\emptyset$   as $r_n^1\rightarrow 0$. Since $supp\;\tilde{u}_n^1\subset \frac{\Omega}{r_n^1}, supp\;\tilde{v}_n^1\subset \frac{\Omega}{r_n^1}$ and $\big(\tilde{u}_n^1, \tilde{v}_n^1\big)\rightarrow (U^1, V^1)$ a.e. in $\R^N$,  $(U^1, V^1)=(0,0)$ a.e. on $K$ follows, and, therefore, $supp\;U^1\subset \R_+^N, supp\;V^1\subset \R_+^N$.

\vskip0.123in

By (\ref{2014-6-10-pe2}) and Lemma \ref{2014-6-10-l1}, we obtain that $A'(U^1, V^1)=0$. Hence, $(U^1,V^1)$ is a critical point of $A$.  Moreover, by Lemma \ref{2014-6-10-l1} again, we see  that the sequence $\big(u_n^2(x), v_n^2(x)\big):=\big(u_n^1(x)-(r_n^1)^{\frac{2-N}{2}}U^1(\frac{x}{r_n^1}), v_n^1(x)-(r_n^1)^{\frac{2-N}{2}}V^1(\frac{x}{r_n^1})\big)$ also satisfies
\be\lab{2014-6-16-we3}
\begin{cases}
\|(u_n^2, v_n^2)\|_{\mathscr{P}}^{2}=\|(u_n,v_n)\|_{\mathscr{P}}^{2}-\|(U^0,V^0)\|_{\mathscr{P}}^{2}-\|(U^1,V^1)\|_{\mathscr{P}}^{2}+o(1),\\
A'(u_n^2, v_n^2)\rightarrow 0\;\hbox{in}\;\mathscr{D}^*,\\
A(u_n^2, v_n^2)\rightarrow c-\Phi(U^0, V^0)-A(U^1,V^1).
\end{cases}
\ee
We also note that $(U^1, V^1)\neq (0,0)$,
$$A(U^1,V^1)\geq c_\infty.$$
By iterating the above procedure, we construct critical points $(U^j,V^j)$ of $A$ and sequences $(r_n^j)$ with the above properties, and, since $\Phi(u_n, v_n) \rightarrow c,$  the iteration must terminate after a finite number of steps.
\ep

\br
Recalling that
$$\|(u_n^i, v_n^i)\|^2:=\int_\Omega \left(|\nabla u_n^i|^2+|\nabla v_n^i|^2\right)dx$$
and
$$
\|(u_n^i, v_n^i)\|_{\mathscr{D}}^2:=\int_\Omega \left(|\nabla u_n^i|^2+\lambda^* \frac{|u_n^i|^2}{|x|^{\sigma_0}}+|\nabla v_n^i|^2+\mu^*\frac{|v_n^i|^2}{|x|^{\eta_0}}\right)dx,
$$
we using the fact that for $i\geq 1$, $\|(u_n^i, v_n^i)\|_{\mathscr{D}}^2=\|(u_n^i, v_n^i)\|^2+o(1)$ in the proof above since that
$$
\int_\Omega \frac{|u_n^i|^2}{|x|^{\sigma_0}}dx=o(1), \int_\Omega \frac{|v_n^i|^2}{|x|^{\eta_0}}dx=o(1)
$$
by Lemma \ref{2014-6-20-l3}.
\er

\s{Variational identity and its applications}
\renewcommand{\theequation}{5.\arabic{equation}}
\renewcommand{\theremark}{5.\arabic{remark}}
\renewcommand{\thedefinition}{5.\arabic{definition}}
For any $(u,v)\in \mathscr{D}$, we denote that
\be\lab{2014-4-9-xe1}
\begin{cases}
a(u,v):=\|(u,v)\|_{\mathscr{D}}^{2},\\
b(u,v):=\lambda\int_\Omega \frac{|u|^{2^*(s_1)}}{|x|^{s_1}}dx+\mu\int_\Omega \frac{|v|^{2^*(s_1)}}{|x|^{s_1}}dx,\\
c(u,v):=\int_{\Omega} \frac{|u|^\alpha |v|^{\beta}}{|x|^{s_2}}dx,
\end{cases}
\ee
and define
\begin{align*}
\Phi(u,v)&:=\frac{1}{2}a(u,v)-\frac{1}{2^*(s_1)}b(u,v)
-\kappa c(u,v)
\end{align*}
for all $(u,v)\in\mathscr{D}$.
\bl
Any solution $(u,v)$ of
$$\begin{cases}
 -\Delta u+\lambda^*\frac{u}{|x|^{\sigma_0}}-\lambda \frac{|u|^{2^*(s_1)-2}u}{|x|^{s_1}}=\kappa\alpha \frac{1}{|x|^{s_2}}|u|^{\alpha-2}u|v|^\beta\quad &\hbox{in}\;\Omega,\\
-\Delta v+\mu^*\frac{v}{|x|^{\eta_0}}-\mu \frac{|v|^{2^*(s_1)-2}v}{|x|^{s_1}}=\kappa\beta \frac{1}{|x|^{s_2}}|u|^{\alpha}|v|^{\beta-2}v\quad &\hbox{in}\;\Omega,\\
(u,v)\in \mathscr{D},
\end{cases}$$
satisfies
\begin{align}\lab{2014-6-4-xe1}
&\frac{1}{N-2}\int_{\partial\Omega}|\nabla (u,v)|^2 x\cdot\nu d\sigma\nonumber\\
& =-\frac{N-\sigma_0}{N-2}\lambda^*\int_\Omega \frac{|u|^2}{|x|^{\sigma_0}}dx-\frac{N-\eta_0}{N-2}\mu^*\int_\Omega \frac{|v|^2}{|x|^{\eta_0}}dx+\lambda\int_\Omega \frac{|u|^{2^*(s_1)}}{|x|^{s_1}}dx\nonumber\\
&\quad +\mu\int_\Omega \frac{|v|^{2^*(s_1)}}{|x|^{s_1}}dx+2^*(s_2)\kappa\int_\Omega \frac{|u|^\alpha|v|^\beta}{|x|^{s_2}}dx-\int_\Omega |\nabla (u,v)|^2dx.
\end{align}
In particular,  if $ \Omega=\R_+^N$ or  $\R^N$, then
\begin{align}\lab{2014-6-4-xe2}
&\int_\Omega |\nabla (u,v)|^2dx+\frac{N-\sigma_0}{N-2}\lambda^*\int_\Omega \frac{|u|^2}{|x|^{\sigma_0}}dx+\frac{N-\eta_0}{N-2}\mu^*\int_\Omega \frac{|v|^2}{|x|^{\eta_0}}dx\nonumber\\
&-\lambda\int_\Omega \frac{|u|^{2^*(s_1)}}{|x|^{s_1}}dx -\mu\int_\Omega \frac{|v|^{2^*(s_1)}}{|x|^{s_1}}dx-2^*(s_2)\kappa\int_\Omega \frac{|u|^\alpha|v|^\beta}{|x|^{s_2}}dx=0
\end{align}
\el
\bp  We take
\begin{align*}
&G(x, u, v)=\\
&-\frac{1}{2}\lambda^*\frac{|u|^2}{|x|^{\sigma_0}}-\frac{1}{2}\mu^*\frac{|v|^2}{|x|^{\eta_0}}
+\frac{1}{2^*(s_1)}\lambda \frac{|u|^{2^*(s_1)}}{|x|^{s_1}}+\frac{1}{2^*(s_1)}\mu\frac{|v|^{2^*(s_1)}}{|x|^{s_1}}
+\kappa \frac{|u|^\alpha|v|^\beta}{|x|^{s_2}}.
\end{align*}
A direct calculation shows that
\begin{align*}
\sum_{i=1}^{N}x_iG_{x_i}(x, u, v)=&\frac{\sigma_0}{2}\lambda^*\frac{|u|^2}{|x|^{\sigma_0}}+\frac{\eta_0}{2}\mu^*\frac{|v|^2}{|x|^{\eta_0}}-\frac{s_1}{2^*(s_1)}\lambda \frac{|u|^{2^*(s_1)}}{|x|^{s_1}}\\
&-\frac{s_1}{2^*(s_1)}\mu\frac{|v|^{2^*(s_1)}}{|x|^{s_1}}-\kappa s_2 \frac{|u|^\alpha |v|^\beta}{|x|^{s_2}}.
\end{align*}
Then by \cite[Proposition 7.1]{ZhongZou.2015arXiv:1504.01005v1[math.AP]4Apr2015},
we obtain (\ref{2014-6-4-xe1}) and  (\ref{2014-6-4-xe2}).
\ep

\bc\lab{2014-6-4-xcro1}
Any solution $(u,v)$ of
$$\begin{cases}
-\Delta u+\lambda^*\frac{u}{|x|^{\sigma_0}}-\lambda \frac{|u|^{2^*(s_1)-2}u}{|x|^{s_1}}=\kappa\alpha \frac{1}{|x|^{s_2}}|u|^{\alpha-2}u|v|^\beta\quad &\hbox{in}\;\Omega,\\
-\Delta v+\mu^*\frac{v}{|x|^{\eta_0}}-\mu \frac{|v|^{2^*(s_1)-2}v}{|x|^{s_1}}=\kappa\beta \frac{1}{|x|^{s_2}}|u|^{\alpha}|v|^{\beta-2}v\quad &\hbox{in}\;\Omega,\\
(u,v)\in \mathscr{D},
\end{cases}$$
satisfies
\begin{align}\lab{2014-6-4-xe3}
&\int_{\partial\Omega}|\nabla (u,v)|^2 x\cdot \nu d\sigma +(2-\sigma_0)\lambda^*\int_\Omega \frac{|u|^2}{|x|^{\sigma_0}}dx+(2-\eta_0)\mu^*\int_\Omega \frac{|v|^2}{|x|^{\eta_0}}dx\nonumber\\
&-(N-2)\big(2^*(s_2)-\alpha-\beta\big)\kappa \int_\Omega \frac{|u|^\alpha|v|^\beta}{|x|^{s_2}}dx=0.
\end{align}
Moreover,  if $\Omega=\R_+^N$ or $\R^N$, then
\begin{align}\lab{2014-6-4-xe4}
&(2-\sigma_0)\lambda^*\int_\Omega \frac{|u|^2}{|x|^{\sigma_0}}dx+(2-\eta_0)\mu^*\int_\Omega \frac{|v|^2}{|x|^{\eta_0}}dx\nonumber\\
&-(N-2)\big(2^*(s_2)-\alpha-\beta\big)\kappa \int_\Omega \frac{|u|^\alpha|v|^\beta}{|x|^{s_2}}dx=0.
\end{align}
\ec
\bp
Since $(u,v)$ is a solution, we have $\langle\Phi'(u,v), (u,v)\rangle=0$. That is
\be\lab{2014-6-4-xe5}
a(u,v)-b(u,v)-\kappa(\alpha+\beta)c(u,v)=0.
\ee
On the other hand, by (\ref{2014-6-4-xe1}), we have that
\begin{align}\lab{2014-6-4-xe6}
& a(u,v)-b(u,v)-\kappa(\alpha+\beta)c(u,v)\nonumber\\
& =-\frac{2-\sigma_0}{N-2}\lambda^*\int_\Omega \frac{|u|^2}{|x|^{\sigma_0}}dx-\frac{2-\eta_0}{N-2}\mu^*\int_\Omega\frac{|v|^2}{|x|^{\eta_0}}dx\nonumber\\
&\quad +\big(2^*(s_2)-\alpha-\beta\big)\kappa \int_\Omega \frac{|u|^\alpha|v|^\beta}{|x|^{s_2}}dx-\frac{1}{N-2}\int_{\partial\Omega}|\nabla (u,v)|^2 x\cdot \nu d\sigma.
\end{align}
By (\ref{2014-6-4-xe5}) and (\ref{2014-6-4-xe6}), we obtain (\ref{2014-6-4-xe3}). When $\Omega=\R_+^N$ or $\R^N$, by (\ref{2014-6-4-xe2}) and (\ref{2014-6-4-xe5}), we obtain (\ref{2014-6-4-xe4}).
\ep

\bc\lab{2014-6-4-xcro2}
If $\Omega$ is a star-shaped domain of $\R^N$ around $0$, especially $\Omega=\R_+^N$ or $\R^N$, $0\leq \sigma_0, \eta_0, s_2<2, 0<s_1<2, \lambda^*>0, \mu^*>0$, then $(0,0)$ is the unique solution of
$$\begin{cases}
-\Delta u+\lambda^*\frac{u}{|x|^{\sigma_0}}-\lambda \frac{|u|^{2^*(s_1)-2}u}{|x|^{s_1}}=\kappa\alpha \frac{1}{|x|^{s_2}}|u|^{\alpha-2}u|v|^\beta\quad &\hbox{in}\;\Omega,\\
-\Delta v+\mu^*\frac{v}{|x|^{\eta_0}}-\mu \frac{|v|^{2^*(s_1)-2}v}{|x|^{s_1}}=\kappa\beta \frac{1}{|x|^{s_2}}|u|^{\alpha}|v|^{\beta-2}v\quad &\hbox{in}\;\Omega,\\
(u,v)\in \mathscr{D},
\end{cases}$$
if furthermore $\big(2^*(s_2)-\alpha-\beta\big)\kappa\leq 0$.
\ec
\bp
Since $\Omega$ is a star-shaped domain around $0$, we have $x\cdot \nu>0$ on $\partial\Omega\backslash\{0\}$, where $\nu$ denote the outward unit normal to $\partial\Omega$. Then by (\ref{2014-6-4-xe3}), we can obtain the result. And when $\Omega=\R_+^N$ or $\R^N$, by (\ref{2014-6-4-xe4}), we can also obtain the result.
\ep

We prefer to give some this kind result about a single equation to close this section. The following results we refer to \cite{LehrerMaia.2013}.
\bo\lab{2014-6-21-prop1}(cf.\cite[Proposition 2.1]{LehrerMaia.2013})
Let $u\in H^1(\Omega)\backslash \{0\}$ be a solution of equation $-\Delta u=g(x, u)$ and $G(x, u)=\int_0^u g(x, s)ds$ is such that $G\big(\cdot, u(\cdot)\big)$ and $x_iG_{x_i}\big(\cdot, u(\cdot)\big)$ are in $L^1(\Omega)$, then $u$ satisfies:
$$\int_{\partial\Omega}|\nabla u|^2 x\cdot \eta dS_x=2N\int_\Omega G(x, u)dx+2\sum_{i=1}^{N}\int_\Omega x_iG_{x_i}(x, u)dx-(N-2)\int_\Omega|\nabla u|^2dx,$$
where $\Omega$ is a regular domain in $\R^N$ and $\eta$ denotes the unitary exterior normal vector to $\partial\Omega$. Moreover, if $\Omega=\R^N$, then
$$2N\int_{\R^N}G(x, u)dx+2\sum_{i=1}^{N}\int_{\R^N}x_iG_{x_i}(x, u)dx=(N-2)\int_{\R^N}|\nabla u|^2 dx.$$
\eo
\br\lab{2014-6-21-r1}
Since the processes are standard now, we give the following results without proof:\\
Let $\Omega=\R_+^N$ and $u\in D_{0}^{1,2}(\Omega)\backslash \{0\}$ be a solution of equation $-\Delta u=g(x, u)$ and $G(x, u)=\int_0^u g(x, s)ds$ is such that $G\big(\cdot, u(\cdot)\big)$ and $x_iG_{x_i}\big(\cdot, u(\cdot)\big)$ are in $L^1(\Omega)$, then $u$ satisfies
$$2N\int_{\R_+^N}G(x, u)dx+2\sum_{i=1}^{N}\int_{\R_+^N}x_iG_{x_i}(x, u)dx=(N-2)\int_{\R_+^N}|\nabla u|^2 dx.$$
\er\hfill$\Box$

\bc\lab{2014-6-21-cro1}
Let $\Omega=\R_+^N$ or $\R^N$, then problem
\be\lab{2014-6-21-xe2}
-\Delta u+\lambda \frac{|u|^{p-1}u}{|x|^{s_1}}=\frac{|u|^{2^*(s_2)-2}u}{|x|^{s_2}}, u\in E:=D_{0}^{1,2}(\Omega)\cap L^{p+1}\left(\Omega, \frac{dx}{|x|^{s_2}}\right).
\ee
has no nontrivial solution if $\lambda\neq 0$ and $p\neq 2^*(s_1)-1$.
\ec
\bp
Let
$g(x, u)=-\lambda\frac{|u|^{p-1}u}{|x|^{s_1}}+\frac{|u|^{2^*(s_2)-2}u}{|x|^{s_2}}$ and $u$ be a solution of (\ref{2014-6-21-xe2}). A direct calculation shows that
$$\sum_{i=1}^{N}x_iG_{x_i}(x, u)=\frac{\lambda s_1}{p+1}\frac{|u|^{p+1}}{|x|^{s_1}}-\frac{s_2}{2^*(s_2)}\frac{|u|^{2^*(s_2)}}{|x|^{s_2}}.$$
Then by Remark\ref{2014-6-21-r1} and Proposition\ref{2014-6-21-prop1}, we have
\be\lab{2014-6-21-xe3}
\int_{\Omega}|\nabla u|^2 dx-\int_\Omega \frac{|u|^{2^*(s_2)}}{|x|^{s_2}}dx=\frac{-2^*(s_1)\lambda}{p+1}\int_\Omega \frac{|u|^{p+1}}{|x|^{s_1}}dx.
\ee
On the other hand, tested by $u$, we have
\be\lab{2014-6-21-xe4}
\int_{\Omega}|\nabla u|^2 dx-\int_\Omega \frac{|u|^{2^*(s_2)}}{|x|^{s_2}}dx=-\lambda\int_\Omega \frac{|u|^{p+1}}{|x|^{s_1}}dx.
\ee
Hence, if $\lambda\neq 0$ and $p\neq 2^*(s_1)-1$, $u=0$ is the unique solution.
\ep


\s{Estimation on the least energy $m_0$}
\renewcommand{\theequation}{6.\arabic{equation}}
\renewcommand{\theremark}{6.\arabic{remark}}
\renewcommand{\thedefinition}{6.\arabic{definition}}
\br\lab{2015-4-26-r1}
Under the assumptions of this section, it is easy to prove the corresponding Nehari manifold $\mathcal{N}$ is well defined, and $\mathcal{N}$ is closed and bounded away from $(0,0)$. We can also prove that any $(PS)_c$ sequence of $\Phi$ is bounded in $\mathscr{D}$. And any bounded $PS$ sequence of $\Phi\big|_{\mathcal{N}}$ is also a bounded $PS$ sequence of $\Phi$. Since these proofs are very standard, we omit the details and refer to \cite[Section 4]{ZhongZou.2015arXiv:1504.01005v1[math.AP]4Apr2015}.
\er

\br\lab{2015-4-26-r2}
We always assume that $(H_1)$ and $(H_2)$ are satisfied, especially when $|\Omega|=\infty$, we assume further that $(A_{\sigma_0}^{*})$ and $(A_{\eta_0}^{*})$ hold(see subsection 4.2), which guarantee that $\displaystyle c(u,v):=\int_\Omega \frac{|u|^\alpha|v|^\beta}{|x|^{s_2}}dx$ is well defined for all $(u,v)\in \mathscr{D}:=E\times F$. We can obtain the following results which are parallel to those in \cite[subsection 6.2]{ZhongZou.2015arXiv:1504.01005v1[math.AP]4Apr2015} and \cite[subsection 2.2]{ZhongZou.2015arXiv:submit/1229202[math.AP]12Apr2015},
see also \cite{AbdellaouiFelliPeral.2009,ZhongZou.2015,ZhongZou.2015arXiv:1503.08917v1[math.AP]31Mar2015}.etc.
Although the framework spaces are different, the arguments depends only on the values of $\alpha$ and $\beta$, these proofs can be modified smoothly. Since the proofs are long but standard, we omit the details.
\er
\bl\lab{2015-4-26-zl1}
Assume $(H_1),(H_2), 1<\alpha,1<\beta, \alpha+\beta\leq 2^*(s_2)$, and if $|\Omega|=\infty,\alpha+\beta<2^*(s_2)$, we assume further that $(A_{\sigma_0}^{*})$ and $(A_{\eta_0}^{*})$ hold. Let $0<U\in E$ is a least energy solution of \eqref{2015-4-27-cle1} and define
\be\lab{2015-4-26-ze1}
\tilde{\eta}_1:=\inf_{v\in F\backslash \{0\}}\frac{\|v\|_F^2}{\int_\Omega \frac{|U|^{\alpha}|v|^2}{|x|^{s_2}}dx}.
\ee
Then
\begin{itemize}
\item[(1)] if $\kappa<0$ or $\beta>2$ or $\begin{cases}\beta=2\\ \kappa<\tilde{\eta}_1 \end{cases}$, $(U,0)$ is a local minimal of $\Phi$ in $\mathcal{N}$.
\item[(2)] if $\beta<2,\kappa>0$ or $\beta=2,\kappa>\tilde{\eta}_1$,
\be\lab{2015-4-26-ze2}
m_0:=\inf_{(\phi,\varphi)\in \mathcal{N}}\Phi(\phi,\varphi)<\Phi(U,0)=\Psi_\lambda(U)=c_{\sigma_0,\lambda^*,\lambda}.
\ee
\end{itemize}
\el
\bp
We omit it.
\ep

\bl\lab{2015-4-26-zl2}
Assume $(H_1),(H_2), 1<\alpha,1<\beta, \alpha+\beta\leq 2^*(s_2)$, and if $|\Omega|=\infty,\alpha+\beta<2^*(s_2)$, we assume further that $(A_{\sigma_0}^{*})$ and $(A_{\eta_0}^{*})$ hold. Let $0<V\in F$ is a least energy solution of \eqref{2015-4-25-xe1} and define
\be\lab{2015-4-26-ze3}
\tilde{\eta}_2:=\inf_{u\in E\backslash \{0\}}\frac{\|u\|_E^2}{\int_\Omega \frac{|V|^{\beta}|u|^2}{|x|^{s_2}}dx}.
\ee
Then
\begin{itemize}
\item[(1)] if $\kappa<0$ or $\alpha>2$ or $\begin{cases}\alpha=2\\ \kappa<\tilde{\eta}_2 \end{cases}$, $(0,V)$ is a local minimal of $\Phi$ in $\mathcal{N}$.
\item[(2)] if $\alpha<2,\kappa>0$ or $\alpha=2,\kappa>\tilde{\eta}_2$,
\be\lab{2015-4-26-ze4}
m_0:=\inf_{(\phi,\varphi)\in \mathcal{N}}\Phi(\phi,\varphi)<\Phi(0,V)=\Upsilon_\mu(V)=c_{\eta_0,\mu^*,\mu}.
\ee
\end{itemize}
\el
\bp
We omit it.
\ep

\bc\lab{2015-4-26-cor1}
Assume $(H_1),(H_2), 1<\alpha,1<\beta, \alpha+\beta\leq 2^*(s_2)$, and if $|\Omega|=\infty,\alpha+\beta<2^*(s_2)$, we assume further that $(A_{\sigma_0}^{*})$ and $(A_{\eta_0}^{*})$ hold. If $\beta<2$ or $\beta=2$ with $\kappa>\tilde{\eta}_1$, and if $\alpha<2$ or $\alpha=2$ with $\kappa>\tilde{\eta}_2$, then we have
\be\lab{2015-4-26-ze5}
m_0<\min\{c_{\sigma_0,\lambda^*,\lambda}, c_{\eta_0,\mu^*,\mu}\}.
\ee
\ec
\bp
By Lemma \ref{2015-4-26-zl1} and Lemma \ref{2015-4-26-zl2}.
\ep
\s{Proofs of Theorem \ref{2015-4-26-xth1} and Theorem \ref{2015-4-26-xth2}}
\renewcommand{\theequation}{7.\arabic{equation}}
\renewcommand{\theremark}{7.\arabic{remark}}
\renewcommand{\thedefinition}{7.\arabic{definition}}
\noindent{\bf Proof of Theorem \ref{2015-4-26-xth1}:}
We note that under the assumptions, it is easy to see that the corresponding Nehari manifold is well defined. Let $\{(u_n,v_n)\}\subset \mathcal{N}$ be a minimizing sequence, then it is standard to prove that $\{(u_n,v_n)\}$ is a bounded $(PS)_{m_0}$ sequence for $\Phi$. And by Corollary \ref{2015-4-26-cor1}, we have
\be\lab{2015-4-26-ze6}
m_0<\min\{c_{\sigma_0,\lambda^*,\lambda}, c_{\eta_0,\mu^*,\mu}\}.
\ee
Hence, by Corollary \ref{2014-4-11-Cro2}, we obtain that there exists some $(u,v)\in \mathscr{D}$ such that $(u_n,v_n)\rightarrow (u,v)$ and $u\neq 0, v\neq 0$. Thus, $(u,v)$ is a nontrivial ground state solution.
\hfill$\Box$

\vskip 0.2in
For the case of $\alpha+\beta=2^*(s_2)$, it is not easy to see that $u_n\rightarrow u$ or $v_n\rightarrow v$ strongly in $\displaystyle L^{2^*(s_2)}\left(\Omega, \frac{dx}{|x|^{s_2}}\right)$ (see Remark \ref{2015-4-26-xbur1}). Hence, for the critical couple case, we can note easily obtain the compactness result basing on Corollary \ref{2015-4-26-cor1} and Corollary \ref{2014-4-11-Cro2}. We also note that
\be\lab{2015-4-26-ze7}
\min\{c_{\sigma_0,\lambda^*,\lambda}, c_{\eta_0,\mu^*,\mu}\}<\left[\frac{1}{2}-\frac{1}{2^*(s_1)}\right]\left(\mu_{s_1}(\R_+^N)\right)^{\frac{2^*(s_1)}{2^*(s_1)-2}}\left(\max\{\lambda,\mu\}\right)^{-\frac{2}{2^*(s_2)-2}}.
\ee
We let $\tilde{\Theta}$ denote  the least energy of the limit problem, then by \cite[Theorem 7.2]{ZhongZou.2015arXiv:1504.01005v1[math.AP]4Apr2015} or \cite[Theorem 1.1]{ZhongZou.2015arXiv:submit/1229202[math.AP]12Apr2015}, we also have
\be\lab{2015-4-26-ze8}
\tilde{\Theta}<\left[\frac{1}{2}-\frac{1}{2^*(s_1)}\right]\left(\mu_{s_1}(\R_+^N)\right)^{\frac{2^*(s_1)}{2^*(s_1)-2}}\left(\max\{\lambda,\mu\}\right)^{-\frac{2}{2^*(s_2)-2}}.
\ee
Hence, we can not judge the relationship of size between $\min\{c_{\sigma_0,\lambda^*,\lambda}, c_{\eta_0,\mu^*,\mu}\}$ and $\tilde{\Theta}$ intuitively. To apply the splitting Theorem \ref{2014-6-16-th2}, we need a further detailed estimation. For the technical reasons, we need the  assumptions $(\tilde{H}_1)$ and $(\tilde{H}_2)$.

\vskip 0.2in
\noindent{\bf Proof of Theorem \ref{2015-4-26-xth2}:}
Similar to the proof of Theorem \ref{2015-4-26-xth1}, there exists a minimizing sequence $\{(u_n,v_n)\}\subset \mathcal{N}$ which is a bounded $(PS)_{m_0}$ sequence of $\Phi$. Moreover,
\be
m_0<\min\{c_{\sigma_0,\lambda^*,\lambda}, c_{\eta_0,\mu^*,\mu}\}.
\ee
Under the assumptions of $(\tilde{H}_1)$ and $(\tilde{H}_2)$, we shall prove that $m_0<\tilde{\Theta}$ and thus
\be
m_0<\min\{\tilde{\Theta},c_{\sigma_0,\lambda^*,\lambda}, c_{\eta_0,\mu^*,\mu}\}.
\ee
Then, applying the splitting Theorem \ref{2014-6-16-th2}, we see that there exists some $u\neq 0,v\neq 0$ such that $(u_n,v_n)\rightarrow (u,v)$. Hence, $(u,v)$ is a nontrivial ground state solution.

\vskip 0.02in
\noindent{\bf The proof of $m_0<\tilde{\Theta}$ under the assumptions $(\tilde{H}_1)$ and $(\tilde{H}_2)$:} Denote the functional corresponding to the limit equation by $A$,  under our assumptions, by \cite[Theorem 7.2]{ZhongZou.2015arXiv:1504.01005v1[math.AP]4Apr2015} or \cite[Theorem 1.1]{ZhongZou.2015arXiv:submit/1229202[math.AP]12Apr2015}, we see that the limite equation possesses a nontrivial ground state solution, i.e., there exists some $u\neq 0,v\neq 0$ such that $A(u,v)=\tilde{\Theta}$ and $(u,v)$ satisfies the following equation:
\be\lab{2015-4-26-we1}
\begin{cases}
-\Delta u-\lambda \frac{|u|^{2^*(s_1)-2}u}{|x|^{s_1}}=\kappa \alpha \frac{1}{|x|^{s_2}}|u|^{\alpha-2}u|v|^\beta \quad &\hbox{in}\;\R_+^N,\\
-\Delta v-\mu\frac{|v|^{2^*(s_2)-2}v}{|x|^{s_1}}=\kappa \beta \frac{1}{|x|^{s_2}}|u|^\alpha |v|^{\beta-2}v&\hbox{in}\;\R_+^N,\\
\kappa>0, (u,v)\in \mathcal{D}:=D_{0}^{1,2}(\R_+^N)\times D_{0}^{1,2}(\R_+^N).&
\end{cases}
\ee
Without any loss of generality, we may assume that, in a suitable neighborhood of $0$, $\partial\Omega$ can be represented as $x_N=\varphi(x')$, where $x'=(x_1,\cdots, x_{N-1}), \ \varphi(0)=0,\ \nabla'\varphi(0)=0, \ \nabla'=(\partial_1, \cdots, \partial_{N-1}),\ $ and that the outer normal to $\partial\Omega$ at $0$ is $-e_N=(0,\cdots, 0, -1)$. Define $\phi(x)=(x', x_N-\varphi(x'))$ to ``flatten out" the boundary. We can choose a small $r_0>0$ and neighborhoods of $0$, $U$ and $\tilde{U}$, such that
$\displaystyle\phi(U)=\mathbb{B}_{r_0}(0),\;\phi(U\cap \Omega)=\mathbb{B}_{r_0}^{+}(0),
\phi(\tilde{U})=\mathbb{B}_{\frac{r_0}{2}}(0)$ and $\phi(\tilde{U}\cap \Omega)=\mathbb{B}_{\frac{r_0}{2}}^{+}(0),$ where, for any $r_0>0,\ \mathbb{B}_{r_0}(0) \subset \R^N$ is an open  ball of radius $r_0$  centered at $0$ and
$\mathbb{B}_{r_0}^{+}(0)=\ \mathbb{B}_{r_0}(0)\cap \R_+^N .$
By the smoothness assumption on $\partial\Omega$, $\varphi$ can be written as
\be\lab{2014-3-1-e0}
\varphi(y')=\sum_{i=1}^{N-1}\alpha_i y_i^2+o(|y'|^2),
\ee
then
$$H(0)=\frac{1}{N-1}\sum_{i=1}^{N-1}\alpha_i.$$
For any $\varepsilon>0$, define
\be u_\varepsilon(x):=\varepsilon^{-\frac{N-2}{2}}u\left(\frac{\phi(x)}{\varepsilon}\right),
v_\varepsilon(x):=\varepsilon^{-\frac{N-2}{2}}v\left(\frac{\phi(x)}{\varepsilon}\right)\;\;\;x\in \Omega\cap U,\ee
and let
\be \hat{u}_\varepsilon(x)=\eta u_\varepsilon(x), \hat{v}_\varepsilon(x)=\eta v_\varepsilon(x).\ee
For the estimation on the gradient terms, similar to \cite[Lemma 3.4]{ZhongZou.2015arXiv:1504.00730[math.AP]}, we have
\be\lab{2015-3-9-e5}
\int_\Omega \left|\nabla \hat{u}_\varepsilon(x)\right|^2dx=\int_{\R_+^N}\left|\nabla u(y)\right|^2dy-2\int_{\mathbb{B}_{\frac{r_0}{\varepsilon}}^{+}}\eta\left(\phi^{-1}(\varepsilon y)\right)^2\partial_N u(y) \nabla' u(y)\cdot (\nabla'\varphi)(\varepsilon y')dy+O(\varepsilon^2)
\ee
and
\be\lab{2015-3-9-e5-1}
\int_\Omega \left|\nabla \hat{v}_\varepsilon(x)\right|^2dx=\int_{\R_+^N}\left|\nabla v(y)\right|^2dy-2\int_{\mathbb{B}_{\frac{r_0}{\varepsilon}}^{+}}\eta\left(\phi^{-1}(\varepsilon y)\right)^2\partial_N v(y) \nabla' v(y)\cdot (\nabla'\varphi)(\varepsilon y')dy+O(\varepsilon^2)
\ee
Using the decay property of $u$ and $v$ (see \cite[Theorem 7.2]{ZhongZou.2015arXiv:1504.01005v1[math.AP]4Apr2015} or \cite[Theorem 1.1]{ZhongZou.2015arXiv:submit/1229202[math.AP]12Apr2015}), by \eqref{2014-3-1-e0} we have (the details we refer to  \cite[Lemma 3.4]{ZhongZou.2015arXiv:1504.00730[math.AP]}):
\begin{align}\lab{2015-3-10-e1}
&-2\int_{\mathbb{B}_{\frac{r_0}{\varepsilon}}^{+}}\eta\left(\phi^{-1}(\varepsilon y)\right)^2\partial_N u(y) \nabla' u(y)\cdot (\nabla'\varphi)(\varepsilon y')dy\nonumber\\
=&\frac{2}{\varepsilon}\int_{\mathbb{B}_{\frac{r_0}{\varepsilon}}^{+}}\eta\left(\phi^{-1}(\varepsilon y)\right)^2\partial_N u(y)\sum_{i=1}^{N-1}\partial_{ii}u(y)\varphi(\varepsilon y')dy+O(\varepsilon^2)
\end{align}
and
\begin{align}\lab{2015-3-10-e2}
&-2\int_{\mathbb{B}_{\frac{r_0}{\varepsilon}}^{+}}\eta\left(\phi^{-1}(\varepsilon y)\right)^2\partial_N v(y) \nabla' v(y)\cdot (\nabla'\varphi)(\varepsilon y')dy\nonumber\\
=&\frac{2}{\varepsilon}\int_{\mathbb{B}_{\frac{r_0}{\varepsilon}}^{+}}\eta\left(\phi^{-1}(\varepsilon y)\right)^2\partial_N v(y)\sum_{i=1}^{N-1}\partial_{ii}v(y)\varphi(\varepsilon y')dy+O(\varepsilon^2)
\end{align}
By equation \eqref{2015-4-26-we1} and the formula of integration by parts, we have
\begin{align}\lab{2015-3-10-e3}
&\frac{2}{\varepsilon}\int_{\mathbb{B}_{\frac{r_0}{\varepsilon}}^{+}}\eta\left(\phi^{-1}(\varepsilon y)\right)^2\partial_N u(y)\sum_{i=1}^{N-1}\partial_{ii}u(y)\varphi(\varepsilon y')dy\nonumber\\
=&\frac{2}{\varepsilon}\int_{\mathbb{B}_{\frac{r_0}{\varepsilon}}^{+}}\eta\left(\phi^{-1}(\varepsilon y)\right)^2\partial_N u(y)\left[\Delta u(y)-\partial_{NN}u(y)\right]\varphi(\varepsilon y')dy\nonumber\\
=&-\frac{2}{\varepsilon}\int_{\mathbb{B}_{\frac{r_0}{\varepsilon}}^{+}}\eta\left(\phi^{-1}(\varepsilon y)\right)^2\partial_N u(y)\left[\lambda \frac{u(y)^{2^*(s_1)-1}}{|y|^{s_1}}+\kappa \alpha \frac{u(y)^{\alpha-1}v(y)^\beta}{|y|^{s_2}}\right]\varphi(\varepsilon y')dy\nonumber\\
&-\frac{1}{\varepsilon}\int_{\mathbb{B}_{\frac{r_0}{\varepsilon}}^{+}}\eta\left(\phi^{-1}(\varepsilon y)\right)^2\partial_N \left[\partial_{N}u(y)\right]^2\varphi(\varepsilon y')dy\nonumber\\
=&-\frac{2}{\varepsilon}\frac{\lambda}{2^*(s_1)}\int_{\mathbb{B}_{\frac{r_0}{\varepsilon}}^{+}}\eta\left(\phi^{-1}(\varepsilon y)\right)^2\frac{\partial_N \left[|u(y)|^{2^*(s_1)}\right]}{|y|^{s_1}}\varphi(\varepsilon y')dy\nonumber\\
&-\frac{2}{\varepsilon}\kappa\int_{\mathbb{B}_{\frac{r_0}{\varepsilon}}^{+}}\eta\left(\phi^{-1}(\varepsilon y)\right)^2\frac{\partial_N\left(u(y)^\alpha\right)v(y)^\beta}{|y|^{s_2}}\varphi(\varepsilon y')dy\nonumber\\
&+\frac{1}{\varepsilon}\int_{\mathbb{B}_{\frac{r_0}{\varepsilon}}^{+}\cap \partial \R_+^N}\eta\left(\phi^{-1}(\varepsilon y)\right)^2\left(\partial_N u(y)\right)^2\varphi(\varepsilon y')dS_y+O(\varepsilon^2)
\end{align}
Similarly, we have
\begin{align}\lab{2015-3-10-e4}
&\frac{2}{\varepsilon}\int_{\mathbb{B}_{\frac{r_0}{\varepsilon}}^{+}}\eta\left(\phi^{-1}(\varepsilon y)\right)^2\partial_N v(y)\sum_{i=1}^{N-1}\partial_{ii}v(y)\varphi(\varepsilon y')dy\nonumber\\
=&-\frac{2}{\varepsilon}\frac{\mu}{2^*(s_1)}\int_{\mathbb{B}_{\frac{r_0}{\varepsilon}}^{+}}\eta\left(\phi^{-1}(\varepsilon y)\right)^2\frac{\partial_N \left[|v(y)|^{2^*(s_1)}\right]}{|y|^{s_1}}\varphi(\varepsilon y')dy\nonumber\\
&-\frac{2}{\varepsilon}\kappa\int_{\mathbb{B}_{\frac{r_0}{\varepsilon}}^{+}}\eta\left(\phi^{-1}(\varepsilon y)\right)^2\frac{\partial_N\left(v(y)^\beta\right)u(y)^\alpha}{|y|^{s_2}}\varphi(\varepsilon y')dy\nonumber\\
&+\frac{1}{\varepsilon}\int_{\mathbb{B}_{\frac{r_0}{\varepsilon}}^{+}\cap \partial \R_+^N}\eta\left(\phi^{-1}(\varepsilon y)\right)^2\left(\partial_N v(y)\right)^2\varphi(\varepsilon y')dS_y+O(\varepsilon^2)
\end{align}
Then, by  \eqref{2015-3-9-e5}-\eqref{2015-3-10-e4} and the formula of integration by parts, we have
\begin{align}\lab{2015-3-10-e5}
&\int_\Omega |\nabla \hat{u}_\varepsilon(x)|^2dx+\int_\Omega |\nabla \hat{v}_\varepsilon(x)|^2dx\nonumber\\
=&\int_{\R_+^N}|\nabla u(y)|^2dy+\int_{\R_+^N}|\nabla v(y)|^2dy\nonumber\\
&-\frac{2}{\varepsilon}\frac{1}{2^*(s_1)}\int_{\mathbb{B}_{\frac{r_0}{\varepsilon}}^{+}}\eta\left(\phi^{-1}(\varepsilon y)\right)^2\frac{\partial_N \left[\lambda|u(y)|^{2^*(s_1)}+\mu |v(y)|^{2^*(s_1)}\right]}{|y|^{s_1}}\varphi(\varepsilon y')dy\nonumber\\
&-\frac{2}{\varepsilon}\kappa\int_{\mathbb{B}_{\frac{r_0}{\varepsilon}}^{+}}\eta\left(\phi^{-1}(\varepsilon y)\right)^2\frac{\partial_N\left[u(y)^\alpha v(y)^\beta\right]}{|y|^{s_2}}\varphi(\varepsilon y')dy\nonumber\\
&+\frac{1}{\varepsilon}\int_{\mathbb{B}_{\frac{r_0}{\varepsilon}}^{+}\cap \partial \R_+^N}\eta\left(\phi^{-1}(\varepsilon y)\right)^2\left[\left(\partial_N u(y)\right)^2+\left(\partial_N v(y)\right)^2\right]\varphi(\varepsilon y')dS_y+O(\varepsilon^2)\nonumber\\
=&\int_{\R_+^N}|\nabla u(y)|^2dy+\int_{\R_+^N}|\nabla v(y)|^2dy\nonumber\\
&-\frac{2}{\varepsilon}\frac{1}{2^*(s_1)}\int_{\mathbb{B}_{\frac{r_0}{\varepsilon}}^{+}}\eta\left(\phi^{-1}(\varepsilon y)\right)^2\frac{\left(\lambda s_1u(y)^{2^*(s_1)}+\mu s_1v(y)^{2^*(s_1)}\right)y_N}{|y|^{s_1+2}}\varphi(\varepsilon y')dy\nonumber\\
&-\frac{2}{\varepsilon}\kappa\int_{\mathbb{B}_{\frac{r_0}{\varepsilon}}^{+}}\eta\left(\phi^{-1}(\varepsilon y)\right)^2\frac{s_2\left(u(y)^\alpha v(y)^\beta\right)y_N}{|y|^{s_2+2}}\varphi(\varepsilon y')dy\nonumber\\
&+\frac{1}{\varepsilon}\int_{\mathbb{B}_{\frac{r_0}{\varepsilon}}^{+}\cap \partial \R_+^N}\eta\left(\phi^{-1}(\varepsilon y)\right)^2\left[\left(\partial_N u(y)\right)^2+\left(\partial_N v(y)\right)^2\right]\varphi(\varepsilon y')dS_y+O(\varepsilon^2)\nonumber\\
:=&\int_{\R_+^N}|\nabla u(y)|^2dy+\int_{\R_+^N}|\nabla v(y)|^2dy+J_1+J_2+J_3+O(\varepsilon^2).
\end{align}
where
\begin{align}\lab{2015-3-10-e6}
J_1=&-\frac{2}{\varepsilon}\frac{s_1}{2^*(s_1)}\int_{\mathbb{B}_{\frac{r_0}{\varepsilon}}^{+}}\eta\left(\phi^{-1}(\varepsilon y)\right)^2\frac{\left(\lambda u(y)^{2^*(s_1)}+\mu v(y)^{2^*(s_1)}\right)y_N}{|y|^{s_1+2}}\varphi(\varepsilon y')dy\nonumber\\
=&-\frac{2s_1}{2^*(s_1)}\int_{\R_+^N}\frac{\left(\lambda u(y)^{2^*(s_1)}+\mu v(y)^{2^*(s_1)}\right)y_N}{|y|^{s_1+2}}|y'|^2dy H(0)\left(1+o(1)\right)\varepsilon+O(\varepsilon^2)\nonumber\\
:=&-\frac{2s_1}{2^*(s_1)}\tilde{K}_1H(0)\left(1+o(1)\right)\varepsilon+O(\varepsilon^2),
\end{align}
\begin{align}\lab{2015-3-10-e7}
J_2=&-\frac{2}{\varepsilon}\kappa s_2\int_{\mathbb{B}_{\frac{r_0}{\varepsilon}}^{+}}\eta\left(\phi^{-1}(\varepsilon y)\right)^2\frac{\left(u(y)^\alpha v(y)^\beta\right)y_N}{|y|^{s_2+2}}\varphi(\varepsilon y')dy\nonumber\\
=&-2\kappa s_2\int_{\R_+^N}\frac{\left(u(y)^\alpha v(y)^\beta\right)y_N}{|y|^{s_2+2}}|y'|^2dy H(0)\left(1+o(1)\right)\varepsilon+O(\varepsilon^2)\nonumber\\
:=&-2\kappa s_2\tilde{K}_2H(0)\left(1+o(1)\right)\varepsilon+O(\varepsilon^2),
\end{align}
\begin{align}\lab{2015-3-10-e8}
J_3=&\frac{1}{\varepsilon}\int_{\mathbb{B}_{\frac{r_0}{\varepsilon}}^{+}\cap \partial \R_+^N}\eta\left(\phi^{-1}(\varepsilon y)\right)^2\left[\left(\partial_N u(y)\right)^2+\left(\partial_N v(y)\right)^2\right]\varphi(\varepsilon y')dS_y\nonumber\\
=&\int_{\R^{N-1}}\left[\left((\partial_N u)(y',0)\right)^2+\left((\partial_N v)(y',0)\right)^2\right]|y'|^2dy' H(0)\left(1+o(1)\right)\varepsilon+O(\varepsilon^2)\nonumber\\
:=&\tilde{K}_3H(0)\left(1+o(1)\right)\varepsilon+O(\varepsilon^2).
\end{align}
On the other hand, by the standard estimation, we have
\be\lab{2015-3-8-e2}
\int_{\Omega}\frac{|\hat{u}_\varepsilon(x)|^{2^*(s_1)}}{|x|^{s_1}}dx=\int_{\R_+^N}\frac{|u(y)|^{2^*(s_1)}}{|y|^{s_1}}dy-\int_{\R_+^N}\frac{s_1|u(y)|^{2^*(s_1)}|y'|^2y_N}{|y|^{2+s_1}}dy H(0)\left(1+o(1)\right)\varepsilon +O(\varepsilon^2),
\ee

\be\lab{2015-3-9-e1}
\int_{\Omega}\frac{|\hat{v}_\varepsilon(x)|^{2^*(s_1)}}{|x|^{s_1}}dx=\int_{\R_+^N}\frac{|v(y)|^{2^*(s_1)}}{|y|^{s_1}}dy-\int_{\R_+^N}\frac{s_1|v(y)|^{2^*(s_1)}|y'|^2y_N}{|y|^{2+s_1}}dy H(0)\left(1+o(1)\right)\varepsilon +O(\varepsilon^2),
\ee
\be\lab{2015-3-9-e2}
\int_\Omega \frac{|\hat{u}_\varepsilon(x)|^2}{|x|^{\sigma_0}}dx=\varepsilon^{2-\sigma_0}\int_{\R_+^N}\frac{|u(y)|^2}{|y|^{\sigma_0}}dy\left(1+o(1)\right)+O(\varepsilon^N),
\ee
\be\lab{2015-3-9-e3}
\int_\Omega \frac{|\hat{v}_\varepsilon(x)|^2}{|x|^{\eta_0}}dx=\varepsilon^{2-\eta_0}\int_{\R_+^N}\frac{|v(y)|^2}{|y|^{\eta_0}}dy\left(1+o(1)\right)+O(\varepsilon^N),
\ee
\be\lab{2015-3-9-e4}
\int_{\Omega}\frac{|\hat{u}_\varepsilon(x)|^\alpha |\hat{v}_\varepsilon(x)|^\beta}{|x|^{s_2}}dx=\int_{\R_+^N}\frac{|u|^\alpha|v|^\beta}{|y|^{s_2}}dy-\int_{\R_+^N}\frac{s_2|u(y)|^\alpha |v(y)|^\beta |y'|^2y_N}{|y|^{2+s_2}}dy H(0)\left(1+o(1)\right)\varepsilon+O(\varepsilon^2).
\ee
Then, basing on the sequence of estimations above, we have
\begin{align}\lab{2015-3-10-e9}
&\Phi(t\hat{u}_\varepsilon, t\hat{v}_\varepsilon)\nonumber\\
\leq&\frac{t^2}{2}\left[\int_{\R_+^N}|\nabla u(y)|^2dy+\int_{\R_+^N}|\nabla v(y)|^2dy+\left(\tilde{K}_3-\frac{2s_1}{2^*(s_1)}\tilde{K}_1-2\kappa s_2\tilde{K}_2\right)H(0)\left(1+o(1)\right)\varepsilon\right.\nonumber\\
&\left.+\lambda^*\int_{\R_+^N}\frac{|u(y)|^2}{|y|^{\sigma_0}}dy \left(1+o(1)\right)\varepsilon^{2-\sigma_0}+\mu^*\int_{\R_+^N}\frac{|v(y)|^2}{|y|^{\eta_0}}dy \left(1+o(1)\right)\varepsilon^{2-\eta_0}+O(\varepsilon^2)\right]\nonumber\\
&-\frac{t^{2^*(s_1)}}{2^*(s_1)}\left[ \int_{\R_+^N}\frac{\lambda|u(y)|^{2^*(s_1)}+\mu|v(y)|^{2^*(s_1)}}{|y|^{s_1}}dy-s_1\tilde{K}_1H(0)\left(1+o(1)\right)\varepsilon+O(\varepsilon^2)\right]\nonumber\\
&-\kappa t^{2^*(s_2)}\left[\int_{\R_+^N}\frac{|u|^\alpha|v|^\beta}{|y|^{s_2}}dy-s_2\tilde{K}_2 H(0)\left(1+o(1)\right)\varepsilon+O(\varepsilon^2)\right]
\end{align}
and it follows easily that there exists some $T>0$ large enough and some  $\varepsilon_0>0$ small enough such that
\be\lab{2015-3-10-e10}
\Phi(T\hat{u}_\varepsilon, T\hat{v}_\varepsilon)<0\;\hbox{for all $0<\varepsilon<\varepsilon_0$}.
\ee
Especially, under the assumptions of$(\tilde{H}_1$ and  $(\tilde{H}_2)$, we can apply the similar arguments as \cite[Lemma 3.4]{ZhongZou.2015arXiv:1504.00730[math.AP]} and obtain that there exists some $0<\varepsilon_1<\varepsilon_0$ such that
\be\lab{2015-3-10-e11}
\max_{t>0} \Phi(t\hat{u}_\varepsilon, t\hat{v}_\varepsilon)<\max_{t>0}A(tu,tv)=\tilde{\Theta}\;\hbox{for all $0<\varepsilon<\varepsilon_1$}.
\ee
\hfill$\Box$

\vskip0.26in

\end{CJK*}
 \end{document}